\documentclass[11pt]{amsart}
\usepackage{amsmath,amsfonts,amssymb}
\usepackage[cp1251]{inputenc}
\usepackage[russian,english]{babel}
\def\wh{\widehat}
\def\wt{\widetilde}
\def\R{\mathbb R}
\def\C{\mathbb C}
\def\Z{\mathbb Z}

\def\A{\mathcal A}
\def\B{\mathcal B}
\def\O{\mathcal O}

\def\u{\mathbf u}

\def\x{\mathbf x}
\def\v{\mathbf v}
\def\a{\mathbf a}

\def\p{\mathbf p}
\def\r{\mathbf r}
\def\b{\mathbf b}

\def\w{\mathbf w}
\def\z{\mathbf z}

\def\VV{\mathbf V}

\def\D{\mathbf D}
\def\FF{\mathbf F}
\def\1{\mathbf 1}

\def\H{\mathfrak H}

\def\O{\mathcal O}

\def\bxi{\boldsymbol \xi}
\def\bnu{\boldsymbol \nu}

\def\1{\bold 1}

\def\div{\mathrm{div}\,}

\def\rank{\mathrm{rank}\,}

\def\eps{\varepsilon}

\def\leq{\leqslant}
\def\le{\leqslant}
\def\geq{\geqslant}
\def\ge{\geqslant}

\begin{document}

\title[Homogenization of elliptic problems]
{Homogenization of elliptic problems:\\ error estimates in dependence \\ of the spectral parameter}
\thanks{Supported by St.~Petersburg State University (grant no. 11.38.63.2012) and RFBR (grant no. 14-01-00760).}

\author{T.~A.~Suslina}

\address{Department of Physics\\ St.~Petersburg State University \\ Ul'yanovskaya 3 \\ Petrodvorets\\ St.~Petersburg\\ 198504 \\ Russia.}

\email{suslina@list.ru}

\begin{abstract}
We consider a strongly elliptic differential expression of the form $b(\D)^* g(\x/\eps) b(\D)$, $\eps >0$,
where $g(\x)$ is a matrix-valued function in $\R^d$ assumed to be bounded, positive definite and periodic with respect to some lattice;
$b(\D)=\sum_{l=1}^d b_l D_l$ is the first order differential operator with constant coefficients.
The symbol $b({\boldsymbol \xi})$ is subject to some condition ensuring strong ellipticity.
The operator given by $b(\D)^* g(\x/\eps) b(\D)$  in $L_2(\R^d;\C^n)$ is denoted by $\A_\eps$.
Let $\O \subset \R^d$ be a bounded domain of class $C^{1,1}$.
In $L_2(\O;\C^n)$, we consider the operators $\A_{D,\eps}$ and $\A_{N,\eps}$
given by $b(\D)^* g(\x/\eps) b(\D)$ with the Dirichlet or Neumann boundary conditions, respectively.
For the resolvents of the operators $\A_\eps$, $\A_{D,\eps}$, and $\A_{N,\eps}$ in a regular point
$\zeta$ we find approximations in different operator norms with error estimates depending on $\eps$ and the spectral parameter $\zeta$.
\end{abstract}


\maketitle

\section*{Introduction}

The paper concerns homogenization theory of periodic differential operators (DO's).
A broad literature is devoted to homogenization problems.
First, we mention the books [BeLP], [BaPa], [ZhKO].

\noindent\textbf{0.1. The class of operators.}
We study a wide class of matrix (of size $n\times n$) strongly elliptic operators
given by the differential expression $b(\D)^* g(\x/\eps)b(\D)$, $\eps>0$.
Here $g(\x)$ is a Hermitian matrix-valued function in $\R^d$ (of size $m\times m$), it is assumed
to be bounded, positive definite and periodic with respect to a lattice $\Gamma$. The operator $b(\D)$ is an $(m\times n)$-matrix first order DO
with constant coefficients. It is assumed that $m \geq n$; the symbol of $b(\D)$ is subject to
some condition ensuring strong ellipticity of the operator under consideration.

The selfadjoint operator in $L_2(\R^d;\C^n)$ given by the differential expression $b(\D)^* g(\x/\eps)b(\D)$ is denoted by $\A_\eps$.
We also study selfadjoint operators $\A_{D,\eps}$ and $\A_{N,\eps}$ in $L_2(\O;\C^n)$ given by
the same expression with the Dirichlet or Neumann boundary conditions, respectively.
Here $\O \subset \R^d$ is a bounded domain with the boundary of class $C^{1,1}$.

The simplest example of the operator $\A_\eps$ is the acoustics operator  $- {\rm div}\,g(\x/\eps) \nabla$;
the operator of elasticity theory can be also written in the required form. These and other examples
are considered in [BSu2] in detail.

We study homogenization problem: in operator terms, the problem is to approximate the resolvents
of the operators introduced above for small $\eps$ in different operator norms.

\smallskip\noindent\textbf{0.2. A survey of the results on operator error estimates.}
Homogenization problems for the operator $\A_\eps$ in $L_2(\R^d;\C^n)$ have been studied in a series of papers by Birman and Suslina
[BSu1--4]. In [BSu1,2] it was proved that
$$
\| ({\A}_\eps +I)^{-1} - ({\A}^0 +I)^{-1} \|_{L_2(\R^d) \to L_2(\R^d)}
\le C \eps.
\eqno(0.1)
$$
Here $\A^0 = b(\D)^* g^0 b(\D)$ is the effective operator and $g^0$ is the constant positive effective matrix
($g^0$ is defined in \S 1 below.) Next, in [BSu4] approximation of the resolvent $({\A}_\eps +I)^{-1}$
in the norm of operators acting from $L_2(\R^d;\C^n)$ to the Sobolev space $H^1(\R^d;\C^n)$ was found:
$$
\| ({\A}_\eps +I)^{-1} - ({\A}^0 +I)^{-1} - \eps K(\eps) \|_{L_2(\R^d) \to H^1(\R^d)}
\le C \eps.
\eqno(0.2)
$$
Here $K(\eps)$ is a corrector. Estimates of the form (0.1), (0.2) called the \textit{operator error estimates}
are order-sharp; the constants in estimates are controlled in terms of the problem data.
The method of [BSu1--4] is based on the scaling transformation, the Floquet-Bloch theory and the analytic perturbation theory.

A different approach to operator error estimates in homogenization problems was suggested by
Zhikov. In [Zh1,2], [ZhPas], the acoustics operator and the operator of elasticity theory were studied.
Estimates of the form (0.1), (0.2) for the corresponding problems in $\R^d$ were obtained.
The method is based on analysis of the first order approximation to the solution and introduction of an additional parameter.
Besides the problems in $\R^d$, in [Zh1,2], [ZhPas] homogenization problems in a bounded domain $\O \subset \R^d$ with the
Dirichlet or Neumann boundary conditions were studied;
the analogs of estimates (0.1), (0.2) with the error terms $O(\eps^{1/2})$ were obtained.
Error estimates become worse due to the boundary influence.
(In the case of the Dirichlet problem for the acoustics operator, the $(L_2 \to L_2)$-estimate was improved in [ZhPas];
but still it was not order-sharp.)

Similar results for the operator $- \div g(\x/\eps) \nabla$ in a bounded domain with the Dirichlet or Neumann
boundary conditions were obtained in the papers [Gr1,2] by Griso by the unfolding method.
In [Gr2], the analog of sharp order estimate (0.1) (for the same operator) was proved for the first time.

For the matrix operators $\A_{D,\eps}$, $\A_{N,\eps}$ that we consider operator error estimates
were obtained in the recent papers [PSu1,2], [Su1--3].
In [PSu1,2] the Dirichlet problem was studied; it was proved that
$$
\| \A_{D,\eps}^{-1} - (\A_D^0)^{-1} - \eps K_D(\eps) \|_{L_2(\O) \to H^1(\O)} \leq C \eps^{1/2}.
$$
Here $\A_D^0$ is the operator given by $b(\D)^* g^0 b(\D)$ with the Dirichlet condition, and $K_D(\eps)$ is the corresponding corrector.
In [Su1,2], the following sharp order estimate was obtained:
$$
\| \A_{D,\eps}^{-1} - (\A_D^0)^{-1}  \|_{L_2(\O) \to L_2(\O)} \leq C \eps.
\eqno(0.3)
$$
Similar results for the Neumann problem were obtained in [Su3].
(Note that in [Su3] the resolvent $(\A_{N,\eps} - \zeta I)^{-1}$ in an arbitrary regular point
$\zeta \in \C \setminus \R_+$ was considered, but the optimal dependence of the constants
in estimates on the parameter $\zeta$ was not searched out.)
The method of [PSu1,2], \hbox{[Su1--3]} is based on using the results for the problem in $\R^d$,
introduction of the boundary layer correction term and estimates for the norms of this term in $H^1(\O)$ and $L_2(\O)$.
Some technical tools are borrowed from [ZhPas].

Independently, by a different method estimate of the form (0.3) was obtained in [KeLiS]
for uniformly elliptic systems with the Dirichlet or Neumann boundary conditions
under some regularity assumptions on coefficients.

\smallskip\noindent\textbf{0.3. Main results.}
In the present paper, the operators $\A_\eps$, $\A_{D,\eps}$, and $\A_{N,\eps}$ are studied.
\textit{Our goal} is to find approximations of the resolvent in a regular point $\zeta$
in dependence of $\eps$ and the spectral parameter $\zeta$.
Estimates for small $\eps$ and large $|\zeta|$ are of main interest.

Now we describe main results. For the operator $\A_\eps$ in
$L_2(\R^d;\C^n)$ we prove the following estimates for $\zeta = |\zeta|e^{i\varphi}\in \C \setminus \R_+$ and $\eps>0$:
$$
\| (\A_\eps - \zeta I)^{-1} - (\A^0 - \zeta I)^{-1}\|_{L_2(\R^d)\to L_2(\R^d)} \leq
C(\varphi) |\zeta|^{-1/2} \eps,
\eqno(0.4)
$$
$$
\| (\A_\eps - \zeta I)^{-1} - (\A^0 - \zeta I)^{-1} - \eps K(\eps;\zeta) \|_{L_2(\R^d)\to H^1(\R^d)} \leq
C(\varphi) (1+|\zeta|^{-1/2}) \eps.
\eqno(0.5)
$$
For the operators $\A_{D,\eps}$ and $\A_{N,\eps}$ the following estimates are obtained for
$\zeta \in \C \setminus \R_+$ and $|\zeta|\geq 1$:
$$
\| (\A_{\dag,\eps} - \zeta I)^{-1} - (\A^0_\dag - \zeta I)^{-1}\|_{L_2(\O)\to L_2(\O)} \leq
C_1(\varphi) (|\zeta|^{-1/2} \eps + \eps^2),
\eqno(0.6)
$$
$$
\begin{aligned}
&\| (\A_{\dag,\eps} - \zeta I)^{-1} - (\A^0_\dag - \zeta I)^{-1} - \eps K_\dag(\eps;\zeta)
\|_{L_2(\O)\to H^1(\O)}
\cr
&\leq
C_2(\varphi) |\zeta|^{-1/4} \eps^{1/2} + C_3(\varphi) \eps,
\end{aligned}
\eqno(0.7)
$$
where $0< \eps \leq \eps_1$ ($\eps_1$ is a sufficiently small number depending on the domain $\O$ and the lattice $\Gamma$).
Here $\dag =D,N$. The dependence of constants in estimates (0.4)--(0.7) on the angle $\varphi$ is traced.
Estimates (0.4)--(0.7) are two-parametric with respect to $\eps$ and $|\zeta|$;
they are uniform in any sector $\varphi \in [\varphi_0, 2 \pi - \varphi_0]$
with arbitrarily small $\varphi_0>0$.

In the general case, the correctors in (0.5), (0.7) contain a smoothing operator.
We distinguish an additional condition under which the standard corrector can be used.

Besides approximation of the resolvent, we find approximation of the operators
$g^\eps b(\D) (\A_\eps - \zeta I)^{-1}$ (corresponding to the flux) in the $(L_2 \to L_2)$-operator norm.
Also, for a strictly interior subdomain $\O'$ of the domain $\O$ we find approximation of the resolvent
of $\A_{D,\eps}$ and $\A_{N,\eps}$ in the $(L_2(\O) \to H^1(\O'))$-norm with sharp order error estimate (with respect to $\eps$).

For completeness, we find a different type approximation for the resolvent of $\A_{D,\eps}$ and $\A_{N,\eps}$,
which is valid in a wider range of the parameter $\zeta$ and may be preferable for bounded values of $|\zeta|$.
Let us explain these results for the Dirichlet problem. Let $c_* >0$ be a common lower bound of the operators
$\A_{D,\eps}$ and $\A_D^0$. Let $\zeta \in \C \setminus [c_*,\infty)$.
We put $\zeta - c_* = |\zeta - c_*| e^{i \psi}$. Then we have
$$
\| (\A_{D,\eps} - \zeta I)^{-1} - (\A^0_D - \zeta I)^{-1}\|_{L_2(\O)\to L_2(\O)} \leq
\hat{C}(\zeta) \eps,
\eqno(0.8)
$$
$$
\| (\A_{D,\eps} - \zeta I)^{-1} - (\A^0_D - \zeta I)^{-1} - \eps K_D(\eps;\zeta)
\|_{L_2(\O)\to H^1(\O)} \leq
\hat{C}(\zeta)  \eps^{1/2}
\eqno(0.9)
$$
for $0< \eps \leq \eps_1$, where $\hat{C}(\zeta) = C(\psi) |\zeta - c_*|^{-2}$ for $|\zeta - c_*| <1$ and
$\hat{C}(\zeta) = C(\psi)$ for $|\zeta - c_*| \geq 1$. The dependence of constants in estimates (0.8), (0.9)
on the angle $\psi$ is traced.
Estimates (0.8), (0.9) are uniform with respect to $\eps$ and $|\zeta - c_*|$
in any sector $\psi \in [\psi_0, 2\pi - \psi_0]$ with arbitrarily small $\psi_0$.

The author considers estimates (0.6), (0.7) as the main achievements of the paper.

\smallskip\noindent\textbf{0.4. The method.} For the problem in $\R^d$ it is not difficult to
deduce estimates (0.4), (0.5) from the known estimates in the point $\zeta =-1$ (see (0.1), (0.2))
with the help of appropriate resolvent identities and a scaling transformation (see \S 2).
However, this method is not suitable for the problems in a bounded domain.
In order to prove estimates (0.6) and (0.7), we have to repeat the whole scheme from [PSu2], [Su2,3]
tracing the dependence of estimates on the spectral parameter $\zeta$ carefully.
The method is based on using an extension operator $P_\O: H^s(\O) \to H^s(\R^d)$, the study of the associated problem in
$\R^d$, introduction of the boundary layer correction term, and a careful analysis of this correction term.
Using the Steklov smoothing operator (borrowed from [ZhPas]) and estimates in the $\eps$-neighborhood of the boundary play an important technical role.
First we prove estimate (0.7), and next we prove (0.6) using
the already proved inequality (0.7) and the duality arguments.

We have achieved some technical simplifications in the scheme compared with [PSu2], [Su2,3]:
the presentation is consistently given in terms of integral identities,
we avoid consideration of the Sobolev spaces with negative indices, and
in the Neumann problem we avoid consideration of the conormal derivative
and introduce the correction term in a different way compared with [Su3].

Estimates (0.8), (0.9) are deduced from the corresponding estimates with $\zeta =-1$ and
appropriate resolvent identities.

\smallskip\noindent\textbf{0.5. Application of the results.} The present investigation was stimulated by
the study of homogenization of initial boundary value problems for a parabolic equation
$$
\frac{\partial \u_\eps(\x,t)}{\partial t} = - b(\D)^* g(\x/\eps) b(\D) \u_\eps(\x,t),\ \ \x \in \O,
\ \ t\geq 0,
$$
with the initial condition $\u_\eps(\x,0)= {\boldsymbol \phi}(\x)$, ${\boldsymbol \phi} \in L_2(\O;\C^n)$,
and the Dirichlet or Neumann conditions on $\partial \O$. The solution can be written in terms of
the operator exponential:
$\u_\eps = e^{-\A_{\dag,\eps}t} {\boldsymbol \phi}$, and therefore the problem is reduced to
approximation of the operator exponential $e^{-\A_{\dag,\eps}t}$ in different operator norms.
It is natural to use representation of the operator exponential as the integral
of the resolvent over an appropriate contour $\gamma$ in the complex plane:
$$
e^{-\A_{\dag,\eps}t} = - \frac{1}{2\pi i} \intop_\gamma e^{-\zeta t} (\A_{\dag,\eps} - \zeta I)^{-1}\,d\zeta.
$$
It turns out that, in order to obtain two-parametric approximations of the exponential
$e^{-\A_{\dag,\eps}t}$ of right order with respect to $\eps$ and $t$, approximations of the resolvent of the form
(0.6), (0.7) are required.
The parabolic homogenization problems will be studied in a separate paper [MSu2] (see also brief communication [MSu1]).

\smallskip\noindent\textbf{0.6. Plan of the paper.} The paper consists of three chapters.
Chapter 1 (\S 1, 2) is devoted to the problem in $\R^d$.
In \S 1, we introduce the class of operators under consideration, describe the effective operator,
and define the smoothing Steklov operator $S_\eps$. In \S 2,
we deduce estimates (0.4), (0.5) from the known estimates for $\zeta=-1$ with the help of the resolvent identities and the scaling transformation.
In the general case, the corrector in (0.5) contains the operator $S_\eps$. It is shown that under the additional condition
(that the solution $\Lambda(\x)$ of the auxiliary problem (1.7) is bounded) the operator $S_\eps$ can be removed
and the standard corrector can be used.

Chapter 2 (\S 3--8) is devoted to the Dirichlet problem. In \S 3, the statement of the problem is given, the effective operator
is described, and some auxiliary statements concerning estimates in the neighborhood of the boundary are collected.
\S 4 contains formulations of the main results for the Dirichlet problem (Theorems 4.2 and 4.3) and the first two steps
of the proof: the associated problem in $\R^d$ is considered, the boundary layer correction term
$\w_\eps$ is introduced, and the problem is reduced to estimates of $\w_\eps$ in $H^1(\O)$ and $L_2(\O)$.
In \S 5, the required estimates of the correction term are obtained and the proof of Theorems 4.2 and 4.3 is completed.
In \S 6, the case of $\Lambda \in L_\infty$ and some special cases are considered.
\S 7 is devoted to approximation of the solutions of the Dirichlet problem in a strictly interior subdomain of $\O$.
In \S 8, approximation of the resolvent $(\A_{D,\eps}-\zeta I)^{-1}$ of a different type (estimates (0.8), (0.9)) is obtained.

Chapter 3 (\S 9--14) is devoted to the Neumann problem. In \S 9, the statement of the problem is given and the effective
operator is defined. In \S 10, the main results for the Neumann problem (Theorems 10.1 and 10.2) are formulated
 and the first two steps of the proof are made: the associated problem in $\R^d$ is considered,
 the boundary layer correction term $\w_\eps$ is introduced, and the problem is reduced to
 estimates of $\w_\eps$ in $H^1(\O)$ and $L_2(\O)$.
In \S 11, the required estimates of $\w_\eps$ are obtained and the proof of Theorems 10.1 and 10.2 is completed.
In \S 12, the case where $\Lambda \in L_\infty$ and the special cases are considered; also, estimates in a strictly
interior subdomain of $\O$ are obtained.
In \S 13, we find approximation of the resolvent $(\A_{N,\eps} - \zeta I)^{-1}$ of a different type for $\zeta \in \C \setminus \R_+$.
Finally, in \S 14 we consider the operator $\B_{N,\eps}$ which is the part of $\A_{N,\eps}$ in the invariant subspace
orthogonal to $\text{Ker}\,b(\D)$. For the resolvent of $\B_{N,\eps}$ we find approximations
(similar to (0.8), (0.9)) for $\zeta \in \C \setminus [c_\flat,\infty)$,
where $c_\flat>0$ is a common lower bound of the operators $\B_{N,\eps}$ and $\B_N^0$.
Next, these results are applied to the operator $\A_{N,\eps}$:
we find approximations of the resolvent $(\A_{N,\eps} - \zeta I)^{-1}$ in a regular point $\zeta \in \C \setminus [c_\flat,\infty)$, $\zeta \ne 0$.

\noindent\textbf{0.7. Notation.} Let $\H$ and $\H_*$ be complex separable Hilbert spaces.
The symbols $(\cdot,\cdot)_{\H}$ and $\|\cdot\|_{\H}$ stand for the inner product and the norm in $\H$;
the symbol $\|\cdot\|_{\H \to \H_*}$ stands for the norm of a linear continuous operator acting from $\H$ to $\H_*$.

The symbols $\langle \cdot, \cdot \rangle$ and $|\cdot|$ stand for the inner product and the norm in $\C^n$;
$\1 = \1_n$ is the identity $(n\times n)$-matrix.
If $a$ is an $(m\times n)$-matrix, the symbol $|a|$ denotes the norm of $a$ as an operator from $\C^n$ to $\C^m$.
We denote $\x = (x_1,\dots,x_d)\in \R^d$, $iD_j = \partial_j = \partial/\partial x_j$,
$j=1,\dots,d$, $\D = -i \nabla = (D_1,\dots,D_d)$.
The $L_p$-classes of $\C^n$-valued functions in a domain ${\mathcal O} \subset \R^d$ are denoted by
$L_p({\mathcal O};\C^n)$, $1 \le p \leq \infty$.
The Sobolev classes of $\C^n$-valued functions in a domain ${\mathcal O} \subset \R^d$
are denoted by $H^s({\mathcal O};\C^n)$.
Next, $H^1_0(\O;\C^n)$ is the closure of $C_0^\infty(\O;\C^n)$ in $H^1(\O;\C^n)$.
If $n=1$, we write simply $L_p({\mathcal O})$, $H^s({\mathcal O})$, etc., but often we use such abbreviated
notation also for the spaces of vector-valued or matrix-valued functions.

We denote $\R_+ = [0,\infty)$.
Different constants in estimates are denoted by $c$, $C$, ${\mathcal C}$, ${\mathfrak C}$
(possibly, with indices and marks).

A brief communication on the results of the present paper is published in [Su4].

\section*{Chapter 1. Homogenization problem in $\mathbb{R}^d$}

\section*{\S 1. Periodic elliptic operators in $L_2(\mathbb{R}^d;\mathbb{\C}^n)$}

In this section, we describe the class of matrix elliptic operators in $\R^d$ under consideration
and define the effective matrix and the effective operator.

\smallskip\noindent\textbf{1.1. Lattices in $\R^d$.}
Let $\Gamma \subset \R^d$ be a lattice generated by the basis $\a_1,\dots,\a_d \in \R^d$:
$$
\Gamma = \{ \a \in \R^d:\ \a = \sum_{j=1}^d \nu_j \a_j,\ \ \nu_j \in \Z\},
$$
and let $\Omega$ be the (elementary) cell of the lattice $\Gamma$, i.~e.,
$$
\Omega := \{\x \in \R^d:\ \x = \sum_{j=1}^d \tau_j \a_j,\ \ -\frac{1}{2} < \tau_j < \frac{1}{2} \}.
$$
We use the notation $|\Omega| = \text{meas}\, \Omega$.

The basis $\b_1,\dots, \b_d$ in $\R^d$ dual to $\a_1,\dots, \a_d$
is defined by the relations $\langle \b_i, \a_j\rangle  = 2 \pi \delta_{ij}$.
This basis generates a lattice $\wt{\Gamma}$ dual to the lattice $\Gamma$.
Below we use the notation
$$
r_0 = \frac{1}{2} \min_{0 \ne \b \in \wt{\Gamma}} |\b|,\quad r_1 = \frac{1}{2} \text{diam}\,\Omega.
$$

By $\wt{H}^1(\Omega)$ we denote the subspace of all functions in
$H^1(\Omega)$ whose $\Gamma$-periodic extension to $\R^d$ belongs to $H^1_{\text{loc}}(\R^d)$.
If $\varphi(\x)$ is a $\Gamma$-periodic function in $\R^d$, we denote
$$
\varphi^\eps(\x) := \varphi(\eps^{-1}\x),\quad \eps >0.
$$

\smallskip\noindent\textbf{1.2. The class of operators.}
In $L_2(\R^d;\C^n)$, consider a second order DO $\A_\eps$ formally given by the differential expression
$$
\mathcal{A}_\eps = b(\D)^* g^\eps(\x) b (\D),\ \ \eps >0.
\eqno(1.1)
$$
Here $g(\x)$ is a measurable Hermitian $(m \times m)$-matrix-valued function (in general, with complex entries).
It is assumed that $g(\x)$ is periodic with respect to the lattice $\Gamma$, bounded and uniformly
positive definite. Next, $b(\D)$ is an $(m\times n)$-matrix first order DO with constant coefficients given by
$$
b(\D)=\sum_{l=1}^d b_l D_l.
\eqno(1.2)
$$
Here $b_l$ are constant matrices (in general, with complex entries).
The symbol $b(\boldsymbol{\xi}) = \sum_{l=1}^d b_l \xi_l$, $\boldsymbol{\xi} \in \mathbb{R}^d$, corresponds to
the operator $b(\D)$. It is assumed that $m \ge n$ and
$$
\textrm{rank} \, b (\boldsymbol{\xi}) = n, \ 0 \ne \boldsymbol{\xi} \in \R^d.
\eqno(1.3)
$$
This condition is equivalent to the inequalities
$$
\alpha_0 \mathbf{1}_n \leq b(\boldsymbol{\theta})^* b(\boldsymbol{\theta}) \leq \alpha_1 \mathbf{1}_n,
\quad \boldsymbol{\theta} \in \mathbb{S}^{d-1}, \quad 0 < \alpha_0 \leq \alpha_1 < \infty,
\eqno(1.4)
$$
with some positive constants $\alpha_0$ and $\alpha_1$. From (1.4) it follows that
$$
|b_l| \le \alpha_1^{1/2},\quad l=1,\dots,d.
\eqno(1.5)
$$
The precise definition of the operator $\mathcal{A}_\eps$ is given in terms of the quadratic form
$$
a_\eps[\u,\u] = \int_{\mathbb{R}^d} \left\langle g^\eps (\x) b (\D) \u, b(\D) \u \right\rangle \,
d \x, \quad \u \in H^1 (\mathbb{R}^{d}; \mathbb{C}^{n}).
$$
Under the above assumptions this form is closed in $L_2(\R^d;\C^n)$ and non-negative.
Using the Fourier transformation and condition (1.4), it is easy to check that
$$
\begin{aligned}
c_0 \int_{\R^d} |\D \u|^2\, d\x \le a_\eps[\u, \u] \le c_1 \int_{\R^d} |\D \u|^2\, d\x,
\ \ \u \in H^1(\R^d;\C^n),
\cr
c_0 = \alpha_0 \|g^{-1}\|^{-1}_{L_\infty},\ \  c_1 = \alpha_1 \|g\|_{L_\infty}.
\end{aligned}
\eqno(1.6)
$$

The simplest example of the operator (1.1) is the scalar elliptic operator
$
\A_\eps = -\div g^\eps(\x) \nabla = \D^* g^\eps(\x)\D.
$
In this case $n=1$, $m=d$, $b(\D)=\D$. Obviously, condition (1.4) is valid with $\alpha_0=\alpha_1=1$.
The operator of elasticity theory also can be written in the form (1.1)
with $n=d$ and $m=d(d+1)/2$. These and other examples are considered in [BSu2] in detail.

\smallskip\noindent\textbf{1.3. The effective operator.}
In order to formulate the results, we need to introduce the effective operator $\A^0$.
Let an $(n\times m)$-matrix-valued function $\Lambda \in \widetilde{H}^1(\Omega)$ be the (weak) $\Gamma$-periodic solution of the problem
$$
b(\D)^* g (\x)\left( b(\D) \Lambda(\x) + \mathbf{1}_m \right) = 0,
\quad \int_{\Omega} \Lambda(\x) \, d \x = 0.
\eqno(1.7)
$$
The so-called \emph{effective matrix} $g^0$ of size $m \times m$ is defined as follows:
$$
g^0 = |\Omega|^{-1}
\int_{\Omega} \wt{g}(\x)\, d \x,
\eqno(1.8)
$$
where
$$
\wt{g}(\x):= g (\x) \left( b (\D) \Lambda (\x) + \mathbf{1}_m \right).
\eqno(1.9)
$$
It turns out that the matrix $g^0$ is positive definite.
The \textit{effective operator} $\A^0$ for the operator (1.1) is given by the differential expression
$$
\A^0 = b (\D)^* g^0 b (\D)
$$
on the domain $H^2(\R^d;\C^n)$.

We need the following estimates for $\Lambda(\x)$ proved in [BSu3, (6.28) and subsection 7.3]:
$$
\|\D \Lambda\|_{L_2(\Omega)} \le |\Omega|^{1/2} m^{1/2} \alpha_0^{-1/2} \|g\|_{L_\infty}^{1/2}
\|g^{-1}\|_{L_\infty}^{1/2},
\eqno(1.10)
$$
$$
\| \Lambda\|_{L_2(\Omega)} \le |\Omega|^{1/2} m^{1/2} (2r_0)^{-1} \alpha_0^{-1/2} \|g\|_{L_\infty}^{1/2}
\|g^{-1}\|_{L_\infty}^{1/2}.
\eqno(1.11)
$$

\smallskip\noindent\textbf{1.4. Properties of the effective matrix}.
The following properties of the effective matrix were checked in [BSu2, Chapter 3, Theorem 1.5].

\smallskip\noindent\textbf{Proposition 1.1.} \textit{The effective matrix satisfies the estimates}
$$
\underline{g} \le g^0 \le \overline{g}.
\eqno(1.12)
$$
\textit{Here}
$$
\overline{g}= |\Omega|^{-1} \int_\Omega g(\x)\,d\x, \quad
\underline{g}= \left(|\Omega|^{-1} \int_\Omega g(\x)^{-1}\,d\x\right)^{-1}.
$$
\textit{If} $m = n$, \textit{the effective matrix} $g^0$ \textit{coincides with} $\underline{g}$.

\smallskip
In homogenization theory for specific DO's, estimates (1.12) are known as the Voight-Reuss bracketing.
We distinguish the cases where one of the inequalities in (1.12) becomes an identity.
The following statements were obtained in [BSu2, Chapter 3, Propositions 1.6 and 1.7].

\smallskip\noindent\textbf{Proposition 1.2.} \textit{The identity} $g^0 = \overline{g}$
\textit{is equivalent to the relations}
$$
b(\D)^* {\mathbf g}_k(\x) = 0, \ \ k = 1,\dots,m,
\eqno(1.13)
$$
\textit{where} ${\mathbf g}_k(\x)$, $k = 1,\dots, m$, \textit{are the columns of the matrix} $g(\x)$.

\smallskip\noindent\textbf{Proposition 1.3.} \textit{The identity} $g^0 = \underline{g}$
\textit{is equivalent to the relations}
$$
 {\mathbf l}_k(\x) = {\mathbf l}_k^0 + b(\D) \w_k, \ \ {\mathbf l}_k^0 \in \C^m,
 \ \ \w_k \in \wt{H}^1(\Omega;\C^n),\ \ k = 1,\dots,m,
\eqno(1.14)
$$
 \textit{where} ${\mathbf l}_k(\x)$, $k= 1,\dots, m$,  \textit{are the columns of the matrix} $g(\x)^{-1}$.

\smallskip
Obviously, (1.12) implies the following estimates for the norms of the matrices $g^0$ and $(g^0)^{-1}$:
$$
| g^0| \le \|g\|_{L_\infty},\ \ | (g^0)^{-1}| \le \|g^{-1}\|_{L_\infty}.
\eqno(1.15)
$$
Note that by (1.4) and (1.15), the symbol of the effective operator $\A^0$ satisfies the following inequality:
$$
c_0 |\bxi|^2 \1_n \leq b(\bxi)^* g^0 b(\bxi) \leq c_1 |\bxi|^2 \1_n,\quad \bxi \in \R^d,
\eqno(1.16)
$$
where $c_0$ and $c_1$ are as in (1.6).

\smallskip\noindent\textbf{1.5. The Steklov smoothing operator.}
Let $S_\eps$ be the operator in $L_2(\R^d;\C^m)$ defined by
$$
(S_\eps \u)(\x) = |\Omega|^{-1} \int_\Omega \u(\x - \eps \z)\, d\z
\eqno(1.17)
$$
and called the \textit{Steklov smoothing operator}.
Note that
$$
\|S_\eps\|_{L_2(\R^d)\to L_2(\R^d)} \leq 1.
\eqno(1.18)
$$
Obviously, $D^\alpha S_\eps \u = S_\eps D^\alpha \u$ for
$\u \in H^s(\R^d;\C^m)$ and any multiindex $\alpha$ such that $|\alpha| \leq s$.

We need some properties of the operator (1.17); see [ZhPas, Lemmas 1.1 and 1.2]
or [PSu2, Propositions 3.1 and 3.2].

\smallskip\noindent\textbf{Proposition 1.4.}
\textit{For any} $\u \in H^1(\R^d;\C^m)$ \textit{we have}
$$
\| S_\eps \u - \u\|_{L_2(\R^d)} \leq \eps r_1 \| \D \u \|_{L_2(\R^d)}, \ \ \eps>0.
$$

\smallskip\noindent\textbf{Proposition 1.5.} \textit{Let $f(\x)$
be a $\Gamma$-periodic function in} $\R^d$ \textit{such that} $f \in L_2(\Omega)$.
\textit{Let $[f^\eps]$ be the operator of multiplication by the function $f^\eps(\x)$. Then the operator} $[f^\eps]S_\eps$
\textit{is continuous in} $L_2(\R^d;\C^m)$, \textit{and}
$$
\| [f^\eps]S_\eps \|_{L_2(\R^d)\to L_2(\R^d)}
\leq |\Omega|^{-1/2} \| f \|_{L_2(\Omega)},\ \ \eps>0.
$$

\smallskip
From Proposition 1.5 and estimates (1.10), (1.11) it follows that
$$
\| [\Lambda^\eps]S_\eps \|_{L_2(\R^d)\to L_2(\R^d)}
\leq m^{1/2} (2r_0)^{-1} \alpha_0^{-1/2} \|g\|_{L_\infty}^{1/2}
\|g^{-1}\|_{L_\infty}^{1/2}=:M_1,\ \ \eps>0,
\eqno(1.19)
$$
$$
\| [(\D\Lambda)^\eps]S_\eps \|_{L_2(\R^d)\to L_2(\R^d)}
\leq m^{1/2} \alpha_0^{-1/2} \|g\|_{L_\infty}^{1/2}
\|g^{-1}\|_{L_\infty}^{1/2}=:M_2,\ \ \eps>0.
\eqno(1.20)
$$

\section*{\S 2. The results for homogenization problem in $\mathbb{R}^d$}

In this section, we obtain the results on homogenization of the operator $\A_\eps$ in $L_2(\R^d;\C^n)$.
Precisely, we deduce theorems about
approximation of the resolvent $(\A_\eps - \zeta I)^{-1}$ for any $\zeta \in \C \setminus \R_+$
from the known results about approximation of the resolvent $(\A_\eps + I)^{-1}$ obtained in [BSu2], [BSu4], and [PSu2].

\smallskip\noindent\textbf{2.1. Approximation of the resolvent in the operator norm in $L_2(\R^d;\C^n)$.}
A point $\zeta \in \C \setminus \R_+$ is regular both for $\A_\varepsilon$ and $\A^0$.
We put $\zeta = |\zeta|e^{i\varphi}$, $\varphi \in (0,2\pi)$, and denote
$$
c(\varphi) =
\begin{cases}
|\sin \varphi|^{-1}, & \varphi \in (0,\pi/2) \cup (3\pi/2, 2\pi), \\
1, & \varphi \in [\pi/2, 3\pi/2].
\end{cases}
\eqno(2.1)
$$
Consider the following elliptic equation in $\R^d$:
$$
\A_\varepsilon \u_{\varepsilon} - \zeta \u_{\varepsilon} = \FF,
\eqno(2.2)
$$
where $\FF \in L_2(\R^d;\C^n)$.

\smallskip\noindent\textbf{Lemma 2.1.} \textit{The solution} $\u_\eps$
\textit{of equation} (2.2) \textit{satisfies the following estimates for $\eps>0$}:
$$
\| \u_\eps \|_{L_2(\R^d)} \leq c(\varphi) |\zeta|^{-1} \|\FF\|_{L_2(\R^d)},
\eqno(2.3)
$$
$$
\| \D \u_\eps \|_{L_2(\R^d)} \leq {\mathcal C}_0 c(\varphi) |\zeta|^{-1/2} \|\FF\|_{L_2(\R^d)},
\eqno(2.4)
$$
\textit{where ${\mathcal C}_0 = \sqrt{2} c_0^{-1/2}$. In operator terms},
$$
\| (\A_\eps - \zeta I)^{-1} \|_{L_2(\R^d) \to L_2(\R^d)}
\leq c(\varphi) |\zeta|^{-1},
\eqno(2.5)
$$
$$
\| \D (\A_\eps - \zeta I)^{-1} \|_{L_2(\R^d) \to L_2(\R^d)}
\leq {\mathcal C}_0 c(\varphi) |\zeta|^{-1/2}.
$$

\smallskip
\noindent\textbf{Proof.}
 The spectrum of $\A_{\eps}$ is contained in $\R_+$.
The norm of the resolvent $(\A_{\eps} - \zeta I)^{-1}$ does not exceed
the inverse distance from the point $\zeta$ to $\R_+$.
Since this distance is equal to $|\zeta|$ if ${\rm Re}\,\zeta \leq 0$ and is equal to $|\zeta| |\sin \varphi|$
if ${\rm Re}\,\zeta > 0$, we arrive at (2.5).

In order to check (2.4), we write down the integral identity for the solution $\u_\eps$ of equation (2.2):
$$
(g^\eps b(\D)\u_\eps, b(\D) {\boldsymbol \eta})_{L_2(\R^d)} - \zeta (\u_\eps, {\boldsymbol \eta})_{L_2(\R^d)}=
(\FF, {\boldsymbol \eta})_{L_2(\R^d)},\ \  {\boldsymbol \eta}\in H^1(\R^d;\C^n).
\eqno(2.6)
$$
Putting ${\boldsymbol \eta}=\u_\eps$ in (2.6) and taking the real part, we obtain
$$
a_\eps[\u_\eps,  \u_\eps] = {\rm Re}\, \zeta \|\u_\eps\|^2_{L_2(\R^d)} +
{\rm Re}\, (\FF, \u_\eps)_{L_2(\R^d)}.
$$
Together with (2.3) this implies
$$
a_\eps[\u_\eps,  \u_\eps] \leq 2 c(\varphi)^2 |\zeta|^{-1} \|\FF\|^2_{L_2(\R^d)}.
$$
By (1.6), this yields (2.4). $\bullet$

It is known that as $\varepsilon \to 0$, the solution $\u_{\varepsilon}$ converges in
$L_2 (\R^d; \C^n)$ to the solution of the "`homogenized"' equation
$$
{\A}^0 \u_0 - \zeta \u_0 = \FF.
\eqno(2.7)
$$

\smallskip\noindent\textbf{Theorem 2.2.} \textit{Let $\zeta= |\zeta|e^{i\varphi}\in \C \setminus \R_+$. Let
$\u_\eps$ be the solution of}
(2.2), \textit{and let} $\u_0$ \textit{be the solution of} (2.7). \textit{Then}
$$
 \| \u_{\varepsilon} - \u_0 \|_{L_2 (\mathbb{R}^{d})}
 \leq C_1 c(\varphi)^2|\zeta|^{-1/2} \varepsilon  \| \FF \|_{L_2 (\mathbb{R}^{d})},
 \quad \varepsilon >0,
 $$
\textit{or, in operator terms},
$$
\| (\A_{\varepsilon} - \zeta I)^{-1} - (\A^0 - \zeta I)^{-1}
\|_{L_2 (\mathbb{R}^{d}) \rightarrow L_2 (\mathbb{R}^{d})}
\leq \ C_1 c(\varphi)^2|\zeta|^{-1/2} \varepsilon, \quad \varepsilon >0.
\eqno(2.8)
$$
\textit{Here $c(\varphi)$ is defined by} (2.1),
\textit{the constant ${C}_1$ depends only on the norms}
$\| g \|_{L_\infty}$, $\| g^{-1} \|_{L_\infty}$, \textit{the constants} $\alpha_0$, $\alpha_1$ \textit{from}
(1.4), \textit{and the parameters of the lattice} $\Gamma$.

\smallskip
\noindent\textbf{Proof.}
In [BSu2, Chapter 4, Theorem 2.1], estimate (2.8) was proved in the case where
$\zeta =-1$ and $0< \eps \le 1$.

Note that for $\zeta=-1$ and $\eps >1$ the left-hand side of (2.8) does not exceed 2, and then also does not exceed 2$\eps$.
So, we start with the estimate
$$
\| (\A_{\varepsilon} + I)^{-1} - (\A^0 + I)^{-1}
\|_{L_2 (\mathbb{R}^{d}) \rightarrow L_2 (\mathbb{R}^{d})}
\leq \ \check{C}_1 \varepsilon, \quad \varepsilon >0,
\eqno(2.9)
$$
with the constant $\check{C}_1$ depending only on
$\| g \|_{L_\infty}$, $\| g^{-1} \|_{L_\infty}$, $\alpha_0$, $\alpha_1$, and the parameters of the lattice $\Gamma$.

Now we carry over this estimate to the case $\zeta = \widehat{\zeta}=e^{i\varphi}$ using the identity
$$
\begin{aligned}
&(\A_{\varepsilon} -\zeta I)^{-1} - (\A^0 -\zeta I)^{-1}
\cr
&= (\A_\eps +I)(\A_\eps - \zeta I)^{-1}\left((\A_{\varepsilon} + I)^{-1} - (\A^0 + I)^{-1}\right)
(\A^0 +I)(\A^0 - \zeta I)^{-1}.
\end{aligned}
\eqno(2.10)
$$
We have
$$
\|(\A_\eps +I)(\A_\eps - \widehat{\zeta} I)^{-1}\|_{L_2 (\mathbb{R}^{d}) \rightarrow L_2 (\mathbb{R}^{d})}
\leq \sup_{x \ge 0} (x+1)|x-\widehat{\zeta}|^{-1} \leq 2 c(\varphi).
\eqno(2.11)
$$
The norm of the operator $(\A^0 +I)(\A^0 - \widehat{\zeta} I)^{-1}$ satisfies a similar estimate
$$
\|(\A^0 +I)(\A^0 - \widehat{\zeta} I)^{-1}\|_{L_2 (\mathbb{R}^{d}) \rightarrow L_2 (\mathbb{R}^{d})} \leq 2 c(\varphi).
\eqno(2.12)
$$
From (2.9)--(2.12) it follows that
$$
\| (\A_{\varepsilon} -\wh{\zeta} I)^{-1} - (\A^0 - \wh{\zeta} I)^{-1}
\|_{L_2 (\mathbb{R}^{d}) \rightarrow L_2 (\mathbb{R}^{d})}
\leq \ 4 \check{C}_1 c(\varphi)^2 \varepsilon, \quad \varepsilon >0.
\eqno(2.13)
$$

Next, by the scaling transformation, (2.13) is equivalent to
$$
\| (\A - \eps^2 \wh{\zeta} I)^{-1} - (\A^0 - \eps^2 \wh{\zeta} I)^{-1}
\|_{L_2 (\mathbb{R}^{d}) \rightarrow L_2 (\mathbb{R}^{d})}
\leq \ 4 \check{C}_1 c(\varphi)^2 \varepsilon^{-1}, \quad \eps >0.
\eqno(2.14)
$$
Here $\A = b(\D)^* g(\x)b(\D)$. Replacing $\eps$ by $\eps |\zeta|^{1/2}$ in (2.14)
and applying the inverse transformation, we obtain (2.8) with $C_1= 4 \check{C}_1$.
$\ \bullet$

\smallskip\noindent\textbf{2.2. Approximation of the resolvent in the $(L_2 \to H^1)$-norm.}
In order to approximate the solution $\u_\eps$ in $H^1(\R^d;\C^n)$
we have to take the first order corrector into account. We put
$$
K(\eps;\zeta) = [\Lambda^{\varepsilon}] S_{\varepsilon} b (\D) (\A^0 - \zeta I)^{-1}.
\eqno(2.15)
$$
Here $[\Lambda^{\varepsilon}]$ denotes the operator of multiplication by the matrix-valued function $\Lambda(\eps^{-1}\x)$,
and $S_\eps$ is the smoothing operator defined by (1.17).
The operator (2.15) is a continuous mapping of $L_2(\R^d;\C^n)$ into $H^1(\R^d;\C^n)$;
this follows from the identity
$$
K(\eps;\zeta) = K(\eps;-1) (\A^0 +I)(\A^0 - \zeta I)^{-1}
\eqno(2.16)
$$
and the following lemma.

\smallskip\noindent\textbf{Lemma 2.3.}
\textit{For $\eps>0$ and $\zeta=-1$ the operator} (2.15) \textit{is continuous from $L_2(\R^d;\C^n)$ to $H^1(\R^d;\C^n)$, and we have}
$$
\| K(\eps;-1)\|_{L_2(\R^d)\to L_2(\R^d)} \leq C_K^{(1)},
\eqno(2.17)
$$
$$
\|\eps \D K(\eps;-1)\|_{L_2(\R^d)\to L_2(\R^d)} \leq C_K^{(2)} \eps + C_K^{(3)}.
\eqno(2.18)
$$
\textit{The constants $C_K^{(j)}$, $j=1,2,3,$ depend only on $m$, $\alpha_0$, $\alpha_1$, $\|g\|_{L_\infty}$,
$\|g^{-1}\|_{L_\infty}$, and the parameters of the lattice} $\Gamma$.

\smallskip
\noindent\textbf{Proof.}
First, we estimate the $(L_2\to L_2)$-norm of $K(\eps;-1)$:
$$
\| K(\eps;-1) \|_{L_2(\R^d)\to L_2(\R^d)} \leq
\| [\Lambda^\eps] S_\eps\|_{L_2\to L_2} \|b(\D) (\A^0+I)^{-1}\|_{L_2\to L_2}.
\eqno(2.19)
$$
Using the Fourier transformation and (1.4), (1.16), we obtain:
$$
\begin{aligned}
&\|b(\D) (\A^0+I)^{-1}\|_{L_2(\R^d)\to L_2(\R^d)} \leq \sup_{{\boldsymbol\xi} \in \R^d}
|b({\boldsymbol\xi}) \left(b({\boldsymbol\xi})^* g^0 b({\boldsymbol\xi})+ \1 \right)^{-1}|
\cr
&\leq \alpha_1^{1/2}\sup_{{\boldsymbol\xi} \in \R^d} |{\boldsymbol\xi}| (c_0 |{\boldsymbol\xi}|^2+1)^{-1}
\leq \frac{1}{2}\alpha_1^{1/2} c_0^{-1/2}.
\end{aligned}
\eqno(2.20)
$$
Relations (1.19), (2.19), and (2.20) imply (2.17) with
$C^{(1)}_K = \frac{1}{2}\alpha_1^{1/2} c_0^{-1/2} M_1.$

Now, consider the operators
$$
\eps D_j K(\eps;-1) = [(D_j \Lambda)^\eps] S_\eps  b(\D) (\A^0+I)^{-1} +
\eps [\Lambda^\eps] S_\eps D_j b(\D) (\A^0+I)^{-1}.
\eqno(2.21)
$$
Similarly to (2.20), we have
$$
\|\D b(\D) (\A^0+I)^{-1}\|_{L_2(\R^d)\to L_2(\R^d)}
\leq \alpha_1^{1/2}\sup_{{\boldsymbol\xi} \in \R^d} |{\boldsymbol\xi}|^2 (c_0 |{\boldsymbol\xi}|^2+1)^{-1}
\leq \alpha_1^{1/2} c_0^{-1}.
\eqno(2.22)
$$
Relations (1.19), (1.20), and (2.20)--(2.22) imply
(2.18) with $C_K^{(2)}= (2\alpha_1)^{1/2} c_0^{-1}M_1$,
$C_K^{(3)}=  \alpha_1^{1/2} (2c_0)^{-1/2}M_2$.
$\bullet$

\smallskip

"`The first order approximation"' to the solution $\u_\eps$ is given by
$$
\v_{\varepsilon} = \u_0 + \varepsilon \Lambda^{\varepsilon} S_{\varepsilon} b (\D) \u_0=
(\A^0- \zeta I)^{-1} \FF + \eps K(\eps;\zeta) \FF.
\eqno(2.23)
$$

\smallskip\noindent\textbf{Theorem 2.4.} \textit{Suppose that the assumptions of Theorem} 2.2
\textit{are satisfied. Suppose that} $\v_\eps$ \textit{is defined by} (2.23). \textit{Then for $\eps>0$ we have}
$$
 \| \D(\u_{\varepsilon} - \v_\eps) \|_{L_2 (\mathbb{R}^{d})}
 \leq C_2 c(\varphi)^2 \varepsilon  \| \FF \|_{L_2 (\mathbb{R}^{d})},
 \eqno(2.24)
 $$
$$
 \| \u_{\varepsilon} - \v_\eps \|_{L_2 (\mathbb{R}^{d})}
 \leq C_3 c(\varphi)^2 |\zeta|^{-1/2} \varepsilon  \| \FF \|_{L_2 (\mathbb{R}^{d})},
 \eqno(2.25)
 $$
\textit{or, in operator terms},
$$
\| \D \left( (\A_{\varepsilon} - \zeta I)^{-1} - (\A^0 - \zeta I)^{-1} -
\varepsilon  K(\eps;\zeta)\right) \|_{L_2 (\mathbb{R}^{d})
\to L_2 (\mathbb{R}^{d})} \leq \ C_2 c(\varphi)^2 \varepsilon,
\eqno(2.26)
$$
$$
\| (\A_{\varepsilon} - \zeta I)^{-1} - (\A^0 - \zeta I)^{-1} -
\varepsilon K(\eps;\zeta) \|_{L_2 (\mathbb{R}^{d})
\to L_2 (\mathbb{R}^{d})} \leq \ C_3 c(\varphi)^2 |\zeta|^{-1/2}\varepsilon.
\eqno(2.27)
$$
\textit{Here $c(\varphi)$ is defined by} (2.1).
\textit{The constants ${C}_2$, $C_3$ depend only on
$m$, $d$, $\| g \|_{L_\infty}$, $\| g^{-1} \|_{L_\infty}$, $\alpha_0$, $\alpha_1$, and the parameters of the lattice} $\Gamma$.

\smallskip
Theorem 2.4 directly implies the following corollary.

\smallskip\noindent\textbf{Corollary 2.5.} \textit{Under the assumptions of Theorem} 2.4 \textit{we have}
$$
 \| \u_{\varepsilon} - \v_\eps \|_{H^1 (\mathbb{R}^{d})}
 \leq c(\varphi)^2 (C_2 + C_3 |\zeta|^{-1/2}) \varepsilon  \| \FF \|_{L_2 (\mathbb{R}^{d})},\ \ \eps >0.
 $$

\smallskip
\noindent\textbf{Proof of Theorem 2.4.}
In [BSu4, Theorem 10.6], a similar result was proved for $\zeta=-1$, but with a different
smoothing operator instead of $S_\eps$. In [PSu2, Theorem 3.3], it was shown that it is possible to pass to the smoothing operator $S_\eps$,
and inequalities (2.24), (2.25) were proved for $\zeta=-1$ and $0< \eps \le 1$.
Thus, we start with the estimates
$$
\begin{aligned}
\| \D \left((\A_{\varepsilon} + I)^{-1} - (\A^0 + I)^{-1} -
\varepsilon K(\eps;-1) \right) \|_{L_2 (\mathbb{R}^{d})
\to L_2 (\mathbb{R}^{d})} \leq \check{C}_2 \varepsilon,
\cr
 0< \eps\le 1,
\end{aligned}
\eqno(2.28)
$$
$$
\| (\A_{\varepsilon} + I)^{-1} - (\A^0 + I)^{-1} -
\varepsilon K(\eps;-1) \|_{L_2 (\mathbb{R}^{d})
\to L_2 (\mathbb{R}^{d})} \leq \ \check{C}_3 \varepsilon,\ 0< \eps\le 1.
\eqno(2.29)
$$
Here the constants $\check{C}_2$, $\check{C}_3$ depend only on
$m$, $d$, $\alpha_0$, $\alpha_1$, $\| g \|_{L_\infty}$, $\| g^{-1} \|_{L_\infty}$, and the parameters of the lattice $\Gamma$.

For $\eps >1$ estimates are trivial: it suffices to estimate each term under the norm sign in (2.28), (2.29) separately.
By (2.17), the left-hand side of (2.29) does not exceed $2+ C_K^{(1)}\varepsilon$,
and so does not exceed $(2+ C_K^{(1)})\varepsilon$ if $\eps>1$. Combining this with (2.29), we obtain
$$
\| (\A_{\varepsilon} + I)^{-1} - (\A^0 + I)^{-1} -
\varepsilon K (\varepsilon;-1) \|_{L_2 (\mathbb{R}^{d}) \to L_2 (\mathbb{R}^{d})} \leq \
 \widehat{C}_3 \varepsilon ,\ \eps >0,
\eqno(2.30)
$$
where $\widehat{C}_3 = \max\{\check{C}_3, 2+ C_K^{(1)}\}$.

Next, the lower estimate (1.6) implies that
$$
\|\D (\A_\eps +I)^{-1}\|_{L_2(\R^d) \to L_2(\R^d)} \le c_0^{-1/2}.
\eqno(2.31)
$$
By (1.15), the norm of the operator $\D(\A^0 + I)^{-1}$ satisfies a similar estimate:
$$
\|\D(\A^0 +I)^{-1}\|_{L_2(\R^d) \to L_2(\R^d)} \le c_0^{-1/2}.
\eqno(2.32)
$$
From (2.18), (2.31), and (2.32) it is seen that the left-hand side of (2.28) does not exceed
$2c_0^{-1/2} + C_K^{(2)}\eps + C_K^{(3)}$,
and so does not exceed $(2c_0^{-1/2} + C_K^{(2)} + C_K^{(3)})\varepsilon$ if $\eps>1$.
Together with (2.28) this yields
$$
\| \D \left((\A_{\varepsilon} + I)^{-1} - (\A^0 + I)^{-1} -
\varepsilon K(\eps;-1) \right) \|_{L_2 (\mathbb{R}^{d})
\to L_2 (\mathbb{R}^{d})} \leq \ \widehat{C}_2 \varepsilon,\  \eps>0,
\eqno(2.33)
$$
where $\widehat{C}_2= \max\{\check{C}_2, 2c_0^{-1/2} + C_K^{(2)} + C_K^{(3)}\}$.

Now we carry over estimates (2.30), (2.33) to the case where $\zeta = \widehat{\zeta}= e^{i\varphi}$
using the identity
$$
\begin{aligned}
&(\A_\eps - \zeta I)^{-1} - (\A^0 - \zeta I)^{-1} - \eps K(\eps;\zeta)
\cr
&=
(\A_\eps + I) (\A_\eps - \zeta I)^{-1} \left( (\A_\eps + I)^{-1} - (\A^0 + I)^{-1} - \eps K(\eps;-1) \right)
\cr
&\times (\A^0 + I) (\A^0 - \zeta I)^{-1} + \eps (\zeta +1) (\A_\eps - \zeta I)^{-1} K(\eps;\zeta).
\end{aligned}
\eqno(2.34)
$$
By (2.5),
$$
\|(\A_\eps - \widehat{\zeta} I)^{-1} \|_{L_2(\R^d) \to L_2(\R^d)}
\leq c(\varphi).
\eqno(2.35)
$$
From (2.12), (2.16), and (2.17) it follows that
$$
\|  K(\eps; \widehat{\zeta})\|_{L_2(\R^d) \to L_2(\R^d)} \leq 2 C_K^{(1)} c(\varphi).
\eqno(2.36)
$$
Combining (2.30) and (2.34)--(2.36) and taking (2.11), (2.12) into account, we obtain
$$
\| (\A_\eps - \widehat{\zeta} I)^{-1} - (\A^0 - \widehat{\zeta} I)^{-1} - \eps K(\eps; \widehat{\zeta})
\|_{L_2(\R^d) \to L_2(\R^d)}
\leq C_3 c(\varphi)^2 \eps,\ \ \eps>0,
\eqno(2.37)
$$
where $C_3= 4 \widehat{C}_3 + 4 C_K^{(1)}$.

Next, relations (2.11), (2.12), (2.34), and (2.36) imply that
$$
\begin{aligned}
&\| \A_\eps^{1/2}\left((\A_\eps - \widehat{\zeta} I)^{-1} - (\A^0 - \widehat{\zeta} I)^{-1} -
\eps K(\eps; \widehat{\zeta}) \right) \|_{L_2(\R^d) \to L_2(\R^d)}
\cr
&\leq 4 c(\varphi)^2
\| \A_\eps^{1/2}\left((\A_\eps + I)^{-1} - (\A^0 +I)^{-1} - \eps K(\eps;-1) \right) \|_{L_2(\R^d) \to L_2(\R^d)}
\cr
&+ 2 \eps |\widehat{\zeta}+1| C_K^{(1)} c(\varphi) \| \A_\eps^{1/2} (\A_\eps - \widehat{\zeta} I)^{-1}\|_{L_2(\R^d) \to L_2(\R^d)}.
\end{aligned}
\eqno(2.38)
$$
We have
$$
\| \A_\eps^{1/2} (\A_\eps - \widehat{\zeta} I)^{-1}\|_{L_2(\R^d) \to L_2(\R^d)}
\leq \sup_{x \geq 0} x^{1/2} |x - \widehat{\zeta}|^{-1} \leq c(\varphi).
\eqno(2.39)
$$
From (1.6), (2.33), (2.38), and (2.39) it follows that
$$
\| \A_\eps^{1/2} \left((\A_\eps - \widehat{\zeta} I)^{-1} - (\A^0 - \widehat{\zeta} I)^{-1} -
\eps K(\eps; \widehat{\zeta}) \right) \|_{L_2(\R^d) \to L_2(\R^d)}
\leq \widetilde{C}_2 c(\varphi)^2 \eps
\eqno(2.40)
$$
for $\eps>0$, where $\widetilde{C}_2= 4 (c_1^{1/2} \widehat{C}_2 + C_K^{(1)})$.
Combining this with the lower estimate (1.6), we obtain
$$
\| \D \left((\A_\eps - \widehat{\zeta} I)^{-1} - (\A^0 - \widehat{\zeta} I)^{-1} -
\eps K(\eps; \widehat{\zeta}) \right) \|_{L_2(\R^d) \to L_2(\R^d)}
\leq {C}_2 c(\varphi)^2 \eps
\eqno(2.41)
$$
for $\eps>0$, where $C_2 = c_0^{-1/2} \widetilde{C}_2$.

By the scaling transformation, (2.37) is equivalent to
$$
\| (\A - \widehat{\zeta} \eps^2 I)^{-1} - (\A^0 - \widehat{\zeta} \eps^2 I)^{-1} - \Lambda S_1 b(\D)(\A^0 - \widehat{\zeta} \eps^2 I)^{-1}
\|_{L_2 \to L_2}
\leq C_3 c(\varphi)^2 \eps^{-1}
\eqno(2.42)
$$
for $\eps >0$, and (2.41) is equivalent to
$$
\begin{aligned}
\|&\D \left((\A - \widehat{\zeta}\eps^2 I)^{-1} - (\A^0 - \widehat{\zeta} \eps^2 I)^{-1} -
\Lambda S_1 b(\D)(\A^0 - \widehat{\zeta} \eps^2 I)^{-1}\right) \|_{L_2 \to L_2}
\cr
&\leq {C}_2 c(\varphi)^2,\quad \eps >0.
\end{aligned}
\eqno(2.43)
$$
Replacing  $\eps$ by $\eps |\zeta|^{1/2}$ in (2.42) and (2.43) and
applying the inverse transformation, we arrive at (2.26), (2.27).
 $\bullet$

\smallskip
Now we obtain approximation of the flux
$\p_\eps := g^\eps b(\D) \u_\eps$ in $L_2(\R^d;\C^m)$.

\smallskip\noindent\textbf{Theorem 2.6.} \textit{Suppose that the assumptions of Theorem } 2.2
\textit{are satisfied. Let $\p_\eps := g^\eps b(\D) \u_\eps$. Then we have}
$$
 \| \p_{\varepsilon} - \wt{g}^\eps S_\eps b(\D)\u_0 \|_{L_2(\R^{d})}
 \leq C_4 c(\varphi)^2 \varepsilon  \| \FF \|_{L_2 (\mathbb{R}^{d})},
 \quad  \varepsilon >0,
 $$
 \textit{or, in operator terms,}
$$
 \| g^{\varepsilon}b(\D) (\A_\eps - \zeta I)^{-1} - \wt{g}^\eps S_\eps b(\D)(\A^0 - \zeta I)^{-1} \|_{L_2(\R^{d}) \to L_2(\R^{d})}
 \leq C_4 c(\varphi)^2 \varepsilon.
 \eqno(2.44)
 $$
\textit{Here $\wt{g}(\x)$ is the matrix} (1.9).
\textit{The constant} $C_4$ \textit{depends only on}
$d$, $m$, $\alpha_0$, $\alpha_1$, $\| g \|_{L_\infty}$, $\| g^{-1} \|_{L_\infty}$,
\textit{and the parameters of the lattice} $\Gamma$.

\smallskip
\noindent\textbf{Proof.}
We start with the case where $\zeta=\widehat{\zeta}=e^{i\varphi}$. From (2.40) it follows that
$$
 \| g^\eps b(\D) \left((\A_{\varepsilon} - \widehat{\zeta}I)^{-1} -
 (\A^0 - \widehat{\zeta}I)^{-1} - \eps K(\eps; \widehat{\zeta}) \right) \|_{L_2 \to L_2}
 \leq \|g\|_{L_\infty}^{1/2} \widetilde{C}_2 c(\varphi)^2 \eps
 \eqno(2.45)
 $$
for $\varepsilon >0$. Taking (1.2) into account, we have
$$
\begin{aligned}
\eps {g}^\eps b(\D) K(\eps; \widehat{\zeta})  &=
g^\eps(b(\D)\Lambda)^\eps S_\eps b(\D) (\A^0 - \widehat{\zeta}I)^{-1}
\cr
&+
\eps \sum_{l=1}^d g^\eps b_l \Lambda^\eps S_\eps b(\D) D_l (\A^0 - \widehat{\zeta}I)^{-1}.
\end{aligned}
\eqno(2.46)
$$

From (2.12) and (2.22) it follows that
$$
\| \D b(\D) (\A^0 - \widehat{\zeta}I)^{-1} \|_{L_2(\R^d) \to L_2(\R^d)} \leq
2 \alpha_1^{1/2} c_0^{-1} c(\varphi).
\eqno(2.47)
$$
The second term in the right-hand side of (2.46) is estimated with the help of (1.5), (1.19), and (2.47):
$$
\begin{aligned}
&\eps \left\| \sum_{l=1}^d g^\eps b_l \Lambda^\eps S_\eps b(\D) D_l (\A^0 - \widehat{\zeta}I)^{-1} \right\|_{L_2(\R^d)\to L_2(\R^d)}
\cr
&\le 2 \eps \|g\|_{L_\infty} \alpha_1 c_0^{-1} M_1 d^{1/2} c(\varphi).
\end{aligned}
\eqno(2.48)
$$
Next, by Proposition 1.4 and (2.47), we have
$$
\| g^\eps (I - S_\eps) b(\D) (\A^0 - \widehat{\zeta}I)^{-1} \|_{L_2(\R^d)\to L_2(\R^d)}
\le 2 \eps \|g\|_{L_\infty} r_1 \alpha_1^{1/2} c_0^{-1} c(\varphi).
\eqno(2.49)
$$
As a result, relations (2.45), (2.46), (2.48), and (2.49) together with (1.9) imply that
$$
\| g^\eps b(\D) (\A_\eps - \widehat{\zeta}I)^{-1} - \wt{g}^\eps S_\eps b(\D) (\A^0 - \widehat{\zeta}I)^{-1} \|_{L_2(\R^d)\to L_2(\R^d)}
\le C_4 c(\varphi)^2 \eps
\eqno(2.50)
$$
for $\eps >0$, where $C_4=   \|g\|_{L_\infty}^{1/2} \widetilde{C}_2+ 2  \|g\|_{L_\infty} c_0^{-1} \alpha_1^{1/2}
((d\alpha_1)^{1/2} M_1+r_1)$.

By the scaling transformation, (2.50) is equivalent to
$$
\| g b(\D) (\A - \widehat{\zeta} \eps^2 I)^{-1} - \wt{g} S_1 b(\D) (\A^0 - \widehat{\zeta} \eps^2 I)^{-1} \|_{L_2(\R^d)\to L_2(\R^d)}
\le C_4 c(\varphi)^2
\eqno(2.51)
$$
for $\eps >0$. Replacing  $\eps$ by $\eps |\zeta|^{1/2}$ in (2.51) and applying the inverse transformation, we arrive at (2.44).
$\ \bullet$

\smallskip
Now we distinguish the case where the corrector is equal to zero. The next statement follows from Theorem 2.4,
Proposition 1.2, and (1.7).

\smallskip\noindent\textbf{Proposition 2.7.}
\textit{Let} $\u_\eps$ \textit{be the solution of}
(2.2), \textit{and let} $\u_0$ \textit{be the solution of} (2.7).
\textit{If $g^0 = \overline{g}$,
 i.~e., relations} (1.13) \textit{are satisfied, then $\Lambda =0$, $K(\eps;\zeta)=0$,
and we have}
 $$
 \| \D (\u_\eps - \u_0)\|_{L_2(\R^d)} \le C_2 c(\varphi)^2 \eps \|\FF\|_{L_2(\R^d)},\quad \eps >0.
 $$

\smallskip\noindent\textbf{2.3. The results for homogenization problem in $\R^d$ in the case where $\Lambda \in L_\infty$.}
It turns out that under some additional assumptions on the solution of problem (1.7) the smoothing operator
$S_{\varepsilon}$ in (2.15) can be removed.

\smallskip\noindent\textbf{Condition 2.8.}
\textit{Suppose that the} $\Gamma$-\textit{periodic solution}
$\Lambda(\x)$ \textit{of problem} (1.7) \textit{is bounded}:
$\Lambda \in L_\infty$.

We need the following multiplicative property of $\Lambda$, see [PSu2, Corollary 2.4].

\smallskip\noindent\textbf{Proposition 2.9.}
\textit{Under Condition} 2.8 \textit{for any function $u \in H^1(\R^d)$ and $\eps>0$ we have}
$$
\int_{\R^d} |(\D \Lambda)^\eps(\x)|^2 |u|^2\,d\x \le
\beta_1 \|u\|^2_{L_2(\R^d)} + \beta_2 \|\Lambda\|^2_{L_\infty} \eps^2 \int_{\R^d} |\D u|^2 d\x.
$$
\textit{Here}
$$
\begin{aligned}
\beta_1 &= 16 m \alpha_0^{-1} \|g\|_{L_\infty} \|g^{-1}\|_{L_\infty},
\cr
\beta_2 &= 2(1+ 2d \alpha_0^{-1} \alpha_1 + 20 d \alpha_0^{-1} \alpha_1 \|g\|_{L_\infty} \|g^{-1}\|_{L_\infty}).
\end{aligned}
$$

We put
$$
K^0(\eps;\zeta)
= [\Lambda^{\varepsilon}] b (\D) (\A^0 - \zeta I)^{-1}.
\eqno(2.52)
$$
Under Condition 2.8 the operator (2.52) is a continuous mapping of $L_2(\R^d;\C^n)$ to $H^1(\R^d;\C^n)$;
this is seen from the identity
$$
K^0(\eps;\zeta) = K^0(\eps;-1) (\A^0 +I)(\A^0 - \zeta I)^{-1}
\eqno(2.53)
$$
and the following lemma.

\smallskip\noindent\textbf{Lemma 2.10.}
\textit{Suppose that Condition} 2.8 \textit{is satisfied. For $\eps>0$ and $\zeta=-1$ the operator} (2.52) \textit{is continuous from $L_2(\R^d;\C^n)$ to $H^1(\R^d;\C^n)$,
and we have}
$$
\|  K^0(\eps;-1) \|_{L_2(\R^d)\to L_2(\R^d)} \leq \check{C}_K^{(1)},
\eqno(2.54)
$$
$$
\|\eps \D K^0(\eps;-1) \|_{L_2(\R^d)\to L_2(\R^d)} \leq \check{C}_K^{(2)} \eps + \check{C}_K^{(3)}.
\eqno(2.55)
$$
\textit{The constants $\check{C}_K^{(j)}$, $j=1,2,3,$ depend only on $m$, $\alpha_0$, $\alpha_1$, $\|g\|_{L_\infty}$,
$\|g^{-1}\|_{L_\infty}$, the parameters of the lattice} $\Gamma$, \textit{and} $\|\Lambda\|_{L_\infty}$.

\smallskip
\noindent\textbf{Proof.}
By (2.20),
$$
\| K^0(\eps;-1) \|_{L_2 (\mathbb{R}^{d})\to L_2 (\mathbb{R}^{d})}
\leq \frac{1}{2}  \|\Lambda\|_{L_\infty} \alpha_1^{1/2} c_0^{-1/2}=: \check{C}_K^{(1)},
$$
which implies (2.54). Now consider the operators
$$
\eps D_j K^0(\eps;-1) = [(D_j \Lambda)^\eps]  b(\D) (\A^0+I)^{-1} +
\eps [\Lambda^\eps]  D_j b(\D) (\A^0+I)^{-1}.
\eqno(2.56)
$$
The second term in (2.56) is estimated with the help of (2.22):
$$
\sum_{j=1}^d \| \eps \Lambda^\eps  b(\D) D_j (\A^0+I)^{-1}\|^2_{L_2 \to L_2}
\leq \eps^2 \|\Lambda\|^2_{L_\infty} \alpha_1 c_0^{-2}.
\eqno(2.57)
$$
To estimate the first term in (2.56), we apply Proposition 2.9 and (2.20), (2.22):
$$
\begin{aligned}
&\sum_{j=1}^d \| (D_j \Lambda)^\eps  b(\D) (\A^0+I)^{-1}\|^2_{L_2 \to L_2}
\cr
&\leq
\beta_1 \| b(\D) (\A^0 +I)^{-1}\|_{L_2 \to L_2}^2 +
\beta_2 \eps^2 \|\Lambda\|^2_{L_\infty} \|\D b(\D) (\A^0 +I)^{-1}\|_{L_2 \to L_2}^2
\cr
&\leq   \beta_1 \alpha_1 (4c_0)^{-1} + \beta_2 \eps^2 \|\Lambda\|^2_{L_\infty} \alpha_1 c_0^{-2}.
\end{aligned}
\eqno(2.58)
$$
Relations (2.56)--(2.58) imply inequality (2.55)
with $\check{C}_K^{(2)}= (2\alpha_1)^{1/2} c_0^{-1} \left( \beta_2 + 1\right)^{1/2} \|\Lambda\|_{L_\infty}$ and
$\check{C}_K^{(3)}= (\beta_1\alpha_1)^{1/2} (2c_0)^{-1/2}$. $\bullet$

Instead of (2.23), consider another first order approximation of $\u_{\varepsilon}$:
$$
\check{\v}_{\varepsilon} = \u_0 + \varepsilon \Lambda^{\varepsilon} b (\D) \u_0
= (\A^0 - \zeta I)^{-1} \FF + \eps K^0(\eps;\zeta)  \FF.
\eqno(2.59)
$$

\smallskip\noindent\textbf{Theorem 2.11.} \textit{Suppose that the assumptions of Theorem} 2.2 \textit{and Condition} 2.8
\textit{are satisfied. Let} $\check{\v}_\eps$ \textit{be defined by} (2.59).
\textit{Let} $\p_\eps = g^\eps b(\D)\u_\eps$.
\textit{Then for $\eps >0$ we have}
$$
 \| \D(\u_{\varepsilon} - \check{\v}_\eps) \|_{L_2 (\mathbb{R}^{d})}
 \leq C_5 c(\varphi)^2  \varepsilon  \| \FF \|_{L_2 (\mathbb{R}^{d})},
 $$
$$
 \| \u_{\varepsilon} - \check{\v}_\eps \|_{L_2 (\mathbb{R}^{d})}
 \leq C_6 c(\varphi)^2 |\zeta|^{-1/2} \varepsilon  \| \FF \|_{L_2 (\mathbb{R}^{d})},
 $$
$$
 \| \p_{\varepsilon} - \widetilde{g}^\eps b(\D) \u_0 \|_{L_2 (\mathbb{R}^{d})}
 \leq C_7 c(\varphi)^2 \varepsilon  \| \FF \|_{L_2 (\mathbb{R}^{d})},
 $$
\textit{or, in operator terms},
$$
\| \D \left( (\A_{\varepsilon} - \zeta I)^{-1} - (\A^0 - \zeta I)^{-1} -
\varepsilon K^0(\eps;\zeta) \right) \|_{L_2 (\mathbb{R}^{d})
\to L_2 (\mathbb{R}^{d})} \leq \ C_5 c(\varphi)^2 \varepsilon,
\eqno(2.60)
$$
$$
\| (\A_{\varepsilon} - \zeta I)^{-1} - (\A^0 - \zeta I)^{-1} -
\varepsilon  K^0(\eps;\zeta) \|_{L_2 (\mathbb{R}^{d})
\to L_2 (\mathbb{R}^{d})} \leq \ C_6 c(\varphi)^2 |\zeta|^{-1/2} \varepsilon,
\eqno(2.61)
$$
$$
\| g^\eps b(\D)(\A_{\varepsilon} - \zeta I)^{-1} -
\widetilde{g}^\eps b(\D) (\A^0 - \zeta I)^{-1} \|_{L_2 (\mathbb{R}^{d})
\to L_2 (\mathbb{R}^{d})} \leq \ C_7 c(\varphi)^2 \varepsilon.
\eqno(2.62)
$$
\textit{The constants} $C_5$, $C_6$, \textit{and} $C_7$ \textit{depend on}
$d$, $m$, $\alpha_0$, $\alpha_1$, $\| g \|_{L_\infty}$, $\| g^{-1} \|_{L_\infty}$,
\textit{the parameters of the lattice} $\Gamma$,
\textit{and the norm} $\|\Lambda\|_{L_\infty}$.

\smallskip\noindent\textbf{Proof.} The proof is similar to that of Theorem 2.4.
In [BSu4], it was shown that under Condition 2.8 we have
$$
\begin{aligned}
\| \D \left((\A_{\varepsilon} + I)^{-1} - (\A^0 + I)^{-1} -
\varepsilon K^0(\varepsilon;-1)\right) \|_{L_2 (\mathbb{R}^{d})
\to L_2 (\mathbb{R}^{d})} \leq  \check{C}_5 \varepsilon,
\cr
 0< \eps\le 1,
\end{aligned}
\eqno(2.63)
$$
$$
\| (\A_{\varepsilon} + I)^{-1} - (\A^0 + I)^{-1} -
\varepsilon K^0(\varepsilon;-1) \|_{L_2 (\mathbb{R}^{d})\to L_2 (\mathbb{R}^{d})} \leq \ \check{C}_6 \varepsilon,\ 0< \eps\le 1.
\eqno(2.64)
$$
The constants $\check{C}_5$, $\check{C}_6$  depend only on
$m$, $d$, $\alpha_0$, $\alpha_1$, $\| g \|_{L_\infty}$, $\| g^{-1} \|_{L_\infty}$, the parameters of the lattice $\Gamma$,
and the norm $\|\Lambda\|_{L_\infty}$.

For $\eps>1$ estimates are trivial.
By (2.54), the left-hand side of (2.64) does not exceed $2+ \check{C}_K^{(1)} \eps$,
and then does not exceed $(2+ \check{C}_K^{(1)}) \eps$ if $\eps >1$.
Combining this with (2.64), we obtain
$$
\| (\A_{\varepsilon} + I)^{-1} - (\A^0 + I)^{-1} -
\varepsilon K^0(\varepsilon;-1) \|_{L_2 (\mathbb{R}^{d})\to L_2 (\mathbb{R}^{d})} \leq \ \widehat{C}_6 \varepsilon,\ \eps>0,
\eqno(2.65)
$$
where $\widehat{C}_6 = \max \{ \check{C}_6, 2+ \check{C}_K^{(1)}\}$.

By (2.31), (2.32), and (2.55), the left-hand side of (2.63) is estimated by
$2 c_0^{-1/2} + \check{C}_K^{(2)} \eps + \check{C}_K^{(3)}$, which does not exceed
$(2 c_0^{-1/2} + \check{C}_K^{(2)} + \check{C}_K^{(3)})\eps$ for $\eps >1$.
Together with (2.63) this yields 
$$
\| \D \left((\A_{\varepsilon} + I)^{-1} - (\A^0 + I)^{-1} -
\varepsilon K^0(\varepsilon;-1) \right) \|_{L_2 (\mathbb{R}^{d})
\to L_2 (\mathbb{R}^{d})} \leq \ \widehat{C}_5 \varepsilon,\  \eps >0,
\eqno(2.66)
$$
where $\widehat{C}_5 = \max \{ \check{C}_5, 2 c_0^{-1/2} + \check{C}_K^{(2)} + \check{C}_K^{(3)}\}$.

Using the analog of the identity (2.34) with $K(\eps;\zeta)$ replaced by $K^0(\eps;\zeta)$,
we carry over estimates (2.65) and (2.66) to the point $\widehat{\zeta}= e^{i\varphi}$.
Taking (2.11), (2.12), (2.35), (2.53), and (2.54) into account, we deduce the following inequality from (2.65):
$$
\| (\A_\eps - \widehat{\zeta} I)^{-1} - (\A^0 - \widehat{\zeta} I)^{-1} - \eps K^0(\eps;\widehat{\zeta}) \|_{L_2 \to L_2}
\leq C_6 c(\varphi)^2 \eps, \ \ \eps>0,
\eqno(2.67)
$$
where $C_6 = 4 \widehat{C}_6 + 4 \check{C}_K^{(1)}$.

Next, similarly to (2.38)--(2.40) we have:
$$
\| \A_\eps^{1/2} \left( (\A_\eps - \widehat{\zeta} I)^{-1} - (\A^0 - \widehat{\zeta} I)^{-1} -
\eps   K^0(\eps;\widehat{\zeta}) \right) \|_{L_2 \to L_2}
\leq \widetilde{C}_5 c(\varphi)^2 \eps, \ \eps>0,
\eqno(2.68)
$$
where $\widetilde{C}_5 = 4 (c_1^{1/2}\widehat{C}_5 + \check{C}_K^{(1)})$.
Combining this with the lower estimate (1.6), we see that
$$
\| \D \left( (\A_\eps - \widehat{\zeta} I)^{-1} - (\A^0 - \widehat{\zeta} I)^{-1} -
\eps  K^0(\eps;\widehat{\zeta}) \right) \|_{L_2 \to L_2}
\leq {C}_5 c(\varphi)^2 \eps, \ \ \eps>0,
\eqno(2.69)
$$
where $C_5 = \widetilde{C}_5 c_0^{-1/2}$.

Now, by the scaling transformation, estimate (2.61) is deduced from (2.67), and (2.60) follows from (2.69).
(Cf. the proof of Theorem 2.4.)

It remains to check (2.62). By (2.68),
$$
\begin{aligned}
&\| g^\eps b(\D) \left( (\A_\eps - \widehat{\zeta} I)^{-1} - (\A^0 - \widehat{\zeta} I)^{-1} -
\eps  K^0(\eps;\widehat{\zeta}) \right) \|_{L_2 \to L_2}
\cr
&\leq \|g\|^{1/2}_{L_\infty} \widetilde{C}_5 c(\varphi)^2 \eps,\ \ \eps >0.
\end{aligned}
\eqno(2.70)
$$
We have
$$
\begin{aligned}
\eps {g}^\eps b(\D) K^0(\eps;\widehat{\zeta}) &=
g^\eps(b(\D)\Lambda)^\eps b(\D) (\A^0 - \widehat{\zeta}I)^{-1}
\cr
&+
\eps \sum_{l=1}^d g^\eps b_l \Lambda^\eps b(\D) D_l (\A^0 - \widehat{\zeta}I)^{-1}.
\end{aligned}
\eqno(2.71)
$$
The second term in (2.71) is estimated with the help of (1.5) and (2.47):
$$
\eps \left\| \sum_{l=1}^d g^\eps b_l \Lambda^\eps  b(\D) D_l (\A^0 - \widehat{\zeta}I)^{-1} \right\|_{L_2\to L_2}
\leq 2 \eps \|g\|_{L_\infty} \|\Lambda\|_{L_\infty} \alpha_1 c_0^{-1} d^{1/2} c(\varphi).
\eqno(2.72)
$$
Relations (2.70)--(2.72) and (1.9) imply that
$$
\| g^\eps b(\D) (\A_\eps - \widehat{\zeta} I)^{-1}
- \widetilde{g}^\eps b(\D) (\A^0 - \widehat{\zeta} I)^{-1} \|_{L_2\to L_2}
\leq C_7 c(\varphi)^2 \eps,\ \eps>0,
\eqno(2.73)
$$
where $C_7 = \|g\|^{1/2}_{L_\infty} \widetilde{C}_5 + 2 \|g\|_{L_\infty} \|\Lambda\|_{L_\infty} \alpha_1 c_0^{-1} d^{1/2}$.

Inequality (2.62) is deduced from (2.73) by the scaling transformation. (Cf. the proof of Theorem 2.6.) $\bullet$

\smallskip
In some cases Condition 2.8 is valid automatically.
The next statement was proved in [BSu4, Lemma 8.7].

\smallskip\noindent\textbf{Proposition 2.12.} \textit{Condition} 2.8 \textit{is a fortiori valid if at least one of
the following assumptions is fulfilled}:

$1^\circ$. $d \leq 2$;

$2^\circ$. \textit{$d\geq 1$ and} $\A_\eps = \D^* g^\eps(\x) \D$,
\textit{where} $g(\x)$ \textit{has real entries};

$3^\circ$. \textit{dimension is arbitrary, and} $g^0 = \underline{g}$, \textit{i.~e., relations} (1.14) \textit{are satisfied}.

\smallskip
Note that Condition 2.8 is also ensured if $g(\x)$ is smooth enough.

We distinguish the special case where $g^0=\underline{g}$. In this case the matrix
(1.9) is constant: $\wt{g}(\x)=g^0 =\underline{g}$; moreover, Condition 2.8 is satisfied.
Applying the statement of Theorem 2.11 concerning the fluxes, we arrive at the following statement.

\smallskip\noindent\textbf{Proposition 2.13.}
\textit{Suppose that the assumptions of Theorem} 2.2
\textit{are satisfied. Let} $\p_\eps = g^\eps b(\D)\u_\eps$.
\textit{Let $g^0 = \underline{g}$, i.~e., relations} (1.14) \textit{are satisfied. Then}
$$
 \| \p_{\varepsilon} - {g}^0 b(\D) \u_0 \|_{L_2 (\mathbb{R}^{d})}
 \leq C_7 c(\varphi)^2 \varepsilon  \| \FF \|_{L_2 (\mathbb{R}^{d})},
 \quad  \varepsilon >0.
$$

\section*{Chapter 2. The Dirichlet problem}

\section*{\S 3. The Dirichlet problem in a bounded domain. Preliminaries}

\smallskip\noindent\textbf{3.1. The operator $\A_{D,\eps}$.}
Let $\mathcal{O} \subset \mathbb{R}^{d}$ be a bounded domain with the boundary of class $C^{1,1}$.
In $L_2(\O;\C^n)$, we consider the operator $\A_{D,\eps}$ formally given by the differential expression
$b(\D)^* g^\eps(\x) b(\D)$ with the Dirichlet condition on $\partial \O$.
Precisely, $\A_{D,\eps}$ is the selfadjoint operator in
$L_2(\O;\C^n)$ generated by the quadratic form
$$
a_{D,\eps}[\u,\u] = \int_\O \langle g^\eps(\x) b(\D)\u,b(\D)\u \rangle\, d\x,\ \ \u \in H^1_0(\O;\C^n).
$$
This form is closed and positive definite. Indeed, extending $\u$ by zero to $\R^d \setminus \O$
and applying (1.6), we obtain
$$
c_0 \| \D \u \|^2_{L_2(\O)} \leq a_{D,\eps}[\u,\u] \leq c_1 \| \D \u \|^2_{L_2(\O)},\ \ \u \in H^1_0(\O;\C^n).
\eqno(3.1)
$$
It remains to take into account that $\| \D \u \|_{L_2(\O)}$ defines a norm in $H^1_0(\O;\C^n)$ equivalent to the standard one.
By the Friedrichs inequality, (3.1) implies that
$$
a_{D,\eps}[\u,\u] \geq c_2 \|\u\|^2_{L_2(\O)},\ \ \u \in H^1_0(\O;\C^n),\ \ c_2 = c_0 ({\rm diam}\,\O)^{-2}.
\eqno(3.2)
$$

\textit{Our goal} in Chapter 2 is to find approximation for small $\eps$ of the generalized solution $\u_\eps \in H^1_0(\O;\C^n)$
of the Dirichlet problem
$$
b(\D)^* g^\eps(\x) b(\D) \u_\eps(\x) - \zeta \u_\eps(\x) = \FF(\x),\ \ \x \in \O; \ \ \u_\eps\vert_{\partial\O} =0,
\eqno(3.3)
$$
where $\FF \in L_2(\O;\C^n)$. As in \S 2, we assume that $\zeta \in \C \setminus \R_+$.
(The case of other admissible values of $\zeta$ is studied in \S 8.)
We have $\u_\eps = (\A_{D,\eps} - \zeta I)^{-1}\FF$. In operator terms, we study the behavior of the resolvent
$(\A_{D,\eps} - \zeta I)^{-1}$ for small $\eps$.

\smallskip\noindent\textbf{Lemma 3.1.} \textit{Let} $\zeta \in \C \setminus \R_+$.
\textit{Let} $\u_\eps$ \textit{be the generalized solution of the problem} (3.3). \textit{Then for $\eps>0$ we have}
$$
\|\u_\eps \|_{L_2(\O)} \leq c(\varphi) |\zeta|^{-1} \|\FF\|_{L_2(\O)},
\eqno(3.4)
$$
$$
\|\D \u_\eps \|_{L_2(\O)} \leq {\mathcal C}_0 c(\varphi) |\zeta|^{-1/2} \|\FF\|_{L_2(\O)}.
\eqno(3.5)
$$
\textit{Here ${\mathcal C}_0= \sqrt{2} \alpha_0^{-1/2} \|g^{-1}\|_{L_\infty}^{1/2}$. In operator terms,},
$$
\| (\A_{D,\eps} - \zeta I)^{-1} \|_{L_2(\O)\to L_2(\O)} \leq c(\varphi) |\zeta|^{-1},
\eqno(3.6)
$$
$$
\|\D (\A_{D,\eps} - \zeta I)^{-1} \|_{L_2(\O)\to L_2(\O)} \leq {\mathcal C}_0 c(\varphi) |\zeta|^{-1/2}.
$$

 \smallskip\noindent\textbf{Proof.}
 By (3.2), the spectrum of $\A_{D,\eps}$ is contained in $[c_2,\infty) \subset \R_+$.
The norm of the resolvent $(\A_{D,\eps} - \zeta I)^{-1}$ does not exceed the inverse distance from $\zeta$ to $\R_+$.
As a result, we arrive at (3.6).

To check (3.5), we write down the integral identity for the solution $\u_\eps \in H^1_0(\O;\C^n)$ of the problem (3.3):
$$
(g^\eps b(\D)\u_\eps, b(\D) {\boldsymbol \eta})_{L_2(\O)} - \zeta (\u_\eps, {\boldsymbol \eta})_{L_2(\O)}=
(\FF, {\boldsymbol \eta})_{L_2(\O)},\ \  {\boldsymbol \eta}\in H^1_0(\O;\C^n).
\eqno(3.7)
$$
Substituting  ${\boldsymbol \eta}=\u_\eps$ in (3.7) and taking the real part, we obtain
$$
(g^\eps b(\D)\u_\eps, b(\D) \u_\eps)_{L_2(\O)} = {\rm Re}\, \zeta \|\u_\eps\|^2_{L_2(\O)} +
{\rm Re}\, (\FF, \u_\eps)_{L_2(\O)}.
$$
Together with (3.4) this implies that
$$
(g^\eps b(\D)\u_\eps, b(\D) \u_\eps)_{L_2(\O)} \leq 2 c(\varphi)^2 |\zeta|^{-1} \|\FF\|^2_{L_2(\O)}.
$$
Taking  (3.1) into account, we obtain (3.5) with ${\mathcal C}_0 = \sqrt{2}c_0^{-1/2}$.
$\bullet$

\smallskip\noindent\textbf{3.2. The effective operator $\A^0_{D}$.}
In $L_2(\O;\C^n)$, consider the selfadjoint operator $\A^0_{D}$ generated by the quadratic form
$$
a^0_{D}[\u,\u] = \int_\O \langle g^0 b(\D)\u,b(\D)\u \rangle\, d\x,\ \ \u \in H^1_0(\O;\C^n).
\eqno(3.8)
$$
Here $g^0$ is the effective matrix defined by (1.8).
Taking (1.15) into account, we see that the form (3.8) satisfies estimates of the form (3.1) and (3.2) with the same constants.

Since $\partial \O \in C^{1,1}$, the operator $\A_D^0$ is given by
$b(\D)^* g^0 b(\D)$ on the domain $H^2(\O;\C^n) \cap H^1_0(\O;\C^n)$, and we have
$$
\| (\A^0_D)^{-1}\|_{L_2(\O) \to H^2(\O)} \leq \widehat{c}.
\eqno(3.9)
$$
The constant $\widehat{c}$ depends only on $\alpha_0$, $\alpha_1$, $\|g\|_{L_\infty}$, $\|g^{-1}\|_{L_\infty}$, and the domain $\O$.
To justify this fact, we note that the operator $b(\D)^* g^0 b(\D)$ is a \textit{strongly elliptic} matrix operator
with constant coefficients and refer to the theorems about regularity of solutions for strongly elliptic systems
(see, e.~g., [McL, Chapter 4]).

Let $\u_0\in H^1_0(\O;\C^n)$ be the generalized solution of the problem
$$
b(\D)^* g^0 b(\D) \u_0(\x) - \zeta \u_0(\x) = \FF(\x),\ \ \x \in \O; \ \ \u_0\vert_{\partial\O} =0.
\eqno(3.10)
$$
Then $\u_0 = (\A^0_{D} - \zeta I)^{-1}\FF$.

\smallskip\noindent\textbf{Lemma 3.2.}
\textit{Let} $\zeta \in \C \setminus \R_+$.
\textit{Let} $\u_0$ \textit{be the solution of the problem} (3.10). \textit{Then we have}
$$
\|\u_0 \|_{L_2(\O)} \leq c(\varphi) |\zeta|^{-1} \|\FF\|_{L_2(\O)},
$$
$$
\|\D \u_0 \|_{L_2(\O)} \leq {\mathcal C}_0 c(\varphi) |\zeta|^{-1/2} \|\FF\|_{L_2(\O)},
$$
$$
\|\u_0 \|_{H^2(\O)} \leq  \widehat{c} c(\varphi) \|\FF\|_{L_2(\O)},
\eqno(3.11)
$$
\textit{or, in operator terms},
$$
\| (\A^0_{D} - \zeta I)^{-1} \|_{L_2(\O)\to L_2(\O)} \leq c(\varphi) |\zeta|^{-1},
\eqno(3.12)
$$
$$
\|\D (\A^0_{D} - \zeta I)^{-1} \|_{L_2(\O)\to L_2(\O)} \leq {\mathcal C}_0 c(\varphi) |\zeta|^{-1/2},
\eqno(3.13)
$$
$$
\|(\A^0_{D} - \zeta I)^{-1} \|_{L_2(\O)\to H^2(\O)} \leq \widehat{c} c(\varphi).
\eqno(3.14)
$$

 \smallskip\noindent\textbf{Proof.}
Estimates (3.12), (3.13) are checked similarly to the proof of Lemma 3.1.
Estimate (3.14) follows from (3.9) and the inequality
$$
\|\A_D^0 (\A^0_{D} - \zeta I)^{-1} \|_{L_2(\O)\to L_2(\O)} \leq
\sup_{x\geq 0} x |x-\zeta|^{-1} \leq  c(\varphi).
$$
$\bullet$

\smallskip\noindent\textbf{3.3. Auxiliary statements: estimates in the neighborhood of the boundary.}
This subsection contains several auxiliary statements (see, e.~g., [PSu2, \S 5]).

\smallskip\noindent\textbf{Lemma 3.3.}
\textit{Let} $\O \subset \R^d$ \textit{be a bounded domain with the boundary of class $C^1$. Denote
$(\partial \O)_\eps = \{\x \in \R^d: {\rm dist}\,\{\x;\partial {\mathcal O} \}< \eps \}$,
$B_\eps = (\partial \O)_\eps \cap \O$. Suppose that the number $\eps_0>0$ is such that
$(\partial \O)_{\eps_0}$ can be covered by a finite number of open sets admitting diffeomorphisms of class $C^1$
rectifying the boundary. Then the following statements are true.}

\noindent $1^\circ$. \textit{For any function $u \in H^1(\O)$ we have}
$$
\int_{B_\eps} |u|^2\,d\x \leq \beta \eps \|u\|_{H^1(\O)} \|u\|_{L_2(\O)},\ \ 0 <  \eps \leq \eps_0.
$$
$2^\circ$. \textit{For any function $u \in H^1(\R^d)$ we have}
$$
\int_{({\partial \O})_\eps} |u|^2\,d\x \leq \beta \eps \|u\|_{H^1(\R^d)} \|u\|_{L_2(\R^d)},\ \ 0 <  \eps \leq \eps_0.
$$
\textit{The constant $\beta$ depends only on the domain $\O$.}

\smallskip\noindent\textbf{Lemma 3.4.}
\textit{Suppose that the assumptions of Lemma} 3.3 \textit{are satisfied.
Let $f(\x)$ be a $\Gamma$-periodic function in $\R^d$ such that
$f \in L_2(\Omega)$. Let $S_\eps$ be the operator} (1.17).
\textit{Then for $0 <  \eps \leq \eps_1$ and any function $\u \in H^1(\R^d;\C^m)$ we have}
$$
\int_{({\partial \O})_\eps} |f^\eps(\x)|^2 |(S_\eps \u)|^2\,d\x
\leq \beta_*  \eps |\Omega|^{-1} \|f\|_{L_2(\Omega)}^2 \|\u\|_{H^1(\R^d)} \|\u\|_{L_2(\R^d)},
$$
\textit{where $\eps_1 = \eps_0 (1+r_1)^{-1}$, $\beta_* = \beta (1+r_1)$, $2r_1 = {\rm diam}\,\Omega$.}

\section*{\S 4. The results for the Dirichlet problem}

Assume that the numbers $\eps_0, \eps_1 \in (0,1]$ are chosen according to the following condition.

\smallskip\noindent\textbf{Condition 4.1.} \textit{Suppose that the number} $\eps_0 \in (0,1]$ \textit{is such that
$(\partial \O)_{\eps_0}$ can be covered by a finite number of open sets admitting diffeomorphisms of class $C^{1,1}$
rectifying the boundary $\partial \O$. Let $\eps_1 = \eps_0 (1+r_1)^{-1}$, where $2r_1 = {\rm diam}\, \Omega$.}

 \smallskip
Clearly, $\eps_1$ depends only on the domain $\O$ and the lattice $\Gamma$.

\smallskip\noindent\textbf{4.1. Approximation of the resolvent $(\A_{D,\eps}-\zeta I)^{-1}$ for $|\zeta|\geq 1$.}
Now we formulate our main results for the operator $\A_{D,\eps}$.
The case of large $|\zeta|$ is of main interest, so
for a while we assume that $|\zeta|\geq 1$. (The case of other admissible values of $\zeta$ will be studied in \S 8.)

\smallskip\noindent\textbf{Theorem 4.2.} \textit{Let $\zeta = |\zeta| e^{i \varphi} \in {\mathbb C}\setminus {\mathbb R}_+$ and $|\zeta|\ge 1$.
Let $\u_\eps$ be the solution of problem} (3.3), \textit{and let $\u_0$ be the solution of problem} (3.10) \textit{with $\FF \in L_2(\O;\C^n)$.
Suppose that $\varepsilon_1$ is subject to Condition} 4.1.
\textit{Then for $0< \varepsilon \leq \varepsilon_1$ we have}
$$
\| \u_\eps - \u_0 \|_{L_2(\O)} \leq {\mathcal C}_1 c(\varphi)^5 (|\zeta|^{-1/2}\varepsilon + \varepsilon^2) \|\FF\|_{L_2(\O)},
\eqno(4.1)
$$
\textit{where $c(\varphi)$ is defined by} (2.1). \textit{In operator terms},
$$
\| (\A_{D,\varepsilon}- \zeta I)^{-1}  - (\A^0_{D}- \zeta I)^{-1} \|_{L_2({\mathcal O}) \to L_2({\mathcal O})}
\le {\mathcal C}_{1} c(\varphi)^5 (|\zeta|^{-1/2} \varepsilon+ \varepsilon^2).
\eqno(4.2)
$$
\textit{The constant ${\mathcal C}_1$ depends on}
$d$, $m$, $\alpha_0$, $\alpha_1$, $\| g \|_{L_\infty}$, $\| g^{-1} \|_{L_\infty}$,
\textit{the parameters of the lattice} $\Gamma$, \textit{and the domain} $\O$.

\smallskip
In order to approximate the solution in $H^1(\O;\C^n)$, we need to introduce the corrector.
We fix a linear continuous extension operator
$P_{\mathcal{O}}: H^2(\mathcal{O};\mathbb{C}^n) \to H^2(\mathbb{R}^d;\mathbb{C}^n).$
Suppose that simultaneously $P_{\mathcal{O}}$ is continuous from
$L_2(\mathcal{O};\mathbb{C}^n)$ to $L_2(\mathbb{R}^d;\mathbb{C}^n)$ and from
$H^1(\mathcal{O};\mathbb{C}^n)$ to $H^1(\mathbb{R}^d;\mathbb{C}^n)$.
Such an operator exists (see, e.~g., [St]). Denote
$$
\|P_{\mathcal{O}} \|_{H^s(\O) \to H^s(\R^d)} =: C_\O^{(s)},\ \ s=0,1,2.
\eqno(4.3)
$$
The constants $C_\O^{(s)}$ depend only on the domain $\O$.
Next, by $R_{\mathcal{O}}$ we denote the operator of restriction
of functions in $\mathbb{R}^d$ onto the domain $\mathcal{O}$. We introduce the corrector
$$
K_D(\varepsilon;\zeta)
= R_{\mathcal{O}} [\Lambda^\varepsilon] S_\varepsilon b(\mathbf{D}) P_{\mathcal{O}} (\A_D^0 -\zeta I)^{-1}.
\eqno(4.4)
$$
The operator $K_D(\varepsilon;\zeta)$ is a continuous mapping of $L_2(\mathcal{O};\mathbb{C}^n)$ to $H^1(\mathcal{O};\mathbb{C}^n)$.
Indeed, the operator $b(\mathbf{D}) P_{\mathcal{O}} (\A_D^0 -\zeta I)^{-1}$
is continuous from $L_2(\mathcal{O};\mathbb{C}^n)$ to $H^1(\R^d;\C^m)$,
and the operator $[\Lambda^\varepsilon] S_\varepsilon$ is a continuous mapping of $H^1(\R^d;\C^m)$ to $H^1(\R^d;\C^n)$
(this follows from Proposition 1.5 and relation $\Lambda \in \widetilde{H}^1(\Omega)$).

Let $\u_0$ be the solution of problem (3.10). Denote $\widetilde{\u}_0:= P_\O \u_0$.
Consider the following function in $\R^d$:
$$
\widetilde{\v}_\eps(\x) = \widetilde{\u}_0(\x) + \eps \Lambda^\eps(\x) (S_\eps b(\D) \widetilde{\u}_0)(\x),
\eqno(4.5)
$$
and put
$$
 \v_\eps :=\widetilde{\v}_\eps\vert_\O.
 \eqno(4.6)
$$
Then
$$
\v_\eps = \u_0 + \eps K_D(\eps;\zeta) \FF = (\A_D^0 - \zeta I)^{-1}\FF + \eps K_D(\eps;\zeta) \FF.
\eqno(4.7)
$$

\smallskip\noindent\textbf{Theorem 4.3.} \textit{Suppose that the assumptions of Theorem} 4.2
 \textit{are satisfied. Let $\v_\eps$ be defined by} (4.4), (4.7).
\textit{Then for $0< \varepsilon \leq \varepsilon_1$ we have}
$$
\| \u_\eps - \v_\eps \|_{H^1(\O)} \leq \left({\mathcal C}_2 c(\varphi)^2  |\zeta|^{-1/4}\varepsilon^{1/2} +
{\mathcal C}_3 c(\varphi)^4 \varepsilon\right) \|\FF\|_{L_2(\O)},
\eqno(4.8)
$$
\textit{or, in operator terms},
$$
\begin{aligned}
&\| (\A_{D,\varepsilon}- \zeta I)^{-1}  - (\A^0_{D}- \zeta I)^{-1} - \eps K_D(\eps;\zeta) \|_{L_2({\mathcal O}) \to H^1({\mathcal O})}
\cr
&\leq {\mathcal C}_{2} c(\varphi)^2 |\zeta|^{-1/4} \varepsilon^{1/2}
+ {\mathcal C}_3 c(\varphi)^4 \varepsilon.
\end{aligned}
\eqno(4.9)
$$
\textit{For the flux} $\p_\eps := g^\eps b(\D) \u_\eps$ \textit{we have}
$$
\| \p_\eps - \widetilde{g}^\eps S_\eps b(\D) \widetilde{\u}_0 \|_{L_2(\O)}
\leq \left(\widetilde{\mathcal C}_2 c(\varphi)^2 |\zeta|^{-1/4}\varepsilon^{1/2} +
\widetilde{\mathcal C}_3 c(\varphi)^4 \varepsilon\right) \|\FF\|_{L_2(\O)}
\eqno(4.10)
$$
\textit{for $0< \eps \leq \eps_1$. The constants ${\mathcal C}_2$, ${\mathcal C}_3$, $\widetilde{\mathcal C}_2$, and $\widetilde{\mathcal C}_3$ depend on}
$d$, $m$, $\alpha_0$, $\alpha_1$, $\| g \|_{L_\infty}$, $\| g^{-1} \|_{L_\infty}$,
\textit{the parameters of the lattice} $\Gamma$, \textit{and the domain} $\O$.

\smallskip\noindent\textbf{4.2. The first step of the proof. Associated problem in $\R^d$.}
The proof of Theorems 4.2 and 4.3 relies on application of the results for the problem in $\R^d$
and introduction of the boundary layer correction term.

By Lemma 3.2 and (4.3), we have
$$
\| \widetilde{\u}_0 \|_{L_2(\R^d)} \leq C_\O^{(0)} c(\varphi) |\zeta|^{-1} \|\FF\|_{L_2(\O)},
\eqno(4.11)
$$
$$
\| \widetilde{\u}_0 \|_{H^1(\R^d)} \leq C_\O^{(1)} ({\mathcal C}_0+1) c(\varphi) |\zeta|^{-1/2} \|\FF\|_{L_2(\O)},
\eqno(4.12)
$$
$$
\| \widetilde{\u}_0 \|_{H^2(\R^d)} \leq C_\O^{(2)} \widehat{c} c(\varphi) \|\FF\|_{L_2(\O)}.
\eqno(4.13)
$$
We have taken into account that $|\zeta| \geq 1$. We put
$$
\widetilde{\FF}:= \A^0 \widetilde{\u}_0 - \zeta \widetilde{\u}_0.
\eqno(4.14)
$$
Then $\widetilde{\FF}\in L_2(\R^d;\C^n)$ and $\widetilde{\FF}\vert_{\mathcal O} =\FF$.
From (1.16), (4.11), and (4.13) it follows that
$$
\| \widetilde{\FF}\|_{L_2(\R^d)} \leq c_1 \|\widetilde{\u}_0\|_{H^2(\R^d)} + |\zeta| \|\widetilde{\u}_0\|_{L_2(\R^d)}
\leq {\mathcal C}_4 c(\varphi) \|\FF\|_{L_2(\O)},
\eqno(4.15)
$$
where ${\mathcal C}_4= c_1 C_\O^{(2)} \widehat{c} + C_\O^{(0)}$.

Let $\widetilde{\u}_\eps \in H^1(\R^d;\C^n)$ be the solution of the following equation in $\R^d$:
$$
\A_\eps \widetilde{\u}_\eps - \zeta \widetilde{\u}_\eps = \widetilde{\FF},
\eqno(4.16)
$$
i.~e., $\widetilde{\u}_\eps = (\A_\eps - \zeta I)^{-1} \widetilde{\FF}$.
We can apply theorems of \S 2.
From Theorem 2.2 and (4.14)--(4.16) it follows that for $\eps >0$ we have
$$
\| \widetilde{\u}_\eps - \widetilde{\u}_0\|_{L_2(\R^d)} \leq
C_1  c(\varphi)^2 |\zeta|^{-1/2} \eps \|\widetilde{\FF}\|_{L_2(\R^d)}
\leq
C_1 {\mathcal C}_4 c(\varphi)^3 |\zeta|^{-1/2} \eps \|{\FF}\|_{L_2(\O)}.
\eqno(4.17)
$$
By Theorem 2.4 and (4.5), (4.14)--(4.16),  for $\eps>0$ we have
$$
\| \D (\widetilde{\u}_\eps - \widetilde{\v}_\eps) \|_{L_2(\R^d)} \leq
C_2  c(\varphi)^2 \eps \|\widetilde{\FF}\|_{L_2(\R^d)}
\leq
C_2 {\mathcal C}_4 c(\varphi)^3 \eps \|{\FF}\|_{L_2(\O)},
\eqno(4.18)
$$
$$
\| \widetilde{\u}_\eps - \widetilde{\v}_\eps \|_{L_2(\R^d)} \leq
C_3  c(\varphi)^2 |\zeta|^{-1/2} \eps \|\widetilde{\FF}\|_{L_2(\R^d)}
\leq
C_3 {\mathcal C}_4 c(\varphi)^3 |\zeta|^{-1/2} \eps \|{\FF}\|_{L_2(\O)}.
\eqno(4.19)
$$

\smallskip\noindent\textbf{4.3. The second step of the proof. Introduction of the correction term $\w_\eps$.}
The first order approximation $\v_\eps$ of the solution $\u_\eps$ does not satisfy the Dirichlet condition.
We have  $\v_\eps\vert_{\partial \O} = \eps \Lambda^\eps (S_\eps b(\D) \widetilde{\u}_0)\vert_{\partial \O}$.
Consider the "`correction term"' $\w_\eps$ which is the generalized solution of the problem
$$
b(\D)^* g^\eps(\x)b(\D) \w_\eps - \zeta \w_\eps = 0\ \ {\text{in}}\ \O;
\ \ \w_\eps\vert_{\partial \O} = \eps \Lambda^\eps (S_\eps b(\D) \widetilde{\u}_0)\vert_{\partial \O}.
\eqno(4.20)
$$
Here equation is understood in the weak sense: $\w_\eps \in H^1(\O;\C^n)$ satisfies the identity
$$
(g^\eps b(\D)\w_\eps, b(\D){\boldsymbol \eta})_{L_2(\O)} - \zeta(\w_\eps,{\boldsymbol \eta})_{L_2(\O)}=0,\ \ {\boldsymbol \eta}\in H^1_0(\O;\C^n).
\eqno(4.21)
$$

\smallskip\noindent\textbf{Lemma 4.4.} \textit{Let $\u_\eps$ be the solution of problem} (3.3), \textit{and let $\v_\eps$ be given by} (4.7).
\textit{Let $\w_\eps$ be the solution of problem} (4.20). \textit{Then for $\eps>0$ we have}
$$
\| \D (\u_\eps - \v_\eps + \w_\eps) \|_{L_2({\mathcal O})} \leq {\mathcal C}_5 c(\varphi)^4  \eps \|\FF\|_{L_2({\mathcal O})},
\eqno(4.22)
$$
$$
\| \u_\eps - \v_\eps + \w_\eps \|_{L_2({\mathcal O})} \leq {\mathcal C}_6 c(\varphi)^4 |\zeta|^{-1/2} \eps \|\FF\|_{L_2({\mathcal O})}.
\eqno(4.23)
$$
\textit{The constants ${\mathcal C}_5$, ${\mathcal C}_6$
depend on} $d$, $m$, $\alpha_0$, $\alpha_1$, $\| g \|_{L_\infty}$, $\| g^{-1} \|_{L_\infty}$,
\textit{the parameters of the lattice} $\Gamma$, \textit{and the domain} $\O$.

\smallskip\noindent\textbf{Proof.} Denote $\VV_\eps:= \u_\eps - \v_\eps + \w_\eps$. By (3.7), (4.20), and (4.21),
the function $\VV_\eps$ belongs to $H^1_0(\O;\C^n)$ and satisfies the identity
$$
\begin{aligned}
&J_\eps[{\boldsymbol \eta}]:=(g^\eps b(\D)\VV_\eps, b(\D){\boldsymbol \eta})_{L_2(\O)} - \zeta(\VV_\eps,{\boldsymbol \eta})_{L_2(\O)}
\cr
&=
(\FF, {\boldsymbol \eta})_{L_2(\O)} - (g^\eps b(\D)\v_\eps, b(\D){\boldsymbol \eta})_{L_2(\O)} + \zeta(\v_\eps,{\boldsymbol \eta})_{L_2(\O)},
\ \ {\boldsymbol \eta}\in H^1_0(\O;\C^n).
\end{aligned}
\eqno(4.24)
$$
Extend ${\boldsymbol \eta}$ by zero to $\R^d\setminus \O$ keeping the same notation; then ${\boldsymbol \eta}\in H^1(\R^d;\C^n)$.
Recall that $\widetilde{\FF}$ is an extension of $\FF$, and $\widetilde{\v}_\eps$ is an extension of  $\v_\eps$.
Hence,
$$
J_\eps[{\boldsymbol \eta}]=
(\widetilde{\FF},{\boldsymbol \eta})_{L_2(\R^d)} - (g^\eps b(\D)\widetilde{\v}_\eps, b(\D){\boldsymbol \eta})_{L_2(\R^d)}
 +\zeta(\widetilde{\v}_\eps,{\boldsymbol \eta})_{L_2(\R^d)}.
$$
Taking (4.16) into account, we obtain
$$
J_\eps[{\boldsymbol \eta}]=
(g^\eps b(\D)(\widetilde{\u}_\eps - \widetilde{\v}_\eps), b(\D){\boldsymbol \eta})_{L_2(\R^d)} -
\zeta(\widetilde{\u}_\eps -\widetilde{\v}_\eps,{\boldsymbol \eta})_{L_2(\R^d)}.
$$
From (1.4), (4.18), and  (4.19) it follows that
$$
|J_\eps[{\boldsymbol \eta}]| \leq
\eps c(\varphi)^3 \|\FF\|_{L_2(\O)} \left({\mathcal C}_7 \|(g^\eps)^{1/2} b(\D){\boldsymbol \eta}\|_{L_2(\O)} + {\mathcal C}_8 |\zeta|^{1/2} \|{\boldsymbol \eta}\|_{L_2(\O)}\right)
\eqno(4.25)
$$
for $\eps >0$, where ${\mathcal C}_7  = \|g\|_{L_\infty}^{1/2} \alpha_1^{1/2} C_2 {\mathcal C}_4$ and
${\mathcal C}_8= C_3 {\mathcal C}_4$.

Substitute ${\boldsymbol \eta} = \VV_\eps$ in (4.24), take the imaginary part of the corresponding relation, and apply (4.25):
$$
\begin{aligned}
&|{\rm Im}\,\zeta| \|\VV_\eps \|^2_{L_2(\O)} \leq \eps c(\varphi)^3 \|\FF\|_{L_2(\O)}
\cr
&\times \left({\mathcal C}_7
\| (g^\eps)^{1/2} b(\D) \VV_\eps \|_{L_2(\O)} + {\mathcal C}_8 |\zeta|^{1/2} \| \VV_\eps \|_{L_2(\O)}\right).
\end{aligned}
\eqno(4.26)
$$
If ${\rm Re}\,\zeta \geq 0$ (and then ${\rm Im}\,\zeta \ne 0$), this impies the inequality
$$
\begin{aligned}
 \|\VV_\eps \|^2_{L_2(\O)} \leq& 2 \eps |\zeta|^{-1} c(\varphi)^4 {\mathcal C}_7 \|\FF\|_{L_2(\O)} \| (g^\eps)^{1/2} b(\D) \VV_\eps \|_{L_2(\O)}
 \cr
 &+ \eps^2 |\zeta|^{-1} c(\varphi)^8 {\mathcal C}^2_8 \| \FF \|^2_{L_2(\O)}.
\end{aligned}
\eqno(4.27)
$$
If ${\rm Re}\,\zeta < 0$, we take the real part of the corresponding relation and apply (4.25). Then we have
$$
\begin{aligned}
&|{\rm Re}\,\zeta| \|\VV_\eps \|^2_{L_2(\O)} \leq \eps c(\varphi)^3 \|\FF\|_{L_2(\O)}
\cr
&\times \left({\mathcal C}_7
\| (g^\eps)^{1/2} b(\D) \VV_\eps \|_{L_2(\O)} + {\mathcal C}_8 |\zeta|^{1/2} \| \VV_\eps \|_{L_2(\O)}\right).
\end{aligned}
\eqno(4.28)
$$
Summing up (4.26) and (4.28), we deduce the inequality similar to (4.27).
As a result, for all values of $\zeta$ under consideration we obtain
$$
\begin{aligned}
 \|\VV_\eps \|^2_{L_2(\O)} \leq& 4 \eps |\zeta|^{-1} c(\varphi)^4 {\mathcal C}_7 \|\FF\|_{L_2(\O)} \| (g^\eps)^{1/2} b(\D) \VV_\eps \|_{L_2(\O)}
 \cr
 &+ 4 \eps^2 |\zeta|^{-1} c(\varphi)^8 {\mathcal C}^2_8 \| \FF \|^2_{L_2(\O)}.
\end{aligned}
\eqno(4.29)
$$

Now, (4.24) with ${\boldsymbol \eta}=\VV_\eps$, (4.25), and (4.29) yield
$$
\begin{aligned}
&a_{D,\eps}[\VV_\eps, \VV_\eps]
\cr
&\leq
9 {\mathcal C}_7 c(\varphi)^4 \eps \|\FF\|_{L_2(\O)}\| (g^\eps)^{1/2} b(\D) \VV_\eps \|_{L_2(\O)}
+ 9 {\mathcal C}_8^2 c(\varphi)^8 \eps^2\|\FF\|_{L_2(\O)}^2.
\end{aligned}
$$
From here we deduce the inequality
$$
a_{D,\eps}[\VV_\eps, \VV_\eps] \leq \check{\mathcal C}_5^2 c(\varphi)^8 \eps^2\|\FF\|_{L_2(\O)}^2,
\eqno(4.30)
$$
where $\check{\mathcal C}_5^2 = 18 {\mathcal C}_8^2 + 81 {\mathcal C}_7^2$.
By (4.30) and (3.1), we obtain (4.22) with ${\mathcal C}_5 = \check{\mathcal C}_5 c_0^{-1/2}$.
Finally, (4.29) and (4.30) imply (4.23) with
${\mathcal C}_6= 2 \left( {\mathcal C}_7 \check{\mathcal C}_5 + {\mathcal C}_8^2\right)^{1/2}$.
$\bullet$

\smallskip
\textbf{Conclusions.} 1) From (4.22) and (4.23) it follows that
$$
\|\D(\u_\eps - \v_\eps)\|_{L_2(\O)} \leq
{\mathcal C}_5 c(\varphi)^4 \eps \|\FF\|_{L_2(\O)} + \|\D \w_\eps \|_{L_2(\O)},\ \ \eps >0,
\eqno(4.31)
$$
$$
\| \u_\eps - \v_\eps \|_{L_2(\O)} \leq
{\mathcal C}_6 c(\varphi)^4 |\zeta|^{-1/2} \eps \|\FF\|_{L_2(\O)} + \|\w_\eps \|_{L_2(\O)},\ \ \eps >0.
\eqno(4.32)
$$
Thus, for the proof of Theorem 4.3 it suffices to obtain an appropriate estimate for $\|\w_\eps\|_{H^1(\O)}$.

2) Note that the difference $\v_\eps - \u_0 = \eps (\Lambda^\eps S_\eps b(\D) \widetilde{\u}_0)\vert_\O$
can be estimated by using (1.4), (1.19), and (4.12):
$$
\| \v_\eps - \u_0\|_{L_2(\O)} \leq \| \widetilde{\v}_\eps - \widetilde{\u}_0\|_{L_2(\R^d)} \leq
{\mathcal C}_9 c(\varphi) |\zeta|^{-1/2} \eps \|\FF\|_{L_2(\O)},
\eqno(4.33)
$$
where ${\mathcal C}_9 = M_1 \alpha_1^{1/2} C_\O^{(1)} ({\mathcal C}_0 +1)$.
From (4.32) and (4.33) it follows that
$$
\| \u_\eps - \u_0 \|_{L_2(\O)} \leq
({\mathcal C}_6 + {\mathcal C}_9)c(\varphi)^4 |\zeta|^{-1/2} \eps \|\FF\|_{L_2(\O)} + \|\w_\eps \|_{L_2(\O)},\ \ \eps >0.
\eqno(4.34)
$$
Therefore, for the proof of Theorem 4.2 we have to obtain an appropriate estimate for $\|\w_\eps\|_{L_2(\O)}$.

\section*{\S 5. Estimates for the correction term. \\ Proof of Theorems 4.2 and 4.3}

First we estimate the norm of $\w_\eps$ in $H^1(\O;\C^n)$ and complete the proof of Theorem 4.3.
Next, using the already proved Theorem 4.3 and the duality arguments, we estimate the norm of $\w_\eps$ in $L_2(\O;\C^n)$
and prove Theorem 4.2.

\smallskip\noindent\textbf{5.1. Localization near the boundary.}
Suppose that the number $\eps_0 \in (0,1]$ is chosen according to Condition 4.1.
Let $0< \eps \leq \eps_0$.
We fix a smooth cut-off function $\theta_\eps(\x)$ in $\R^d$  such that
$$
\begin{aligned}
&\theta_\eps \in C_0^\infty(\R^d),\ \ {\rm supp}\, \theta_\eps \subset (\partial \O)_\eps,\ \ 0\leq \theta_\eps(\x) \leq 1,
\cr
&\theta_\eps(\x)=1 \ {\text{for}}\ \x \in \partial \O;\ \ \eps |\nabla \theta_\eps(\x)| \leq \kappa = {\rm Const}.
\end{aligned}
\eqno(5.1)
$$
The constant $\kappa$ depends only on the domain $\O$.
Consider the following function in $\R^d$ :
$$
{\boldsymbol \phi}_\eps(\x)= \eps \theta_\eps(\x) \Lambda^\eps(\x) (S_\eps b(\D) \widetilde{\u}_0)(\x).
\eqno(5.2)
$$

\smallskip\noindent\textbf{Lemma 5.1.} \textit{Let $\w_\eps$ be the solution of problem} (4.20). \textit{Let ${\boldsymbol \phi}_\eps$ 
be defined by} (5.2). \textit{Then for $0< \eps \leq \eps_0$ and $\zeta \in \C \setminus \R_+$, $|\zeta|\geq 1$, we have}
$$
\| \w_\eps \|_{H^1({\mathcal O})} \leq c(\varphi)\left({\mathcal C}_{10} |\zeta|^{1/2} \| {\boldsymbol \phi}_\eps\|_{L_2({\mathcal O})}+
{\mathcal C}_{11} \| \D {\boldsymbol \phi}_\eps\|_{L_2({\mathcal O})} \right).
\eqno(5.3)
$$
\textit{The constants ${\mathcal C}_{10}$, ${\mathcal C}_{11}$
depend on} $d$, $m$, $\alpha_0$, $\alpha_1$, $\| g \|_{L_\infty}$, $\| g^{-1} \|_{L_\infty}$,
\textit{the parameters of the lattice} $\Gamma$, \textit{and the domain} $\O$.

\smallskip\noindent\textbf{Proof.}
By (4.20), (5.1), and (5.2), we have
$\w_\eps\vert_{\partial \O}= {\boldsymbol \phi}_\eps\vert_{\partial \O}$.
By (4.21), the function $\w_\eps - {\boldsymbol \phi}_\eps \in H^1_0(\O;\C^n)$ satisfies the identity
$$
\begin{aligned}
&(g^\eps b(\D)(\w_\eps - {\boldsymbol \phi}_\eps), b(\D) {\boldsymbol \eta} )_{L_2(\O)}
- \zeta (\w_\eps - {\boldsymbol \phi}_\eps, {\boldsymbol \eta} )_{L_2(\O)}
\cr
&=
- (g^\eps b(\D) {\boldsymbol \phi}_\eps, b(\D) {\boldsymbol \eta} )_{L_2(\O)}
+ \zeta ( {\boldsymbol \phi}_\eps, {\boldsymbol \eta} )_{L_2(\O)},
\ \ {\boldsymbol \eta} \in H^1_0(\O;\C^n).
\end{aligned}
\eqno(5.4)
$$

We substitute ${\boldsymbol \eta} = \w_\eps - {\boldsymbol \phi}_\eps$ in (5.4) and take the imaginary part of the corresponding
relation. Taking (1.2) and (1.5) into account, we have
$$
\begin{aligned}
|{\rm Im}\,\zeta| \|\w_\eps - {\boldsymbol \phi}_\eps \|^2_{L_2(\O)} &\leq
{\mathcal C}_{12} \| \D {\boldsymbol \phi}_\eps\|_{L_2(\O)}
\|(g^\eps)^{1/2} b(\D)(\w_\eps - {\boldsymbol \phi}_\eps)\|_{L_2(\O)}
\cr
&+ |\zeta|  \| {\boldsymbol \phi}_\eps\|_{L_2(\O)}
\| \w_\eps - {\boldsymbol \phi}_\eps\|_{L_2(\O)},
\end{aligned}
\eqno(5.5)
$$
where ${\mathcal C}_{12}= \|g\|_{L_\infty}^{1/2} (d \alpha_1)^{1/2}$.
If ${\rm Re}\,\zeta \geq 0$ (and then ${\rm Im}\,\zeta \ne 0$), this yields
$$
\begin{aligned}
 \| \w_\eps - {\boldsymbol \phi}_\eps\|^2_{L_2(\O)}
 &\leq 2 {\mathcal C}_{12} c(\varphi) |\zeta|^{-1} \| \D {\boldsymbol \phi}_\eps\|_{L_2(\O)}
\|(g^\eps)^{1/2} b(\D)(\w_\eps - {\boldsymbol \phi}_\eps)\|_{L_2(\O)}
\cr
&+c(\varphi)^2 \| {\boldsymbol \phi}_\eps\|_{L_2(\O)}^2.
\end{aligned}
 \eqno(5.6)
$$
If ${\rm Re}\,\zeta < 0$, we take the real part of the corresponding relation and obtain
$$
\begin{aligned}
|{\rm Re}\,\zeta| \|\w_\eps - {\boldsymbol \phi}_\eps \|^2_{L_2(\O)}
&\leq
{\mathcal C}_{12} \| \D {\boldsymbol \phi}_\eps\|_{L_2(\O)}
\|(g^\eps)^{1/2} b(\D)(\w_\eps - {\boldsymbol \phi}_\eps)\|_{L_2(\O)}
\cr
&+ |\zeta|  \| {\boldsymbol \phi}_\eps\|_{L_2(\O)}
\| \w_\eps - {\boldsymbol \phi}_\eps\|_{L_2(\O)}.
\end{aligned}
\eqno(5.7)
$$
Summing up (5.5) and (5.7), we deduce the inequality similar to (5.6).
As a result, for all values of $\zeta$ under consideration we obtain
$$
\begin{aligned}
 \| \w_\eps - {\boldsymbol \phi}_\eps\|^2_{L_2(\O)}
 &\leq 4 {\mathcal C}_{12} c(\varphi) |\zeta|^{-1} \| \D {\boldsymbol \phi}_\eps\|_{L_2(\O)}
\|(g^\eps)^{1/2} b(\D)(\w_\eps - {\boldsymbol \phi}_\eps)\|_{L_2(\O)}
\cr
&+4 c(\varphi)^2 \| {\boldsymbol \phi}_\eps\|_{L_2(\O)}^2.
\end{aligned}
 \eqno(5.8)
$$

Now, from (5.4) with ${\boldsymbol \eta} = \w_\eps - {\boldsymbol \phi}_\eps$ and (5.8) it follows that
$$
\begin{aligned}
&a_{D,\eps}[\w_\eps - {\boldsymbol \phi}_\eps, \w_\eps - {\boldsymbol \phi}_\eps]
\cr
&\leq
9  c(\varphi)^2 |\zeta| \| {\boldsymbol \phi}_\eps\|_{L_2(\O)}^2
+ 9 {\mathcal C}_{12} c(\varphi) \| \D {\boldsymbol \phi}_\eps\|_{L_2(\O)}
\|(g^\eps)^{1/2} b(\D)(\w_\eps - {\boldsymbol \phi}_\eps)\|_{L_2(\O)}.
\end{aligned}
$$
This implies the inequality
$$
a_{D,\eps}[\w_\eps - {\boldsymbol \phi}_\eps, \w_\eps - {\boldsymbol \phi}_\eps]
\leq 18 c(\varphi)^2 |\zeta| \| {\boldsymbol \phi}_\eps\|_{L_2(\O)}^2
+ 81 {\mathcal C}_{12}^2 c(\varphi)^2 \| \D {\boldsymbol \phi}_\eps\|^2_{L_2(\O)}.
\eqno(5.9)
$$

By (5.9) and (3.1),
$$
\| \D (\w_\eps - {\boldsymbol \phi}_\eps) \|_{L_2(\O)}
\leq c(\varphi) \left( \check{\mathcal C}_{10} |\zeta|^{1/2} \| {\boldsymbol \phi}_\eps\|_{L_2(\O)}
+  \check{\mathcal C}_{11}  \| \D {\boldsymbol \phi}_\eps\|_{L_2(\O)}\right),
\eqno(5.10)
$$
where $\check{\mathcal C}_{10} = \sqrt{18} c_0^{-1/2}$, $\check{\mathcal C}_{11} =9 c_0^{-1/2} {\mathcal C}_{12}$.
Next, by (5.8) and (5.9),
$$
 \| \w_\eps - {\boldsymbol \phi}_\eps\|_{L_2(\O)}
\leq c(\varphi) \left( \sqrt{22} \| {\boldsymbol \phi}_\eps\|_{L_2(\O)}
+  \sqrt{85} \,{\mathcal C}_{12}  |\zeta|^{-1/2}\| \D {\boldsymbol \phi}_\eps\|_{L_2(\O)}\right).
\eqno(5.11)
$$
Relations (5.10), (5.11) together with the condition $|\zeta| \geq 1$ imply (5.3) with
${\mathcal C}_{10}= \check{\mathcal C}_{10}+ \sqrt{22} + 1$ and
${\mathcal C}_{11}= \check{\mathcal C}_{11}+ \sqrt{85} {\mathcal C}_{12}+1$. $\bullet$

\smallskip\noindent\textbf{5.2. Estimate for ${\boldsymbol \phi}_\eps$.}

\smallskip\noindent\textbf{Lemma 5.2.} \textit{Suppose that $\eps_1$ satisfies Condition} 4.1.
\textit{Let ${\boldsymbol \phi}_\eps$ be given by} (5.2). \textit{For $0< \eps \leq \eps_1$ and $\zeta \in \C \setminus \R_+$, $|\zeta|\geq 1$,
we have}
$$
\| {\boldsymbol \phi}_\eps \|_{L_2({\mathcal O})} \leq {\mathcal C}_{13} c(\varphi)
|\zeta|^{-1/2}\eps \| \FF \|_{L_2({\mathcal O})},
 \eqno(5.12)
$$
$$
\| \D {\boldsymbol \phi}_\eps \|_{L_2({\mathcal O})} \leq c(\varphi)
\left({\mathcal C}_{14} |\zeta|^{-1/4}\eps^{1/2} + {\mathcal C}_{15} \eps \right)
 \| \FF \|_{L_2({\mathcal O})}.
 \eqno(5.13)
$$
\textit{The constants ${\mathcal C}_{13}$, ${\mathcal C}_{14}$, and ${\mathcal C}_{15}$
depend on} $d$, $m$, $\alpha_0$, $\alpha_1$, $\| g \|_{L_\infty}$, $\| g^{-1} \|_{L_\infty}$,
\textit{the parameters of the lattice} $\Gamma$, \textit{and the domain} $\O$.

\smallskip\noindent\textbf{Proof.} Assume that $0< \eps \leq \eps_1$.
We start with the estimate for the $L_2$-norm of ${\boldsymbol \phi}_\eps$.
Relations (1.4), (1.19), (4.12), and (5.1) imply (5.12)
with ${\mathcal C}_{13}=M_1 \alpha_1^{1/2} C_\O^{(1)}({\mathcal C}_0+1)$.

Now we consider the derivatives
$$
\begin{aligned}
\partial_j {\boldsymbol \phi}_\eps =& \eps (\partial_j \theta_\eps) \Lambda^\eps S_\eps b(\D) \widetilde{\u}_0
+ \theta_\eps (\partial_j \Lambda)^\eps S_\eps b(\D) \widetilde{\u}_0
\cr
&+ \eps \theta_\eps \Lambda^\eps S_\eps b(\D) \partial_j \widetilde{\u}_0,\ \ j=1,\dots,d.
\end{aligned}
$$
Hence,
$$
\begin{aligned}
&\| \D {\boldsymbol \phi}_\eps \|^2_{L_2(\O)}
\leq 3 \eps^2 \| (\nabla \theta_\eps) \Lambda^\eps S_\eps b(\D) \widetilde{\u}_0 \|^2_{L_2(\R^d)}+
\cr
&+ 3  \| \theta_\eps (\D\Lambda)^\eps S_\eps b(\D) \widetilde{\u}_0 \|^2_{L_2(\R^d)}
+ 3 \eps^2 \sum_{j=1}^d \| \theta_\eps \Lambda^\eps S_\eps b(\D) \partial_j \widetilde{\u}_0 \|^2_{L_2(\R^d)}.
\end{aligned}
\eqno(5.14)
$$
Denote the consecutive terms in the right-hand side of  (5.14) by $J_1(\eps)$, $J_2(\eps)$, $J_3(\eps)$.

By (5.1) and Lemma 3.4,
$$
\begin{aligned}
&J_1(\eps) \leq 3 \kappa^2 \int_{(\partial \O)_\eps} |\Lambda^\eps S_\eps b(\D) \widetilde{\u}_0|^2 \,d\x
\cr
&\leq 3 \kappa^2 \beta_* \eps |\Omega|^{-1} \|\Lambda\|_{L_2(\Omega)}^2
\| b(\D) \widetilde{\u}_0\|_{H^1(\R^d)} \| b(\D) \widetilde{\u}_0\|_{L_2(\R^d)}.
\end{aligned}
$$
Combining this with (1.4), (1.11), (4.12), and (4.13), we arrive at
$$
J_1(\eps) \leq  \gamma_1 c(\varphi)^2 |\zeta|^{-1/2} \eps \|\FF\|^2_{L_2(\O)},
\eqno(5.15)
$$
where $\gamma_1 = 3 \kappa^2 \beta_* M_1^2 \alpha_1 C_\O^{(2)} \widehat{c} C_\O^{(1)}({\mathcal C}_0+1)$.
Similarly, for the second term in (5.14) we have
$$
\begin{aligned}
&J_2(\eps) \leq 3 \int_{(\partial \O)_\eps} |(\D\Lambda)^\eps S_\eps b(\D) \widetilde{\u}_0|^2 \,d\x
\cr
&\leq 3 \beta_* \eps |\Omega|^{-1} \|\D\Lambda\|_{L_2(\Omega)}^2
\| b(\D) \widetilde{\u}_0\|_{H^1(\R^d)} \| b(\D) \widetilde{\u}_0\|_{L_2(\R^d)}.
\end{aligned}
$$
Using (1.4), (1.10), (4.12), and (4.13), from here we deduce that
$$
J_2(\eps) \leq  \gamma_2 c(\varphi)^2 |\zeta|^{-1/2} \eps \|\FF\|^2_{L_2(\O)},
\eqno(5.16)
$$
where $\gamma_2 = 3 \beta_* M_2^2 \alpha_1 C_\O^{(2)} \widehat{c} C_\O^{(1)}({\mathcal C}_0+1)$.
The third term in (5.14) is estimated with the help of (1.4), (1.19), (4.13), and (5.1):
$$
J_3(\eps) \leq  \gamma_3 c(\varphi)^2 \eps^2 \|\FF\|^2_{L_2(\O)},
\eqno(5.17)
$$
where $\gamma_3 = 3 M_1^2 \alpha_1 (C_\O^{(2)} \widehat{c})^2$.
Now, relations (5.14)--(5.17) imply (5.13) with
${\mathcal C}_{14}= (\gamma_1 + \gamma_2)^{1/2}$ and ${\mathcal C}_{15}= \gamma_3^{1/2}$. $\bullet$

\smallskip\noindent\textbf{5.3. Completion of the proof of Theorem 4.3.}
By Lemmas 5.1 and 5.2,
$$
\|\w_\eps\|_{H^1(\O)} \leq c(\varphi)^2 \left( {\mathcal C}_{16} |\zeta|^{-1/4} \eps^{1/2} + {\mathcal C}_{17} \eps \right)\|\FF\|_{L_2(\O)},
\ \ 0< \eps \leq \eps_1,
$$
where ${\mathcal C}_{16} =  {\mathcal C}_{11} {\mathcal C}_{14}$ and
${\mathcal C}_{17} =  {\mathcal C}_{10} {\mathcal C}_{13} + {\mathcal C}_{11} {\mathcal C}_{15}$.
Together with (4.31) and (4.32) this implies the required estimate (4.8) with
${\mathcal C}_2 = \sqrt{2}{\mathcal C}_{16}$ and
${\mathcal C}_3 = {\mathcal C}_5 + {\mathcal C}_6+ \sqrt{2}{\mathcal C}_{17}$.

It remains to check (4.10). From (4.8) and (1.2), (1.5) it follows that
$$
\begin{aligned}
&\| \p_\eps - g^\eps b(\D) \v_\eps \|_{L_2(\O)}
\cr
&\leq
\|g\|_{L_\infty} (d \alpha_1)^{1/2}
\left( {\mathcal C}_{2} c(\varphi)^2 |\zeta|^{-1/4} \eps^{1/2} + {\mathcal C}_{3} c(\varphi)^4 \eps \right)\|\FF\|_{L_2(\O)}.
\end{aligned}
\eqno(5.18)
$$
We have
$$
\begin{aligned}
g^\eps b(\D) \v_\eps &= g^\eps b(\D) \u_0 + g^\eps (b(\D) \Lambda)^\eps S_\eps b(\D) \widetilde{\u}_0
\cr
&+
\eps \sum_{l=1}^d g^\eps b_l \Lambda^\eps S_\eps b(\D) D_l \widetilde{\u}_0.
\end{aligned}
\eqno(5.19)
$$
The third term in (5.19) is estimated with the help of (1.4), (1.5), and (1.19):
$$
\left\| \eps \sum_{l=1}^d g^\eps b_l \Lambda^\eps S_\eps b(\D) D_l \widetilde{\u}_0 \right\|_{L_2(\O)} \leq
{\mathcal C}' \eps \| \widetilde{\u}_0\|_{H^2(\R^d)},
\eqno(5.20)
$$
where ${\mathcal C}' = \|g\|_{L_\infty} \alpha_1 d^{1/2} M_1$.
Note that, by Proposition 1.4 and (1.4), we have
$$
\| g^\eps (I- S_\eps) b(\D) \widetilde{\u}_0 \|_{L_2(\O)}
\leq {\mathcal C}''  \eps \| \widetilde{\u}_0\|_{H^2(\R^d)},
\eqno(5.21)
$$
where ${\mathcal C}'' = \|g\|_{L_\infty} r_1 \alpha_1^{1/2}$.
From (5.19)--(5.21) and (1.9), (4.13) it follows that
$$
\|g^\eps b(\D) \v_\eps - \widetilde{g}^\eps S_\eps b(\D) \widetilde{\u}_0\|_{L_2(\O)}
\leq {\mathcal C}_{18} c(\varphi) \eps \|\FF\|_{L_2(\O)},
\eqno(5.22)
$$
where ${\mathcal C}_{18} =  ({\mathcal C}' + {\mathcal C}'') C_\O^{(2)} \widehat{c}$.

Now,  (5.18) and (5.22) imply (4.10) with
$\widetilde{\mathcal C}_2= \|g\|_{L_\infty} (d\alpha_1)^{1/2}{\mathcal C}_2$ and
$\widetilde{\mathcal C}_3=\|g\|_{L_\infty} (d\alpha_1)^{1/2}{\mathcal C}_3 + {\mathcal C}_{18}$.
$\bullet$

\smallskip\noindent\textbf{5.4. Proof of Theorem 4.2.}

\smallskip\noindent\textbf{Lemma 5.3.} \textit{Let $\w_\eps$ be the solution of problem} (4.20).
\textit{Suppose that the number $\eps_1$ satisfies Condition} 4.1.
\textit{Then for $0< \eps \leq \eps_1$ and $\zeta \in \C \setminus \R_+$, $|\zeta|\geq 1$, we have}
$$
\| \w_\eps \|_{L_2({\mathcal O})} \leq
c(\varphi)^5 \left( {\mathcal C}_{19}|\zeta|^{-1/2} \eps + {\mathcal C}_{20} \eps^2 \right) \|\FF\|_{L_2(\O)}.
 \eqno(5.23)
$$
\textit{The constants ${\mathcal C}_{19}$ and ${\mathcal C}_{20}$
depend on} $d$, $m$, $\alpha_0$, $\alpha_1$, $\| g \|_{L_\infty}$, $\| g^{-1} \|_{L_\infty}$,
\textit{the parameters of the lattice} $\Gamma$, \textit{and the domain} $\O$.

\smallskip\noindent\textbf{Proof.} Assume that $0< \eps \leq \eps_1$.
Consider the identity (5.4) and take ${\boldsymbol \eta}$ equal to
${\boldsymbol \eta}_\eps = (\A_{D,\eps} - \overline{\zeta} I)^{-1}{\boldsymbol \Phi}$, where ${\boldsymbol \Phi}\in L_2(\O;\C^n)$.
 Then the left-hand side of (5.4) takes the form
 $$
 (g^\eps b(\D)(\w_\eps - {\boldsymbol \phi}_\eps), b(\D){\boldsymbol \eta}_\eps)_{L_2(\O)}- \zeta (\w_\eps - {\boldsymbol \phi}_\eps,
 {\boldsymbol \eta}_\eps)_{L_2(\O)} = (\w_\eps - {\boldsymbol \phi}_\eps, {\boldsymbol \Phi})_{L_2(\O)}.
 $$
 Hence,
 $$
 (\w_\eps - {\boldsymbol \phi}_\eps, {\boldsymbol \Phi})_{L_2(\O)} =
 -(g^\eps b(\D) {\boldsymbol \phi}_\eps, b(\D){\boldsymbol \eta}_\eps)_{L_2(\O)} +  \zeta ({\boldsymbol \phi}_\eps,
 {\boldsymbol \eta}_\eps)_{L_2(\O)}.
\eqno(5.24)
 $$
 To approximate  ${\boldsymbol \eta}_\eps$ in $H^1(\O;\C^n)$, we apply the already proved Theorem 4.3.
We put ${\boldsymbol \eta}_0= (\A_D^0 - \overline{\zeta} I)^{-1}{\boldsymbol \Phi}$, $\widetilde{\boldsymbol \eta}_0 = P_\O {\boldsymbol \eta}_0$. Then approximation of ${\boldsymbol \eta}_\eps$ is given by
 ${\boldsymbol \eta}_0 + \eps \Lambda^\eps S_\eps b(\D) \widetilde{\boldsymbol \eta}_0$. Rewrite (5.24) as
 $$
 \begin{aligned}
 &(\w_\eps - {\boldsymbol \phi}_\eps, {\boldsymbol \Phi})_{L_2(\O)} =
 -\left(g^\eps b(\D) {\boldsymbol \phi}_\eps, b(\D)({\boldsymbol \eta}_\eps - {\boldsymbol \eta}_0 - \eps \Lambda^\eps S_\eps b(\D) \widetilde{\boldsymbol \eta}_0)\right)_{L_2(\O)}
 \cr
 &-\left(g^\eps b(\D) {\boldsymbol \phi}_\eps, b(\D){\boldsymbol \eta}_0 \right)_{L_2(\O)}
-\left(g^\eps b(\D) {\boldsymbol \phi}_\eps, b(\D)(\eps \Lambda^\eps S_\eps b(\D) \widetilde{\boldsymbol \eta}_0)  \right)_{L_2(\O)}
 \cr
 &+  \zeta ({\boldsymbol \phi}_\eps,{\boldsymbol \eta}_\eps)_{L_2(\O)}.
\end{aligned}
\eqno(5.25)
 $$
 Denote the consecutive terms in the right-hand side of (5.25) by ${\mathcal I}_j(\eps)$, $j=1,2,3,4$.

It is easy to estimate the fourth term in (5.25), using Lemma 3.1 and (5.12):
$$
|{\mathcal I}_4(\eps)| \leq |\zeta| \|{\boldsymbol \phi}_\eps\|_{L_2(\O)}
\|{\boldsymbol \eta}_\eps\|_{L_2(\O)} \leq {\mathcal C}_{13} c(\varphi)^2 |\zeta|^{-1/2} \eps
\|\FF\|_{L_2(\O)} \|{\boldsymbol \Phi}\|_{L_2(\O)}.
\eqno(5.26)
$$

Now, let us estimate the first term in (5.25). By (1.2) and (1.5), we have
$$
|{\mathcal I}_1(\eps)| \leq \|g\|_{L_\infty} d \alpha_1  \|\D{\boldsymbol \phi}_\eps\|_{L_2(\O)}
\|\D({\boldsymbol \eta}_\eps  - {\boldsymbol \eta}_0 - \eps \Lambda^\eps S_\eps b(\D) \widetilde{\boldsymbol \eta}_0) \|_{L_2(\O)}.
$$
Applying Theorem 4.3 and (5.13), we arrive at
$$
\begin{aligned}
|{\mathcal I}_1(\eps)| &\leq \|g\|_{L_\infty} d \alpha_1  c(\varphi)
\left({\mathcal C}_{14} |\zeta|^{-1/4} \eps^{1/2} + {\mathcal C}_{15} \eps \right)
\cr
&\times \left({\mathcal C}_{2} c(\varphi)^2|\zeta|^{-1/4} \eps^{1/2} + {\mathcal C}_{3} c(\varphi)^4 \eps \right)
\|\FF\|_{L_2(\O)} \|{\boldsymbol \Phi}\|_{L_2(\O)}.
\end{aligned}
$$
Hence,
$$
|{\mathcal I}_1(\eps)| \leq  c(\varphi)^5 \left(\check{\gamma}_1 |\zeta|^{-1/2} \eps + \check{\gamma}_2 \eps^2 \right)
\|\FF\|_{L_2(\O)} \|{\boldsymbol \Phi}\|_{L_2(\O)},
\eqno(5.27)
$$
where $\check{\gamma}_1 = \|g\|_{L_\infty} d \alpha_1\left({\mathcal C}_{2}{\mathcal C}_{14} +
{\mathcal C}_{2}{\mathcal C}_{15} + {\mathcal C}_{3} {\mathcal C}_{14}\right)$
and $\check{\gamma}_2 = \|g\|_{L_\infty} d \alpha_1
\left({\mathcal C}_{3}{\mathcal C}_{14} + {\mathcal C}_{3}{\mathcal C}_{15} + {\mathcal C}_{2} {\mathcal C}_{15}\right)$.

In order to estimate the second term in  (5.25), recall that ${\boldsymbol \phi}_\eps$ is supported in the $\eps$-neighborhood of the boundary and apply
Lemma 3.3($1^\circ$). Taking (1.2) and (1.5) into account, we have
$$
\begin{aligned}
&|{\mathcal I}_2(\eps)| \leq \|g\|_{L_\infty} (d \alpha_1)^{1/2}  \|\D{\boldsymbol \phi}_\eps\|_{L_2(\O)}
\left(\int_{B_\eps} |b(\D){\boldsymbol \eta}_0|^2\,d\x \right)^{1/2}
\cr
&\leq
\|g\|_{L_\infty} (d \alpha_1)^{1/2}  \|\D{\boldsymbol \phi}_\eps\|_{L_2(\O)}
\beta^{1/2} \eps^{1/2} \|b(\D){\boldsymbol \eta}_0\|_{H^1(\O)}^{1/2}
\|b(\D){\boldsymbol \eta}_0\|_{L_2(\O)}^{1/2}.
\end{aligned}
$$
Together with Lemma 3.2 and (5.13) this yields
$$
\begin{aligned}
|{\mathcal I}_2(\eps)| &\leq \|g\|_{L_\infty} d \alpha_1 c(\varphi)
\left({\mathcal C}_{14} |\zeta|^{-1/4} \eps^{1/2} + {\mathcal C}_{15} \eps \right)
\cr
&\times \beta^{1/2} \eps^{1/2} ( \widehat{c} c(\varphi))^{1/2}
({\mathcal C}_0 c(\varphi) |\zeta|^{-1/2})^{1/2}
\|\FF\|_{L_2(\O)} \|{\boldsymbol \Phi}\|_{L_2(\O)}.
\end{aligned}
$$
Hence,
$$
|{\mathcal I}_2(\eps)| \leq  c(\varphi)^2 \left(\check{\gamma}_3 |\zeta|^{-1/2} \eps + \check{\gamma}_4 \eps^2 \right)
\|\FF\|_{L_2(\O)} \|{\boldsymbol \Phi}\|_{L_2(\O)},
\eqno(5.28)
$$
where $\check{\gamma}_3 = \|g\|_{L_\infty} d \alpha_1 (\beta\widehat{c}\, {\mathcal C}_0)^{1/2}
({\mathcal C}_{14} +{\mathcal C}_{15})$ and
$\check{\gamma}_4 = \|g\|_{L_\infty} d \alpha_1 (\beta \widehat{c} \,{\mathcal C}_0)^{1/2} {\mathcal C}_{15}$.

It remains to estimate the third term in (5.25). By (1.2), we have
$$
{\mathcal I}_3(\eps) =
{\mathcal I}_3^{(1)}(\eps)+ {\mathcal I}_3^{(2)}(\eps),
\eqno(5.29)
$$
$$
{\mathcal I}_3^{(1)}(\eps)
= -\left(g^\eps b(\D) {\boldsymbol \phi}_\eps, (b(\D) \Lambda)^\eps S_\eps b(\D) \widetilde{\boldsymbol \eta}_0  \right)_{L_2(\O)},
\eqno(5.30)
$$
$$
{\mathcal I}_3^{(2)}(\eps)
= -\bigl(g^\eps b(\D) {\boldsymbol \phi}_\eps, \eps \sum_{l=1}^d b_l  \Lambda^\eps S_\eps b(\D) D_l \widetilde{\boldsymbol \eta}_0  \bigr)_{L_2(\O)}.
\eqno(5.31)
$$
The term (5.30) can be estimated with the help of (1.2), (1.5), and Lemma 3.4:
$$
\begin{aligned}
&| {\mathcal I}_3^{(1)}(\eps)| \leq
\|g\|_{L_\infty} (d \alpha_1)^{1/2} \|\D{\boldsymbol \phi}_\eps\|_{L_2(\O)}
\left( \int_{(\partial \O)_\eps} |(b(\D) \Lambda)^\eps S_\eps b(\D) \widetilde{\boldsymbol \eta}_0|^2\,d\x \right)^{1/2}
\cr
&\leq
\|g\|_{L_\infty} (d \alpha_1)^{1/2} \|\D{\boldsymbol \phi}_\eps\|_{L_2(\O)}
\cr
&\times
\left( \beta_* \eps |\Omega|^{-1} \| b(\D) \Lambda\|^2_{L_2(\Omega)} \|b(\D) \widetilde{\boldsymbol \eta}_0\|_{H^1(\R^d)}
\|b(\D) \widetilde{\boldsymbol \eta}_0\|_{L_2(\R^d)} \right)^{1/2}.
\end{aligned}
\eqno(5.32)
$$
Similarly to (4.12) and (4.13), we have
$$
\| \widetilde{\boldsymbol \eta}_0\|_{H^1(\R^d)}
\leq {C}_\O^{(1)} ({\mathcal C}_0 +1) c(\varphi) |\zeta|^{-1/2} \| {\boldsymbol \Phi}\|_{L_2(\O)},
\eqno(5.33)
$$
$$
\| \widetilde{\boldsymbol \eta}_0\|_{H^2(\R^d)}
\leq  {C}_\O^{(2)} \widehat{c} c(\varphi)  \| {\boldsymbol \Phi}\|_{L_2(\O)}.
\eqno(5.34)
$$
From (5.32)--(5.34) and (1.10), (5.13) it follows that
$$
\begin{aligned}
&| {\mathcal I}_3^{(1)}(\eps)| \leq
\|g\|_{L_\infty} d \alpha_1^{3/2} M_2
c(\varphi) \left({\mathcal C}_{14} |\zeta|^{-1/4} \eps^{1/2} + {\mathcal C}_{15} \eps \right)
\cr
&\times
\beta_*^{1/2} \eps^{1/2} ( {C}_\O^{(2)} \widehat{c} c(\varphi))^{1/2}
({C}_\O^{(1)} ({\mathcal C}_0 +1) c(\varphi) |\zeta|^{-1/2})^{1/2}
\|\FF\|_{L_2(\O)} \|{\boldsymbol \Phi}\|_{L_2(\O)}.
\end{aligned}
$$
Hence,
$$
| {\mathcal I}_3^{(1)}(\eps)| \leq
c(\varphi)^2 \left(\check{\gamma}_5 |\zeta|^{-1/2} \eps + \check{\gamma}_6 \eps^2 \right)
\|\FF\|_{L_2(\O)} \|{\boldsymbol \Phi}\|_{L_2(\O)},
\eqno(5.35)
$$
where $\check{\gamma}_5 = \|g\|_{L_\infty} d \alpha_1^{3/2} M_2 ( \beta_* {C}_\O^{(2)} \widehat{c}
\,{C}_\O^{(1)} ({\mathcal C}_0 +1))^{1/2} ({\mathcal C}_{14}+{\mathcal C}_{15})$ and
$\check{\gamma}_6 = \|g\|_{L_\infty} d \alpha_1^{3/2} M_2 (\beta_* {C}_\O^{(2)} \widehat{c}\,
{C}_\O^{(1)} ({\mathcal C}_0 +1))^{1/2} {\mathcal C}_{15}$.

Finally, the term (5.31) is estimated by using  (1.2), (1.4), (1.5), and (1.19):
$$
| {\mathcal I}_3^{(2)}(\eps)| \leq
\eps \|g\|_{L_\infty} d \alpha_1^{3/2} M_1 \|\D{\boldsymbol \phi}_\eps\|_{L_2(\O)}
\|\widetilde{\boldsymbol \eta}_0\|_{H^2(\R^d)}.
$$
Combining this with (5.13) and (5.34), we obtain
$$
| {\mathcal I}_3^{(2)}(\eps)| \leq
c(\varphi)^2 \left(\check{\gamma}_7 |\zeta|^{-1/2} \eps + \check{\gamma}_8 \eps^2 \right)
\|\FF\|_{L_2(\O)} \|{\boldsymbol \Phi}\|_{L_2(\O)},
\eqno(5.36)
$$
where $\check{\gamma}_7 = \|g\|_{L_\infty} d \alpha_1^{3/2} M_1 {C}_\O^{(2)} \widehat{c}\,{\mathcal C}_{14}$,
$\check{\gamma}_8 = \|g\|_{L_\infty} d \alpha_1^{3/2} M_1 {C}_\O^{(2)} \widehat{c} ({\mathcal C}_{14} + {\mathcal C}_{15})$.

As a result, relations (5.25)--(5.29), (5.35), and (5.36) imply that
$$
|(\w_\eps - {\boldsymbol \phi}_\eps, {\boldsymbol \Phi})_{L_2(\O)} | \leq
c(\varphi)^5 \left(\widetilde{\gamma} |\zeta|^{-1/2} \eps + \widehat{\gamma} \eps^2 \right)
\|\FF\|_{L_2(\O)} \|{\boldsymbol \Phi}\|_{L_2(\O)}
$$
for any ${\boldsymbol \Phi}\in L_2(\O;\C^n)$,
where $\widetilde{\gamma}= {\mathcal C}_{13} + \check{\gamma}_1 + \check{\gamma}_3 + \check{\gamma}_5 + \check{\gamma}_7$,
$\widehat{\gamma}= \check{\gamma}_2 + \check{\gamma}_4 + \check{\gamma}_6 + \check{\gamma}_8$.
Consequently,
$$
\|\w_\eps - {\boldsymbol \phi}_\eps\|_{L_2(\O)}  \leq
c(\varphi)^5 \left(\widetilde{\gamma} |\zeta|^{-1/2} \eps + \widehat{\gamma} \eps^2 \right)
\|\FF\|_{L_2(\O)}.
\eqno(5.37)
$$
Now, the required estimate (5.23) follows directly from (5.37) and (5.12); herewith,
${\mathcal C}_{19}= \widetilde{\gamma} + {\mathcal C}_{13}$ and  ${\mathcal C}_{20}= \widehat{\gamma}$.
$\bullet$

\smallskip\noindent\textbf{Completion of the proof of Theorem 4.2.}
By (4.34) and (5.23), for $0< \eps \leq \eps_1$ we have
$$
\begin{aligned}
&\| \u_\eps - \u_0\|_{L_2(\O)} \leq
({\mathcal C}_{6}+{\mathcal C}_{9}) c(\varphi)^4 |\zeta|^{-1/2} \eps \|\FF\|_{L_2(\O)}
\cr
&+
c(\varphi)^5 \left( {\mathcal C}_{19}|\zeta|^{-1/2} \eps + {\mathcal C}_{20} \eps^2 \right)
\|\FF\|_{L_2(\O)}.
\end{aligned}
$$
This implies the required estimate (4.1) with
${\mathcal C}_1= \max\{ {\mathcal C}_{6}+{\mathcal C}_{9}+{\mathcal C}_{19}; {\mathcal C}_{20}\}$.
$\bullet$

\section*{\S 6. The results for the Dirichlet problem: \\ the case where $\Lambda \in L_\infty$, special cases}

\smallskip\noindent\textbf{6.1. The case where $\Lambda \in L_\infty$.}
As for the problem in $\R^d$ (see Subsection 2.3), under the additional Condition 2.8
it is possible to remove the smoothing operator in the corrector,
i.~e., instead of the corrector (4.4) one can use the following simpler operator
$$
K^0_D(\eps;\zeta) = [\Lambda^\eps] b(\D) (\A^0_D - \zeta I)^{-1}.
\eqno(6.1)
$$
Using Proposition 2.9, it is easy to check that under Condition 2.8 the operator (6.1)
is a continuous mapping of $L_2(\O;\C^n)$ to $H^1(\O;\C^n)$.
Instead of (4.7), now we use another approximation of the solution $\u_\eps$ of problem (3.3):
$$
\check{\v}_\eps := (\A^0_D - \zeta I)^{-1}\FF +\eps K^0_D(\eps;\zeta)\FF = \u_0 + \eps \Lambda^\eps b(\D)\u_0.
\eqno(6.2)
$$

\smallskip\noindent\textbf{Theorem 6.1.} \textit{Suppose that the assumptions of Theorem} 4.2 \textit{and Condition} 2.8
 \textit{are satisfied. Let $\check{\v}_\eps$ be defined by} (6.2).
\textit{Then for $0< \varepsilon \leq \varepsilon_1$ we have}
$$
\| \u_\eps - \check{\v}_\eps \|_{H^1(\O)} \leq \left({\mathcal C}_2 c(\varphi)^2 |\zeta|^{-1/4} \varepsilon^{1/2} +
{\mathcal C}_3^\circ c(\varphi)^4 \varepsilon\right) \|\FF\|_{L_2(\O)},
\eqno(6.3)
$$
\textit{or, in operator terms},
$$
\begin{aligned}
&\| (\A_{D,\varepsilon}- \zeta I)^{-1}  - (\A^0_{D}- \zeta I)^{-1} - \eps K^0_D(\eps;\zeta) \|_{L_2({\mathcal O}) \to H^1({\mathcal O})}
\cr
&\leq {\mathcal C}_{2} c(\varphi)^2  |\zeta|^{-1/4}\varepsilon^{1/2}
+ {\mathcal C}_3^\circ c(\varphi)^4 \varepsilon.
\end{aligned}
$$
\textit{For the flux} $\p_\eps := g^\eps b(\D) \u_\eps$ \textit{we have}
$$
\| \p_\eps - \widetilde{g}^\eps b(\D) {\u}_0 \|_{L_2(\O)}
\leq \left(\widetilde{\mathcal C}_2 c(\varphi)^2 |\zeta|^{-1/4}\varepsilon^{1/2} +
\widetilde{\mathcal C}_3^\circ c(\varphi)^4 \varepsilon\right) \|\FF\|_{L_2(\O)}
\eqno(6.4)
$$
\textit{for $0< \eps \leq \eps_1$. The constants ${\mathcal C}_2$ and $\widetilde{\mathcal C}_2$ are the same as in Theorem} 4.3.
\textit{The constants ${\mathcal C}_3^\circ$ and $\widetilde{\mathcal C}_3^\circ$ depend on}
$d$, $m$, $\alpha_0$, $\alpha_1$, $\| g \|_{L_\infty}$, $\| g^{-1} \|_{L_\infty}$,
\textit{the parameters of the lattice} $\Gamma$, \textit{the domain $\O$, and the norm} $\|\Lambda\|_{L_\infty}$.

\smallskip\noindent\textbf{Proof.}
In order to deduce (6.3) from (4.8), we need to estimate the $H^1(\O;\C^n)$-norm of the function
$\v_\eps - \check{\v}_\eps = \eps  \left(\Lambda^\eps(S_\eps-I) b(\D)\widetilde{\u}_0\right)\vert_\O.$
We start with the $L_2$-norm. By Condition 2.8 and estimates (1.4), (1.18), we have
$$
\|\v_\eps - \check{\v}_\eps\|_{L_2(\O)} \leq \eps \|\Lambda\|_{L_\infty} \|(S_\eps -I)b(\D)\widetilde{\u}_0\|_{L_2(\R^d)}
\leq 2 \|\Lambda\|_{L_\infty} \alpha_1^{1/2}\eps \| \widetilde{\u}_0 \|_{H^1(\R^d)}.
\eqno(6.5)
$$
Now, consider the derivatives
$$
\partial_j(\v_\eps - \check{\v}_\eps) = \left((\partial_j \Lambda)^\eps (S_\eps-I) b(\D)\widetilde{\u}_0\right)\vert_\O
+ \eps \left( \Lambda^\eps (S_\eps-I) b(\D)\partial_j \widetilde{\u}_0\right)\vert_\O,
$$
$j=1,\dots,d$. Hence,
$$
\begin{aligned}
\|\D (\v_\eps - \check{\v}_\eps)\|^2_{L_2(\O)}
&\leq 2 \|(\D \Lambda)^\eps (S_\eps-I) b(\D)\widetilde{\u}_0 \|^2_{L_2(\R^d)}
\cr
&+2 \eps^2 \sum_{j=1}^d \|\Lambda^\eps (S_\eps-I) b(\D)\partial_j \widetilde{\u}_0 \|^2_{L_2(\R^d)}.
\end{aligned}
\eqno(6.6)
$$
By Proposition 2.9, this yields
$$
\begin{aligned}
&\|\D (\v_\eps - \check{\v}_\eps)\|^2_{L_2(\O)}
\leq 2 \beta_1 \|(S_\eps-I) b(\D)\widetilde{\u}_0 \|^2_{L_2(\R^d)}
\cr
&+2 \eps^2 \|\Lambda\|^2_{L_\infty} (\beta_2 +1) \sum_{j=1}^d \| (S_\eps-I) b(\D)\partial_j \widetilde{\u}_0 \|^2_{L_2(\R^d)}.
\end{aligned}
\eqno(6.7)
$$
To estimate the first term, we apply Proposition 1.4. The second term is estimated with the help of
(1.18). Taking (1.4) into account, we arrive at
$$
\|\D (\v_\eps - \check{\v}_\eps)\|_{L_2(\O)} \leq {\mathcal C}_{21} \eps \| \widetilde{\u}_0 \|_{H^2(\R^d)},
\eqno(6.8)
$$
where ${\mathcal C}_{21} = \alpha_1^{1/2} (2 \beta_1 r_1^2 + 8 (\beta_2+1) \|\Lambda\|^2_{L_\infty})^{1/2}$.

As a result, from (6.5) and (6.8) it follows that
$$
\|\v_\eps - \check{\v}_\eps\|_{H^1(\O)} \leq {\mathcal C}'''  \eps \| \widetilde{\u}_0 \|_{H^2(\R^d)},
\eqno(6.9)
$$
where ${\mathcal C}'''=\alpha_1^{1/2} (2 \beta_1 r_1^2 + (8\beta_2+12) \|\Lambda\|^2_{L_\infty})^{1/2}$.
Now (4.8), (4.13), and (6.9) imply the required estimate (6.3) with
${\mathcal C}_{3}^\circ = {\mathcal C}_{3}+ {\mathcal C}''' C_{\O}^{(2)} \widehat{c}$.

It remains to check (6.4). From (6.3), (1.2), and (1.5) it follows that
$$
\begin{aligned}
&\| \p_\eps - g^\eps b(\D) \check{\v}_\eps \|_{L_2(\O)}
\cr
&\leq
\|g\|_{L_\infty} (d \alpha_1)^{1/2}
\left( {\mathcal C}_{2} c(\varphi)^2 |\zeta|^{-1/4} \eps^{1/2} + {\mathcal C}_{3}^\circ c(\varphi)^4 \eps \right)\|\FF\|_{L_2(\O)}
\end{aligned}
\eqno(6.10)
$$
for $0< \eps \leq \eps_1$. We have
$$
\begin{aligned}
g^\eps b(\D) \check{\v}_\eps &= g^\eps b(\D) \u_0 + g^\eps (b(\D) \Lambda)^\eps b(\D) {\u}_0
\cr
&+
\eps \sum_{l=1}^d g^\eps b_l \Lambda^\eps b(\D) D_l {\u}_0.
\end{aligned}
\eqno(6.11)
$$
The third term in (6.11) can be estimated by using (1.2) and (1.5):
$$
\left\| \eps \sum_{l=1}^d g^\eps b_l \Lambda^\eps b(\D) D_l {\u}_0 \right\|_{L_2(\O)} \leq
\widetilde{\mathcal C}' \eps \| \u_0 \|_{H^2(\O)},
\eqno(6.12)
$$
where $\widetilde{\mathcal C}' = \|g\|_{L_\infty} \alpha_1 d \|\Lambda\|_{L_\infty}$.
From (6.11), (6.12), and (1.9) it follows that
$$
\|g^\eps b(\D) \check{\v}_\eps - \widetilde{g}^\eps b(\D) {\u}_0\|_{L_2(\O)}
\leq \widetilde{\mathcal C}' \eps \|\u_0 \|_{H^2(\O)}.
\eqno(6.13)
$$

Relations (3.11), (6.10), and (6.13) imply (6.4) with
$\widetilde{\mathcal C}_3^\circ =\|g\|_{L_\infty} (d\alpha_1)^{1/2}{\mathcal C}_3^\circ + \widetilde{\mathcal C}' \widehat{c}$.
$\bullet$

\smallskip\noindent\textbf{6.2. The case of zero corrector.}
Next statement follows from Theorem 4.3, Proposition 1.2, and equation (1.7).
(Cf. Proposition 2.7.)

\smallskip\noindent\textbf{Proposition 6.2.} \textit{Suppose that the assumptions of Theorem} 4.3
 \textit{are satisfied. If $g^0 = \overline{g}$, i.~e., relations} (1.13) \textit{are satisfied, then
 $\Lambda=0$, $\v_\eps = \u_0$, and for $0< \varepsilon \leq \varepsilon_1$ we have}
$$
\| \u_\eps - {\u}_0 \|_{H^1(\O)} \leq \left({\mathcal C}_2 c(\varphi)^2 |\zeta|^{-1/4}\varepsilon^{1/2} +
{\mathcal C}_3  c(\varphi)^4 \varepsilon\right) \|\FF\|_{L_2(\O)}.
$$

\smallskip\noindent\textbf{6.3. The case where $g^0 = \underline{g}$.}
As has already been mentioned (see Subsection 2.3), under the condition $g^0 = \underline{g}$ the matrix (1.9)
is constant: $\widetilde{g}(\x)= g^0 = \underline{g}$. Besides, in this case Condition 2.8 is satisfied
(cf. Proposition 2.12). Applying the statement of Theorem 6.1 concerning the fluxes, we arrive at
the following result.

\smallskip\noindent\textbf{Proposition 6.3.} \textit{Suppose that the assumptions of Theorem} 4.3
 \textit{are satisfied. If $g^0 = \underline{g}$, i.~e., relations} (1.14) \textit{are satisfied, then
for $0< \varepsilon \leq \varepsilon_1$ we have}
$$
\| \p_\eps - g^0 b(\D){\u}_0 \|_{L_2(\O)} \leq \left(\widetilde{\mathcal C}_2 c(\varphi)^2 |\zeta|^{-1/4}\varepsilon^{1/2} +
\widetilde{\mathcal C}_3^\circ c(\varphi)^4 \varepsilon\right) \|\FF\|_{L_2(\O)}.
$$

\section*{\S 7. Approximation of solutions of the Dirichlet problem \\ in a strictly interior subdomain}

\smallskip\noindent\textbf{7.1. The general case.}
Using Theorem 4.2 and the results for homogenization problem in $\R^d$, it is not difficult to improve error estimates
in $H^1(\O';\C^n)$ for a strongly interior subdomain $\O'$ of $\O$.

\smallskip\noindent\textbf{Theorem 7.1.} \textit{Suppose that the assumptions of Theorem} 4.3
\textit{are satisfied. Let $\O'$ be a strictly interior subdomain of the domain $\O$. Denote}
$\delta:= {\rm dist}\,\{ \O'; \partial \O \}$. \textit{Then for $0< \eps \leq \eps_1$ we have}
$$
\|\u_\eps - \v_\eps \|_{H^1(\O')} \leq ({\mathcal C}'_{22} \delta^{-1} + {\mathcal C}''_{22}) c(\varphi)^6 \eps \| \FF \|_{L_2(\O)},
\eqno(7.1)
$$
\textit{or, in operator terms,}
$$
\begin{aligned}
&\| (\A_{D,\eps} - \zeta I)^{-1} - (\A_{D}^0 - \zeta I)^{-1} - \eps K_D(\eps;\zeta) \|_{L_2(\O) \to H^1(\O')}
\cr
&\leq ({\mathcal C}'_{22} \delta^{-1} + {\mathcal C}''_{22}) c(\varphi)^6 \eps.
\end{aligned}
$$
\textit{For the flux $\p_\eps = g^\eps b(\D) \u_\eps$ we have}
$$
\| \p_\eps - \widetilde{g}^\eps S_\eps b(\D) \widetilde{\u}_0
\|_{L_2(\O')} \leq (\widetilde{\mathcal C}'_{22} \delta^{-1} + \widetilde{\mathcal C}''_{22}) c(\varphi)^6 \eps \| \FF \|_{L_2(\O)}
\eqno(7.2)
$$
\textit{for $0< \eps \leq \eps_1$. The constants ${\mathcal C}'_{22}$, ${\mathcal C}''_{22}$, $\widetilde{\mathcal C}'_{22}$, and $\widetilde{\mathcal C}''_{22}$ depend on}
$d$, $m$, $\alpha_0$, $\alpha_1$, $\| g \|_{L_\infty}$, $\| g^{-1} \|_{L_\infty}$,
\textit{the parameters of the lattice} $\Gamma$, \textit{and the domain $\O$}.

\smallskip\noindent\textbf{Proof.}
We fix a smooth cut-off function $\chi(\x)$ such that
$$
\begin{aligned}
&\chi \in C_0^\infty(\O),\ \ 0\leq \chi(\x) \leq 1;
\cr
&\chi(\x)=1 \ \text{for}\ \x \in \O';\ \ |\nabla \chi(\x)| \leq \kappa' \delta^{-1}.
\end{aligned}
\eqno(7.3)
$$
The constant $\kappa'$ depends only on the domain $\O$.
Let $\u_\eps$ be the solution of problem (3.3), and let $\widetilde{\u}_\eps$ be the solution of equation (4.16).
Then \hbox{$(\A_\eps - \zeta)(\u_\eps - \widetilde{\u}_\eps)=0$} in the domain $\O$.
Hence, we have
$$
(g^\eps b(\D)(\u_\eps - \widetilde{\u}_\eps), b(\D) {\boldsymbol \eta})_{L_2(\O)}
- \zeta (\u_\eps - \widetilde{\u}_\eps, {\boldsymbol \eta})_{L_2(\O)} =0,\ \ \forall
{\boldsymbol \eta} \in H^1_0(\O;\C^n).
\eqno(7.4)
$$
We substitute ${\boldsymbol \eta} = \chi^2 (\u_\eps - \widetilde{\u}_\eps)$ in (7.4) and denote
$$
{\mathfrak A}(\eps):= (g^\eps b(\D)(\chi(\u_\eps - \widetilde{\u}_\eps)), b(\D)(\chi(\u_\eps - \widetilde{\u}_\eps)))_{L_2(\O)}.
\eqno(7.5)
$$
The corresponding identity can be easily transformed to the form
$$
\begin{aligned}
&{\mathfrak A}(\eps) - \zeta \| \chi(\u_\eps - \widetilde{\u}_\eps) \|^2_{L_2(\O)}
=- (g^\eps b(\D)(\chi(\u_\eps - \widetilde{\u}_\eps)), \z_\eps)_{L_2(\O)}
\cr
&+ (g^\eps \z_\eps, b(\D)(\chi(\u_\eps - \widetilde{\u}_\eps)))_{L_2(\O)}
+ (g^\eps \z_\eps,\z_\eps)_{L_2(\O)},
\end{aligned}
\eqno(7.6)
$$
where
$$
\z_\eps := \sum_{l=1}^d b_l (D_l \chi) (\u_\eps - \widetilde{\u}_\eps).
\eqno(7.7)
$$
The right-hand side of (7.6) is estimated by
$2 {\mathfrak A}(\eps)^{1/2} \|g\|^{1/2}_{L_\infty} \|\z_\eps \|_{L_2(\O)} + \|g\|_{L_\infty} \|\z_\eps \|_{L_2(\O)}^2$.
Taking the imaginary part in (7.6), we obtain
$$
|{\rm Im}\,\zeta| \| \chi(\u_\eps - \widetilde{\u}_\eps) \|^2_{L_2(\O)}
\leq 2 {\mathfrak A}(\eps)^{1/2} \|g\|^{1/2}_{L_\infty} \|\z_\eps \|_{L_2(\O)} + \|g\|_{L_\infty} \|\z_\eps \|_{L_2(\O)}^2.
\eqno(7.8)
$$
If ${\rm Re}\,\zeta \geq 0$ (and then ${\rm Im}\,\zeta \ne 0$), it follows that
$$
 \| \chi(\u_\eps - \widetilde{\u}_\eps) \|^2_{L_2(\O)}
\leq c(\varphi) |\zeta|^{-1} \left(2 {\mathfrak A}(\eps)^{1/2} \|g\|^{1/2}_{L_\infty} \|\z_\eps \|_{L_2(\O)} + \|g\|_{L_\infty} \|\z_\eps \|_{L_2(\O)}^2\right).
\eqno(7.9)
$$
If ${\rm Re}\,\zeta < 0$, we take the real part in  (7.6) and obtain
$$
 |{\rm Re}\,\zeta| \| \chi(\u_\eps - \widetilde{\u}_\eps) \|^2_{L_2(\O)}
\leq  2 {\mathfrak A}(\eps)^{1/2} \|g\|^{1/2}_{L_\infty} \|\z_\eps \|_{L_2(\O)} + \|g\|_{L_\infty} \|\z_\eps \|_{L_2(\O)}^2.
\eqno(7.10)
$$
Summing up (7.8) and (7.10), we deduce the inequality similar to (7.9).
As a result, for all values of $\zeta$ under consideration we have
$$
\begin{aligned}
 &\| \chi(\u_\eps - \widetilde{\u}_\eps) \|^2_{L_2(\O)}
\cr
&\leq 2 c(\varphi) |\zeta|^{-1} \left(2 {\mathfrak A}(\eps)^{1/2} \|g\|^{1/2}_{L_\infty} \|\z_\eps \|_{L_2(\O)} + \|g\|_{L_\infty} \|\z_\eps \|_{L_2(\O)}^2\right).
\end{aligned}
\eqno(7.11)
$$
Now, (7.6) and (7.11) imply that
$$
{\mathfrak A}(\eps) \leq 42 c(\varphi)^2 \|g\|_{L_\infty} \| \z_\eps \|^2_{L_2(\O)}.
\eqno(7.12)
$$
By (1.5) and (7.3), the function (7.7) satisfies
$$
\| \z_\eps \|_{L_2(\O)} \leq (d \alpha_1)^{1/2} \kappa' \delta^{-1} \| \u_\eps - \widetilde{\u}_\eps \|_{L_2(\O)}.
\eqno(7.13)
$$

From (7.5), (7.12), (7.13), and (3.1) it follows that
$$
\| \D (\chi(\u_\eps - \widetilde{\u}_\eps)) \|_{L_2(\O)}
\leq {\mathcal C}_{23} c(\varphi) \delta^{-1} \| \u_\eps - \widetilde{\u}_\eps \|_{L_2(\O)},
\eqno(7.14)
$$
where ${\mathcal C}_{23} = c_0^{-1/2} \|g\|_{L_\infty}^{1/2} (42 d \alpha_1)^{1/2} \kappa'$.

Relations (4.1) and (4.17) imply that
$$
\| \u_\eps - \widetilde{\u}_\eps \|_{L_2(\O)} \leq {\mathcal C}_{24} c(\varphi)^5
(|\zeta|^{-1/2}\eps +\eps^2) \|\FF\|_{L_2(\O)},\ \ 0< \eps \leq \eps_1,
\eqno(7.15)
$$
where ${\mathcal C}_{24} = {C}_{1} {\mathcal C}_{4} + {\mathcal C}_{1}$.
By (7.14) and (7.15),
$$
\| \D (\chi(\u_\eps - \widetilde{\u}_\eps)) \|_{L_2(\O)}
\leq {\mathcal C}_{23} {\mathcal C}_{24} \delta^{-1} c(\varphi)^6
(|\zeta|^{-1/2}\eps +\eps^2) \|\FF\|_{L_2(\O)},\  0< \eps \leq \eps_1.
$$

Hence,
$$
\| \u_\eps - \widetilde{\u}_\eps \|_{H^1(\O')} \leq {\mathcal C}_{24} ({\mathcal C}_{23} \delta^{-1}+1)c(\varphi)^6
(|\zeta|^{-1/2}\eps +\eps^2) \|\FF\|_{L_2(\O)},\ 0< \eps \leq \eps_1.
\eqno(7.16)
$$
From (4.6), (4.18), and (4.19) it follows that
$$
\|  \widetilde{\u}_\eps - \v_\eps \|_{H^1(\O)} \leq (C_2+C_3) {\mathcal C}_{4} c(\varphi)^3
\eps \|\FF\|_{L_2(\O)}.
\eqno(7.17)
$$
Relations (7.16) and (7.17) imply the required estimate (7.1) with
${\mathcal C}'_{22}= 2 {\mathcal C}_{23} {\mathcal C}_{24}$,
${\mathcal C}''_{22} = 2 {\mathcal C}_{24} + (C_2+C_3){\mathcal C}_{4}$.

Now we check (7.2). From (7.1) and (1.2), (1.5) it follows that
$$
\| \p_\eps - g^\eps b(\D)\v_\eps\|_{L_2(\O')} \leq
\|g\|_{L_\infty} (d \alpha_1)^{1/2} ({\mathcal C}'_{22} \delta^{-1} + {\mathcal C}''_{22}) c(\varphi)^6 \eps \| \FF \|_{L_2(\O)}.
$$
Combining this with (5.22), we obtain (7.2) with $\widetilde{\mathcal C}'_{22} = \|g\|_{L_\infty} (d \alpha_1)^{1/2} {\mathcal C}'_{22}$,
$\widetilde{\mathcal C}''_{22} = \|g\|_{L_\infty} (d \alpha_1)^{1/2} {\mathcal C}''_{22} + {\mathcal C}_{18}$. $\bullet$

\smallskip\noindent\textbf{7.2. The case where $\Lambda \in L_\infty$.}
Similarly, in the case where Condition 2.8 is satisfied, we obtain the following result.

\smallskip\noindent\textbf{Theorem 7.2.} \textit{Suppose that the assumptions of Theorem} 6.1
\textit{are satisfied. Let $\O'$ be a strictly interior subdomain of the domain $\O$, and let}
$\delta:= {\rm dist}\,\{ \O'; \partial \O \}$. \textit{Then for $0< \eps \leq \eps_1$ we have}
$$
\|\u_\eps - \check{\v}_\eps \|_{H^1(\O')} \leq ({\mathcal C}'_{22} \delta^{-1} + \check{\mathcal C}''_{22}) c(\varphi)^6 \eps \| \FF \|_{L_2(\O)},
\eqno(7.18)
$$
\textit{or, in operator terms,}
$$
\| (\A_{D,\eps} - \zeta I)^{-1} - (\A_{D}^0 - \zeta I)^{-1} - \eps K_D^0(\eps;\zeta) \|_{L_2(\O) \to H^1(\O')}
\leq ({\mathcal C}'_{22} \delta^{-1} + \check{\mathcal C}''_{22}) c(\varphi)^6 \eps.
$$
\textit{For the flux $\p_\eps = g^\eps b(\D) \u_\eps$ we have}
$$
\| \p_\eps - \widetilde{g}^\eps b(\D) {\u}_0
\|_{L_2(\O')} \leq (\widetilde{\mathcal C}' \delta^{-1} + \widehat{\mathcal C}'') c(\varphi)^6 \eps \| \FF \|_{L_2(\O)}
\eqno(7.19)
$$
\textit{for $0< \eps \leq \eps_1$. The constants ${\mathcal C}'_{22}$ and $\widetilde{\mathcal C}'_{22}$ are the same as in Theorem} 7.1.
\textit{The constants $\check{\mathcal C}''_{22}$ and
$\widehat{\mathcal C}''_{22}$ depend on}
$d$, $m$, $\alpha_0$, $\alpha_1$, $\| g \|_{L_\infty}$, $\| g^{-1} \|_{L_\infty}$,
\textit{the parameters of the lattice} $\Gamma$, \textit{the domain $\O$, and the norm $\|\Lambda\|_{L_\infty}$}.

\smallskip\noindent\textbf{Proof.}
Relations (4.13), (6.9), and (7.1) imply (7.18) with
$\check{\mathcal C}''_{22} = {\mathcal C}''_{22} + {\mathcal C}''' C_{\O}^{(2)} \widehat{c}$.

From (7.18) it follows that
$$
\| \p_\eps - g^\eps b(\D) \check{\v}_\eps\|_{L_2(\O')} \leq
\|g\|_{L_\infty} (d \alpha_1)^{1/2} ({\mathcal C}'_{22} \delta^{-1} + \check{\mathcal C}''_{22}) c(\varphi)^6 \eps \| \FF \|_{L_2(\O)}
$$
for $0< \eps \leq \eps_1$.
Combining this with (3.11) and (6.13), we obtain (7.19) with
$\widehat{\mathcal C}''_{22} = \|g\|_{L_\infty} (d\alpha_1)^{1/2} \check{\mathcal C}''_{22} + \widetilde{\mathcal C}' \widehat{c}$. $\bullet$

\section*{\S 8. Another approximation of the resolvent $(\A_{D,\eps} - \zeta I)^{-1}$}

\smallskip\noindent\textbf{8.1. Approximation of the resolvent $(\A_{D,\eps} - \zeta I)^{-1}$ for $\zeta\in \C \setminus [c_*,\infty)$.}
Theorems 4.2, 4.3, and 6.1 give two-parametric estimates with respect to $\eps$ and $|\zeta|$
uniform in the domain $\{\zeta \in \C: |\zeta|\geq 1,\ \varphi \in [\varphi_0, 2\pi - \varphi_0] \}$
with arbitrarily small $\varphi_0>0$.

For completeness, we obtain one more result about approximation of the resolvent
$(\A_{D,\eps} - \zeta I)^{-1}$ which is valid for a wider domain of the parameter $\zeta$.
This result may be preferable for bounded $|\zeta|$, and for points $\zeta$ with small $\varphi$
or $2\pi - \varphi$.

\smallskip\noindent\textbf{Theorem 8.1.} \textit{Let $\zeta \in \C \setminus [c_*,\infty)$, where $c_*>0$ is a common
lower bound of the operators $\A_{D,\eps}$ and $\A_D^0$. We put $\zeta - c_* = |\zeta - c_*|e^{i\psi}$ and denote}
$$
\rho_*(\zeta) =
\begin{cases}
c(\psi)^2 |\zeta -c_*|^{-2}, & |\zeta - c_*|<1, \\
c(\psi)^2, & |\zeta - c_*| \ge 1.
\end{cases}
\eqno(8.1)
$$
\textit{Let $\u_\eps$ be the solution of problem} (3.3),
\textit{and let $\u_0$ be the solution of problem} (3.10). \textit{Let $\v_\eps$ be defined by} (4.7).
\textit{Suppose that the number $\varepsilon_1$ satisfies Condition} 4.1.
\textit{Then for $0< \varepsilon \leq \varepsilon_1$ we have}
$$
\| \u_\eps - \u_0 \|_{L_2(\O)} \leq {\mathcal C}_{25} \rho_*(\zeta) \varepsilon \|\FF\|_{L_2(\O)},
\eqno(8.2)
$$
$$
\| \u_\eps - \v_\eps \|_{H^1(\O)} \leq
{\mathcal C}_{26} \rho_*(\zeta) \varepsilon^{1/2} \|\FF\|_{L_2(\O)},
\eqno(8.3)
$$
\textit{or, in operator terms,}
$$
\| (\A_{D,\varepsilon}- \zeta I)^{-1}  - (\A^0_{D}- \zeta I)^{-1} \|_{L_2({\mathcal O}) \to L_2({\mathcal O})}
\leq    {\mathcal C}_{25} \rho_*(\zeta) \varepsilon,
\eqno(8.4)
$$
$$
\| (\A_{D,\varepsilon}- \zeta I)^{-1}  - (\A^0_{D}- \zeta I)^{-1} - \eps K_D(\eps;\zeta)\|_{L_2({\mathcal O}) \to H^1({\mathcal O})}
\leq    {\mathcal C}_{26} \rho_*(\zeta) \varepsilon^{1/2}.
\eqno(8.5)
$$
\textit{For the flux $\p_\eps= g^\eps b(\D)\u_\eps$ we have}
$$
\| \p_\eps - \widetilde{g}^\eps S_\eps b(\D) \widetilde{\u}_0 \|_{L_2(\O)}
\leq    {\mathcal C}_{27} \rho_*(\zeta) \varepsilon^{1/2} \|\FF\|_{L_2(\O)}, \ \ 0< \eps \leq \eps_1.
\eqno(8.6)
$$
\textit{The constants ${\mathcal C}_{25}$, ${\mathcal C}_{26}$, and ${\mathcal C}_{27}$ depend on}
$d$, $m$, $\alpha_0$, $\alpha_1$, $\| g \|_{L_\infty}$, $\| g^{-1} \|_{L_\infty}$,
\textit{the parameters of the lattice} $\Gamma$, \textit{and the domain} $\O$.

\smallskip\noindent\textbf{Remark 8.2.}
1) The expression $c(\psi)^2 |\zeta -c_*|^{-2}$ in (8.1) is equal to
$(\text{dist}\,\{\zeta,[c_*,\infty)\})^{-2}$.
2) One can take $c_*=c_2$, where $c_2$ is defined by (3.2).
3) Let $\nu >0$ be an arbitrarily small number. If $\eps$ is sufficiently small,
one can take  $c_*=\lambda_1^0(D) - \nu$,
where $\lambda_1^0(D)$ is the first eigenvalue of $\A_D^0$.
4) It is easy to find the upper bound for $c_*$: from (3.1) it is seen that
$c_* \leq c_1 \mu_1^0(D)$, where $\mu_1^0(D)$ is the first eigenvalue of the operator
$-\Delta$ with the Dirichlet condition.
Therefore, $c_*$ is bounded by the number depending only on $\alpha_1$, $\|g\|_{L_\infty}$, and the domain $\O$.

\smallskip\noindent\textbf{Proof.}
We apply Theorem 4.2 with $\zeta =-1$. By (4.2),
$$
\| (\A_{D,\varepsilon}+ I)^{-1}  - (\A^0_{D}+ I)^{-1} \|_{L_2({\mathcal O}) \to L_2({\mathcal O})}
\leq    {\mathcal C}_1 (\varepsilon + \varepsilon^2)\leq 2{\mathcal C}_1 \varepsilon,\ 0< \eps \leq \eps_1.
$$
Using the analog of the identity (2.10) (with $\A_\eps$ replaced by $\A_{D,\eps}$ and $\A^0$ replaced by $\A^0_D$),
we see that
$$
\| (\A_{D,\varepsilon}- \zeta I)^{-1}  - (\A^0_{D}- \zeta I)^{-1} \|_{L_2({\mathcal O}) \to L_2({\mathcal O})}
\leq    2{\mathcal C}_1 \varepsilon \sup_{x\geq c_*} (x+1)^2 |x-\zeta|^{-2}
\eqno(8.7)
$$
for $0< \eps \leq \eps_1$. A calculation shows that
$$
\sup_{x\geq c_*} (x+1)^2 |x-\zeta|^{-2} \leq \check{c} \rho_*(\zeta), \ \ \zeta \in \C \setminus [c_*,\infty),
\eqno(8.8)
$$
where $\check{c}= c_*^2 + 4 c_* + 5$. By Remark 8.2(4), $\check{c}$
is bounded by the number depending only on $\alpha_1$, $\|g\|_{L_\infty}$, and the domain $\O$.
Now (8.7) and (8.8) imply (8.4) with ${\mathcal C}_{25} = 2 {\mathcal C}_1 \check{c}$.

Now we apply Theorem 4.3 with $\zeta=-1$. By (4.9),
$$
\begin{aligned}
&\| (\A_{D,\varepsilon}+ I)^{-1}  - (\A^0_{D} + I)^{-1} - \eps K_D(\eps;-1)\|_{L_2({\mathcal O}) \to H^1({\mathcal O})}
\cr
&\leq    {\mathcal C}_2 \varepsilon^{1/2}+ {\mathcal C}_3 \varepsilon \leq ({\mathcal C}_2 + {\mathcal C}_3) \varepsilon^{1/2}
\end{aligned}
\eqno(8.9)
$$
for $0< \eps \leq \eps_1$. From Lemma  5.2 with $\zeta=-1$ it follows that 
$$
\|  \eps \theta_\eps K_D(\eps;-1)\|_{L_2({\mathcal O}) \to H^1({\mathcal O})}
\leq ({\mathcal C}_{13}+ {\mathcal C}_{14} +{\mathcal C}_{15})\eps^{1/2},\ \ 0< \eps \leq \eps_1.
\eqno(8.10)
$$
By (8.9) and (8.10), we have
$$
\| (\A_{D,\varepsilon}+ I)^{-1}  - (\A^0_{D} + I)^{-1} - \eps (1 - \theta_\eps) K_D(\eps;-1)\|_{L_2({\mathcal O}) \to H^1({\mathcal O})}
\leq    {\mathcal C}_{28} \varepsilon^{1/2}
\eqno(8.11)
$$
for $0< \eps \leq \eps_1$, with ${\mathcal C}_{28} = {\mathcal C}_2 + {\mathcal C}_3 + {\mathcal C}_{13}+ {\mathcal C}_{14} +{\mathcal C}_{15}$.
We use the following analog of the identity (2.34):
$$
\begin{aligned}
&(\A_{D,\varepsilon}-\zeta I)^{-1}  - (\A^0_{D} - \zeta I)^{-1} - \eps (1-\theta_\eps) K_D(\eps;\zeta)
\cr
 &= (\A_{D,\varepsilon}+ I)(\A_{D,\varepsilon}-\zeta I)^{-1}
\cr
&\times \left((\A_{D,\varepsilon}+ I)^{-1}  - (\A^0_{D} + I)^{-1} - \eps (1-\theta_\eps)K_D(\eps;-1) \right)
\cr
&\times (\A_{D}^0 + I)(\A_{D}^0 -\zeta I)^{-1} + \eps (\zeta+1) (\A_{D,\varepsilon}-\zeta I)^{-1}(1-\theta_\eps) K_D(\eps;\zeta).
\end{aligned}
\eqno(8.12)
$$
Since the range of the operators in  (8.12) is contained in $H^1_0(\O;\C^n)$, we can multiply (8.12) by $\A_{D,\eps}^{1/2}$ from the left.
Taking (8.8) into account, we obtain
 $$
\begin{aligned}
&\| \A_{D,\varepsilon}^{1/2} \left((\A_{D,\varepsilon}-\zeta I)^{-1}  - (\A^0_{D} - \zeta I)^{-1} - \eps (1-\theta_\eps)K_D(\eps;\zeta) \right) \|_{L_2 \to L_2}
\cr
&\leq \check{c} \rho_*(\zeta) \| \A_{D,\varepsilon}^{1/2} \left((\A_{D,\varepsilon}+ I)^{-1}  - (\A^0_{D} + I)^{-1} - \eps (1-\theta_\eps)K_D(\eps;-1) \right)\|_{L_2 \to L_2}
\cr
&+ \eps |\zeta+1| \sup_{x\geq c_*} x^{1/2} |x-\zeta|^{-1} \|(1-\theta_\eps)K_D(\eps;\zeta)\|_{L_2({\mathcal O}) \to L_2({\mathcal O})}.
\end{aligned}
\eqno(8.13)
$$
Denote the summands in the right-hand side of (8.13) by ${\mathcal L}_1(\eps)$ and ${\mathcal L}_2(\eps)$.
The first term is estimated with the help of (8.11) and (3.1):
$$
{\mathcal L}_1(\eps) \leq c_1^{1/2} \check{c}\, {\mathcal C}_{28} \rho_*(\zeta) \eps^{1/2},\ \ 0< \eps \leq \eps_1.
\eqno(8.14)
$$
The operator $K_D(\eps;\zeta)$ can be represented as
$K_D(\eps;\zeta) = R_\O [\Lambda^\eps] S_\eps b(\D) P_\O (\A_D^0)^{-1/2}(\A_D^0)^{1/2}(\A_D^0 - \zeta I)^{-1}$.
Hence, by (1.4), (1.19), (4.3), and (5.1), we obtain
$$
{\mathcal L}_2(\eps) \leq \eps |\zeta+1| \left(\sup_{x\geq c_*} x |x-\zeta|^{-2}\right)
M_1 \alpha_1^{1/2} C_\O^{(1)} \| (\A_D^0)^{-1/2}\|_{L_2({\mathcal O}) \to H^1({\mathcal O})}.
\eqno(8.15)
$$
Using the analogs of the estimates (3.1) and (3.2) for $\A_D^0$, we have
$$
\| (\A_D^0)^{-1/2}\|_{L_2({\mathcal O}) \to H^1({\mathcal O})} \leq (c_0^{-1} + c_2^{-1})^{1/2}.
\eqno(8.16)
$$
A calculation shows that
$$
\sup_{x\geq c_*} x |x-\zeta|^{-2} \leq
\begin{cases}
(c_* +1) c(\psi)^2 |\zeta -c_*|^{-2}, & |\zeta - c_*|<1, \\
(c_* +1) c(\psi)^2 |\zeta -c_*|^{-1}, & |\zeta - c_*| \ge 1.
\end{cases}
$$
Note that $|\zeta+1| \leq 2+c_*$ for $|\zeta -c_*| < 1$,
and $|\zeta+1| |\zeta -c_*|^{-1} \leq 2+c_*$ for $|\zeta -c_*| \geq 1$, whence
$$
|\zeta+1| \sup_{x\geq c_*} x |x-\zeta|^{-2} \leq (c_* +2)(c_* +1) \rho_*(\zeta).
\eqno(8.17)
$$
From (8.15)--(8.17) it follows that
$$
{\mathcal L}_2(\eps) \leq {\mathcal C}_{29} \rho_*(\zeta) \eps,
\eqno(8.18)
$$
where ${\mathcal C}_{29} = (c_* +2)(c_* +1) M_1 \alpha_1^{1/2} C_\O^{(1)} (c_0^{-1} + c_2^{-1})^{1/2}.$

As a result, inequalities (8.13), (8.14), and (8.18) yield
$$
\begin{aligned}
&\| \A_{D,\varepsilon}^{1/2} \left((\A_{D,\varepsilon}-\zeta I)^{-1}  - (\A^0_{D} - \zeta I)^{-1} - \eps (1-\theta_\eps)K_D(\eps;\zeta) \right) \|_{L_2({\mathcal O}) \to L_2({\mathcal O})}
\cr
&\leq (c_1^{1/2}\check{c}\, {\mathcal C}_{28} + {\mathcal C}_{29}) \rho_*(\zeta) \eps^{1/2},\ \ 0< \eps \leq \eps_1.
\end{aligned}
$$
Together with (3.1) and (3.2) this implies that
$$
\begin{aligned}
&\| (\A_{D,\varepsilon}-\zeta I)^{-1}  - (\A^0_{D} - \zeta I)^{-1} - \eps (1-\theta_\eps)K_D(\eps;\zeta) \|_{L_2({\mathcal O}) \to H^1({\mathcal O})}
\cr
&\leq (c_0^{-1}+ c_2^{-1})^{1/2} (c_1^{1/2}\check{c}\,{\mathcal C}_{28} + {\mathcal C}_{29}) \rho_*(\zeta) \eps^{1/2},\ \ 0< \eps \leq \eps_1.
\end{aligned}
\eqno(8.19)
$$
Finally, by (8.10) and (8.8), we have
$$
\begin{aligned}
&\| \eps \theta_\eps K_D(\eps;\zeta) \|_{L_2({\mathcal O}) \to H^1({\mathcal O})}
\leq \| \eps \theta_\eps K_D(\eps;-1) \|_{L_2({\mathcal O}) \to H^1({\mathcal O})}
\cr
&\times \| (\A_D^0 +I) (\A_D^0 - \zeta I)^{-1} \|_{L_2({\mathcal O}) \to L_2({\mathcal O})}
\leq ({\mathcal C}_{13}+ {\mathcal C}_{14} + {\mathcal C}_{15}) \check{c}^{1/2}\rho_*(\zeta)^{1/2} \eps^{1/2}
\end{aligned}
\eqno(8.20)
$$
for $0< \eps \leq \eps_1$. As a result, relations (8.19) and (8.20) imply (8.5) with
${\mathcal C}_{26} = (c_0^{-1} + c_2^{-1})^{1/2} (c_1^{1/2}\check{c}\,{\mathcal C}_{28} + {\mathcal C}_{29}) +
({\mathcal C}_{13}+ {\mathcal C}_{14} + {\mathcal C}_{15}) \check{c}^{1/2}$.

It remains to check (8.6). From (8.3) and (1.2), (1.5) it is seen that
$$
\| \p_\eps - g^\eps b(\D) \v_\eps \|_{L_2(\O)} \leq
\|g\|_{L_\infty} (d \alpha_1)^{1/2} {\mathcal C}_{26} \rho_*(\zeta) \eps^{1/2}\|\FF\|_{L_2(\O)}
\eqno(8.21)
$$
 for $0< \eps \leq \eps_1$. Next, similarly to (5.19)--(5.21) we obtain
$$
\| g^\eps b(\D) \v_\eps - \widetilde{g}^\eps S_\eps b(\D) \widetilde{\u}_0 \|_{L_2(\O)}
\leq ({\mathcal C}' +{\mathcal C}'') \eps \| \widetilde{\u}_0 \|_{H^2(\R^d)}.
\eqno(8.22)
$$
From (3.9) and (8.8) it follows that
$$
\| (\A_D^0 - \zeta I)^{-1}\|_{L_2(\O) \to H^2(\O)} \leq \widehat{c} \sup_{x\geq c_*} x |x-\zeta|^{-1} \leq
\widehat{c} \,\check{c}^{1/2} \rho_*(\zeta)^{1/2}.
$$
Hence,
$$
\|\u_0\|_{H^2(\O)} \leq \widehat{c} \,\check{c}^{1/2} \rho_*(\zeta)^{1/2} \|\FF\|_{L_2(\O)}.
\eqno(8.23)
$$
Together with (4.3) and (8.22) this yields
$$
\| g^\eps b(\D) \v_\eps - \widetilde{g}^\eps S_\eps b(\D) \widetilde{\u}_0 \|_{L_2(\O)}
\leq \widetilde{\mathcal C}_{27}  \rho_*(\eps)^{1/2} \eps \|\FF \|_{L_2(\O)},
\eqno(8.24)
$$
where $\widetilde{\mathcal C}_{27}= ({\mathcal C}' +{\mathcal C}'') C_\O^{(2)} \widehat{c} \,\check{c}^{1/2}$.
Combining this with (8.21), we obtain (8.6) with
${\mathcal C}_{27} = \|g\|_{L_\infty} (d \alpha_1)^{1/2} {\mathcal C}_{26}
+ \widetilde{\mathcal C}_{27}$. $\bullet$

\smallskip\noindent\textbf{8.2. The case where $\Lambda \in L_\infty$.}

\smallskip\noindent\textbf{Theorem 8.3.} \textit{Suppose that the assumptions of Theorem} 8.1 \textit{and Condition} 2.8
\textit{are satisfied. Let $\check{\v}_\eps$ be defined by} (6.2).
\textit{Then for $0< \varepsilon \leq \varepsilon_1$ we have}
$$
\| \u_\eps - \check{\v}_\eps \|_{H^1(\O)} \leq
{\mathcal C}_{26}^\circ \rho_*(\zeta) \varepsilon^{1/2} \|\FF\|_{L_2(\O)},
\eqno(8.25)
$$
\textit{or, in operator terms,}
$$
\| (\A_{D,\varepsilon}- \zeta I)^{-1}  - (\A^0_{D}- \zeta I)^{-1} - \eps K_D^0(\eps;\zeta)\|_{L_2({\mathcal O}) \to H^1({\mathcal O})}
\leq    {\mathcal C}_{26}^\circ \rho_*(\zeta) \varepsilon^{1/2}.
$$
\textit{For the flux $\p_\eps= g^\eps b(\D)\u_\eps$ we have}
$$
\| \p_\eps - \widetilde{g}^\eps b(\D) {\u}_0 \|_{L_2(\O)}
\leq    {\mathcal C}_{27}^\circ \rho_*(\zeta) \varepsilon^{1/2} \|\FF\|_{L_2(\O)}, \ \ 0< \eps \leq \eps_1.
\eqno(8.26)
$$
\textit{The constants ${\mathcal C}_{26}^\circ$ and ${\mathcal C}_{27}^\circ$ depend on}
$d$, $m$, $\alpha_0$, $\alpha_1$, $\| g \|_{L_\infty}$, $\| g^{-1} \|_{L_\infty}$,
\textit{the parameters of the lattice} $\Gamma$, \textit{the domain} $\O$, \textit{and the norm} $\|\Lambda \|_{L_\infty}$.

\smallskip\noindent\textbf{Proof.}
Similarly to the proof of Theorem 6.1 (see (6.5)--(6.9)), it is easy to check that
$$
\| \v_\eps - \check{\v}_\eps\|_{H^1(\O)} \leq
{\mathcal C}''' \eps \| \widetilde{\u}_0\|_{H^2(\R^d)}.
$$
Together with (4.3) and (8.23) this yields
$$
\| \v_\eps - \check{\v}_\eps\|_{H^1(\O)} \leq
{\mathcal C}''' C_{\O}^{(2)} \widehat{c}\, \check{c}^{1/2} \rho_*(\zeta)^{1/2} \eps \|\FF \|_{L_2(\O)}.
\eqno(8.27)
$$
Now relations (8.3) and (8.27) imply the required estimate (8.25) with
${\mathcal C}_{26}^\circ = {\mathcal C}_{26} + {\mathcal C}''' C_{\O}^{(2)} \widehat{c}\, \check{c}^{1/2}$.

It remains to check (8.26). From (8.25) and (1.2), (1.5) it follows that
$$
\| \p_\eps - g^\eps b(\D) \check{\v}_\eps \|_{L_2(\O)} \leq
\|g\|_{L_\infty} (d \alpha_1)^{1/2} {\mathcal C}_{26}^\circ \rho_*(\zeta) \eps^{1/2}\|\FF\|_{L_2(\O)},
\ \ 0< \eps \leq \eps_1.
\eqno(8.28)
$$
Similarly to  (6.11)--(6.13), we have
$$
\|  g^\eps b(\D) \check{\v}_\eps - \widetilde{g}^\eps b(\D) \u_0 \|_{L_2(\O)} \leq
\widetilde{\mathcal C}' \eps \|\u_0 \|_{H^2(\O)}.
\eqno(8.29)
$$
As a result, relations (8.23), (8.28), and (8.29) imply (8.26) with
${\mathcal C}_{27}^\circ = \|g\|_{L_\infty} (d \alpha_1)^{1/2} {\mathcal C}_{26}^\circ
+ \widetilde{\mathcal C}' \widehat{c} \check{c}^{1/2}$.
$\bullet$

\smallskip\noindent\textbf{8.3. Special cases.}
The following statement is similar to Proposition 6.2.

\smallskip\noindent\textbf{Proposition 8.4.} \textit{Suppose that the assumptions of Theorem} 8.1
 \textit{are satisfied. If $g^0 = \overline{g}$, i.~e., relations} (1.13) \textit{are satisfied, then
 $\Lambda=0$, $\v_\eps = \u_0$, and for $0< \varepsilon \leq \varepsilon_1$ we have}
$$
\| \u_\eps - \u_0 \|_{H^1(\O)} \leq
{\mathcal C}_{26} \rho_*(\zeta) \varepsilon^{1/2} \|\FF\|_{L_2(\O)}.
$$

\smallskip
Next statement is similar to Proposition 6.3.

\smallskip\noindent\textbf{Proposition 8.5.} \textit{Suppose that the assumptions of Theorem} 8.1
 \textit{are satisfied. If $g^0 = \underline{g}$, i.~e., relations} (1.14) \textit{are satisfied, then
for $0< \varepsilon \leq \varepsilon_1$ we have}
$$
\| \p_\eps - g^0 b(\D){\u}_0 \|_{L_2(\O)}
\leq    {\mathcal C}_{27}^\circ \rho_*(\zeta) \varepsilon^{1/2} \|\FF\|_{L_2(\O)}.
$$

\smallskip\noindent\textbf{8.4. Estimates in a strictly interior subdomain.}
Similarly to Theorem 7.1
we obtain an error estimate of order $\eps$ in $H^1(\O')$ in a strictly interior subdomain $\O'$ of the domain $\O$,
using Theorem 8.1 and the results for the problem in $\R^d$.

\smallskip\noindent\textbf{Theorem 8.6.} \textit{Suppose that the assumptions of Theorem} 8.1
\textit{are satisfied. Let $\O'$ be a strictly interior subdomain of the domain $\O$, and let $\delta:= {\rm dist}\,\{ \O'; \partial \O \}$.}
\textit{Then for $0< \varepsilon \leq \varepsilon_1$ we have}
$$
\begin{aligned}
&\| \u_\eps - \v_\eps \|_{H^1(\O')}
\cr
&\leq
\left({\mathcal C}'_{30}\delta^{-1} (c(\psi) \rho_*(\zeta) +c(\psi)^{5/2} \rho_*(\zeta)^{3/4}) +
{\mathcal C}''_{30} c(\psi)^{1/2} \rho_*(\zeta)^{5/4}\right) \varepsilon \|\FF\|_{L_2(\O)},
\end{aligned}
\eqno(8.30)
$$
\textit{or, in operator terms,}
$$
\begin{aligned}
&\| (\A_{D,\varepsilon}- \zeta I)^{-1}  - (\A^0_{D}- \zeta I)^{-1} - \eps K_D(\eps;\zeta)\|_{L_2({\mathcal O}) \to H^1({\mathcal O}')}
\cr
&\leq
\left({\mathcal C}'_{30}\delta^{-1} (c(\psi) \rho_*(\zeta) +c(\psi)^{5/2} \rho_*(\zeta)^{3/4}) +
{\mathcal C}''_{30} c(\psi)^{1/2} \rho_*(\zeta)^{5/4}\right) \varepsilon.
\end{aligned}
$$
\textit{For the flux $\p_\eps= g^\eps b(\D)\u_\eps$ we have}
$$
\begin{aligned}
&\| \p_\eps - \widetilde{g}^\eps S_\eps b(\D) \widetilde{\u}_0 \|_{L_2(\O')}
\cr
&\leq
\left(\widetilde{\mathcal C}'_{30}\delta^{-1} (c(\psi) \rho_*(\zeta) +c(\psi)^{5/2} \rho_*(\zeta)^{3/4})
+ \widetilde{\mathcal C}''_{30} c(\psi)^{1/2} \rho_*(\zeta)^{5/4}\right
) \varepsilon \|\FF\|_{L_2(\O)}
\end{aligned}
\eqno(8.31)
$$
\textit{for $0< \eps \leq \eps_1$. The constants ${\mathcal C}'_{30}$, ${\mathcal C}''_{30}$,
$\widetilde{\mathcal C}'_{30}$, and $\widetilde{\mathcal C}''_{30}$ depend on} $d$, $m$, $\alpha_0$, $\alpha_1$, $\| g \|_{L_\infty}$,
$\| g^{-1} \|_{L_\infty}$,
\textit{the parameters of the lattice} $\Gamma$, \textit{and the domain} $\O$.

\smallskip\noindent\textbf{Proof.}
In general, the proof is similar to that of Theorem 7.1. However, the associated problem in $\R^d$ is chosen in a different way.
Let $\u_\eps$ be the solution of problem (3.3), and let $\u_0$ be the solution of problem (3.10).
Let $\widetilde{\u}_0 = P_\O \u_0$. Since $\| (\A_D^0 -\zeta I)^{-1}\|_{L_2(\O)\to L_2(\O)} \leq c(\psi) |\zeta -c_*|^{-1}$,
by (4.3), we have
$$
\| \widetilde{\u}_0 \|_{L_2(\R^d)} \leq C_\O^{(0)} \| {\u}_0 \|_{L_2(\O)} \leq
C_\O^{(0)} c(\psi) |\zeta - c_*|^{-1} \| \FF \|_{L_2(\O)}.
\eqno(8.32)
$$
By (4.3) and (8.23),
$$
\| \widetilde{\u}_0 \|_{H^2(\R^d)} \leq C_\O^{(2)} \widehat{c} \, \check{c}^{1/2} \rho_*(\zeta)^{1/2} \| \FF \|_{L_2(\O)}.
\eqno(8.33)
$$
We put
$$
\widehat{\FF} := \A^0 \widetilde{\u}_0 - (\zeta - c_*) \widetilde{\u}_0.
\eqno(8.34)
$$
Similarly to (4.15), from (8.32) and (8.33) we deduce that
$$
\| \widehat{\FF} \|_{L_2(\R^d)} \leq {\mathcal C}_{31} \rho_*(\zeta)^{1/2} \| \FF \|_{L_2(\O)},
\eqno(8.35)
$$
where ${\mathcal C}_{31} = c_1 C_\O^{(2)} \widehat{c} \, \check{c}^{1/2} + C_\O^{(0)}.$
Note that $(\widehat{\FF} - \FF)\vert_{\O} = c_* \u_0$.

Under our assumptions, $\zeta \in \C \setminus [c_*,\infty)$.
Then the point $\zeta - c_* \in \C \setminus [0,\infty)$ is regular for the operator $\A_\eps$.
Let $\widehat{\u}_\eps$ be the solution of equation
$$
\A_\eps \widehat{\u}_\eps - (\zeta -c_*)\widehat{\u}_\eps = \widehat{\FF}
\eqno(8.36)
$$
in $\R^d$. Then the function $\u_\eps - \widehat{\u}_\eps$ satisfies the identity
$$
\begin{aligned}
&(g^\eps b(\D)(\u_\eps - \widehat{\u}_\eps), b(\D) {\boldsymbol \eta})_{L_2(\O)}
- (\zeta -c_*) (\u_\eps - \widehat{\u}_\eps, {\boldsymbol \eta})_{L_2(\O)}
\cr
&= c_* ( \u_\eps - \u_0, {\boldsymbol \eta})_{L_2(\O)},
\ \ \forall
{\boldsymbol \eta} \in H^1_0(\O;\C^n).
\end{aligned}
\eqno(8.37)
$$
Let $\chi(\x)$ be a cut-off function satisfying  (7.3).
We substitute ${\boldsymbol \eta} = \chi^2 (\u_\eps - \widehat{\u}_\eps)$ in (8.37)
and denote
$$
{\mathfrak B}(\eps):= (g^\eps b(\D)(\chi(\u_\eps - \widehat{\u}_\eps)), b(\D)(\chi(\u_\eps - \widehat{\u}_\eps)))_{L_2(\O)}.
$$
Similarly to (7.6), the corresponding relation can be rewritten as
$$
\begin{aligned}
&{\mathfrak B}(\eps) - (\zeta-c_*) \| \chi(\u_\eps - \widehat{\u}_\eps) \|^2_{L_2(\O)}
=c_* (\u_\eps - \u_0, \chi^2 (\u_\eps - \widehat{\u}_\eps))_{L_2(\O)}
\cr
&- (g^\eps b(\D)(\chi(\u_\eps - \widehat{\u}_\eps)), \r_\eps)_{L_2(\O)}
\cr
&+ (g^\eps \r_\eps, b(\D)(\chi(\u_\eps - \widehat{\u}_\eps)))_{L_2(\O)}
+ (g^\eps \r_\eps,\r_\eps)_{L_2(\O)},
\end{aligned}
\eqno(8.38)
$$
where
$$
\r_\eps := \sum_{l=1}^d b_l (D_l \chi) (\u_\eps - \widehat{\u}_\eps).
\eqno(8.39)
$$
Take the imaginary part in (8.38). Then
$$
\begin{aligned}
&|{\rm Im}\,\zeta| \| \chi(\u_\eps - \widehat{\u}_\eps) \|^2_{L_2(\O)}
\leq
c_* \|\u_\eps - \u_0\|_{L_2(\O)}  \| \chi(\u_\eps - \widehat{\u}_\eps) \|_{L_2(\O)}
\cr
&+
2 {\mathfrak B}(\eps)^{1/2} \|g\|^{1/2}_{L_\infty} \|\r_\eps \|_{L_2(\O)} + \|g\|_{L_\infty} \|\r_\eps \|_{L_2(\O)}^2.
\end{aligned}
\eqno(8.40)
$$
If ${\rm Re}\,\zeta \geq c_*$ (and then ${\rm Im}\,\zeta \ne 0$), we deduce that
$$
\begin{aligned}
&\| \chi(\u_\eps - \widehat{\u}_\eps) \|^2_{L_2(\O)}
\leq c(\psi)^2 |\zeta -c_*|^{-2} c_*^2 \|\u_\eps - \u_0\|_{L_2(\O)}^2
\cr
&+  c(\psi) |\zeta -c_*|^{-1} \left(4 {\mathfrak B}(\eps)^{1/2} \|g\|^{1/2}_{L_\infty} \|\r_\eps \|_{L_2(\O)} +
2 \|g\|_{L_\infty} \|\r_\eps \|_{L_2(\O)}^2\right).
\end{aligned}
\eqno(8.41)
$$
If ${\rm Re}\,\zeta < c_*$, we take the real part in (8.38) and obtain
$$
\begin{aligned}
&|{\rm Re}\,\zeta - c_*| \| \chi(\u_\eps - \widehat{\u}_\eps) \|^2_{L_2(\O)}
\leq
c_* \|\u_\eps - \u_0\|_{L_2(\O)}  \| \chi(\u_\eps - \widehat{\u}_\eps) \|_{L_2(\O)}
\cr
&+ 2 {\mathfrak B}(\eps)^{1/2} \|g\|^{1/2}_{L_\infty} \|\r_\eps \|_{L_2(\O)} + \|g\|_{L_\infty} \|\r_\eps \|_{L_2(\O)}^2.
\end{aligned}
\eqno(8.42)
$$
Summing up (8.40) and (8.42), we deduce the inequality similar to (8.41).
As a result, for all values of $\zeta$ under consideration we obtain
$$
\begin{aligned}
 &\| \chi(\u_\eps - \widehat{\u}_\eps) \|^2_{L_2(\O)}
\leq
4 c(\psi)^2 |\zeta -c_*|^{-2} c_*^2 \|\u_\eps - \u_0\|_{L_2(\O)}^2
\cr
&+
c(\psi) |\zeta -c_*|^{-1} \left(8 {\mathfrak B}(\eps)^{1/2} \|g\|^{1/2}_{L_\infty} \|\r_\eps \|_{L_2(\O)} + 4 \|g\|_{L_\infty} \|\r_\eps \|_{L_2(\O)}^2\right).
\end{aligned}
\eqno(8.43)
$$
Relations (8.38) and (8.43) imply that
$$
{\mathfrak B}(\eps) \leq 342 c(\psi)^2 \|g\|_{L_\infty} \| \r_\eps \|^2_{L_2(\O)}
+ 18 c_*^2 c(\psi)^2 |\zeta - c_*|^{-1} \|\u_\eps - \u_0\|^2_{L_2(\O)}.
\eqno(8.44)
$$
By (1.5) and (7.3), the function (8.39) satisfies the estimate
$$
\| \r_\eps \|_{L_2(\O)} \leq (d \alpha_1)^{1/2} \kappa' \delta^{-1} \| \u_\eps - \widehat{\u}_\eps \|_{L_2(\O)}.
\eqno(8.45)
$$
From (8.44), (8.45), and (3.1) it follows that
$$
\begin{aligned}
&\|\D (\chi(\u_\eps - \widehat{\u}_\eps)) \|_{L_2(\O)} \leq {\mathcal C}_{32} c(\psi) \delta^{-1} \|\u_\eps - \widehat{\u}_\eps \|_{L_2(\O)}
\cr
&+ \sqrt{18} c_0^{-1/2} c_* c(\psi) |\zeta - c_*|^{-1/2} \|\u_\eps - {\u}_0 \|_{L_2(\O)},
\end{aligned}
\eqno(8.46)
$$
where ${\mathcal C}_{32}^2 = 342 c_0^{-1} \|g\|_{L_\infty} d\alpha_1 (\kappa')^2$.

The norm $\|\u_\eps - {\u}_0 \|_{L_2(\O)}$ satisfies (8.2).
In order to estimate \hbox{$\|\widehat{\u}_\eps - \widetilde{\u}_0 \|_{L_2(\R^d)}$}, we apply Theorem
2.2, taking (8.34) and (8.36) into account. Using also (8.35), we obtain
$$
\begin{aligned}
&\|\widehat{\u}_\eps - {\u}_0 \|_{L_2(\O)}  \leq
\|\widehat{\u}_\eps - \widetilde{\u}_0 \|_{L_2(\R^d)} \leq
C_1 c(\psi)^2 |\zeta - c_*|^{-1/2} \eps \| \widehat{\FF} \|_{L_2(\R^d)}
\cr
&\leq C_1 {\mathcal C}_{31} c(\psi)^2 \rho_*(\zeta)^{1/2} |\zeta - c_*|^{-1/2} \eps \| {\FF} \|_{L_2(\O)}.
\end{aligned}
\eqno(8.47)
$$
By (8.1), (8.2), and (8.47),
$$
\|\u_\eps - \widehat{\u}_\eps \|_{L_2(\O)} \leq \left({\mathcal C}_{25} \rho_*(\zeta) + C_1 {\mathcal C}_{31}  c(\psi)^{3/2} \rho_*(\zeta)^{3/4}\right)
\eps \|\FF\|_{L_2(\O)}.
\eqno(8.48)
$$
Now relations (8.2), (8.46), and (8.48) yield
$$
\begin{aligned}
&\|\D (\chi(\u_\eps - \widehat{\u}_\eps)) \|_{L_2(\O)}
\cr
&\leq
\left({\mathcal C}'_{30} \delta^{-1} (c(\psi) \rho_*(\zeta) + c(\psi)^{5/2} \rho_*(\zeta)^{3/4}) +
{\mathcal C}_{33} c(\psi)^{1/2} \rho_*(\zeta)^{5/4}\right) \eps \|\FF\|_{L_2(\O)},
\end{aligned}
$$
where ${\mathcal C}'_{30} = {\mathcal C}_{32} \max \{ {\mathcal C}_{25}, C_1 {\mathcal C}_{31} \}$ and
${\mathcal C}_{33} = \sqrt{18} c_0^{-1/2} c_* {\mathcal C}_{25}$.
Together with (8.48) this implies that
$$
\begin{aligned}
&\|\u_\eps - \widehat{\u}_\eps \|_{H^1(\O')}
\cr
&\leq \left({\mathcal C}_{30}' \delta^{-1} (c(\psi) \rho_*(\zeta) + c(\psi)^{5/2} \rho_*(\zeta)^{3/4})
+ {\mathcal C}_{34}  c(\psi)^{1/2} \rho_*(\zeta)^{5/4}\right) \eps \|\FF\|_{L_2(\O)},
\end{aligned}
\eqno(8.49)
$$
where ${\mathcal C}_{34} = {\mathcal C}_{33} + {\mathcal C}_{25} + C_1 {\mathcal C}_{31}$.

By Corollary 2.5 and (8.35),
$$
\begin{aligned}
&\|\widehat{\u}_\eps - \widetilde{\u}_0 - \eps \Lambda^\eps S_\eps b(\D)\widetilde{\u}_0\|_{H^1(\R^d)}
\leq  c(\psi)^2 (C_2 + C_3 |\zeta - c_*|^{-1/2}) \eps \|\widehat{\FF}\|_{L_2(\R^d)}
\cr
&\leq {\mathcal C}_{31} (C_2+C_3)c(\psi)^{3/2} \rho_*(\zeta)^{3/4} \eps \|{\FF}\|_{L_2(\O)}.
\end{aligned}
\eqno(8.50)
$$
As a result, relations (8.49) and (8.50) imply the required estimate (8.30) with ${\mathcal C}_{30}'' =
{\mathcal C}_{34} + {\mathcal C}_{31} (C_2+C_3)$.

It remains to check (8.31). From (8.30) and (1.2), (1.5) it follows that
$$
\begin{aligned}
&\| \p_\eps - g^\eps b(\D)\v_\eps \|_{L_2(\O')} \leq
\|g\|_{L_\infty} (d \alpha_1)^{1/2}
\cr
&\times \left({\mathcal C}'_{30} \delta^{-1} (c(\psi) \rho_*(\zeta) + c(\psi)^{5/2} \rho_*(\zeta)^{3/4}) +
{\mathcal C}''_{30} c(\psi)^{1/2} \rho_*(\zeta)^{5/4}\right) \varepsilon \|\FF\|_{L_2(\O)}.
\end{aligned}
$$
Combining this with (8.24), we obtain (8.31) with
$\widetilde{\mathcal C}'_{30} = \|g\|_{L_\infty} (d \alpha_1)^{1/2} {\mathcal C}'_{30}$ and
$\widetilde{\mathcal C}''_{30} = \|g\|_{L_\infty} (d \alpha_1)^{1/2} {\mathcal C}''_{30} +
\widetilde{\mathcal C}_{27}$. $\bullet$

\smallskip
Similarly, in the case where Condition 2.8 is satisfied, we obtain the following result.

\smallskip\noindent\textbf{Theorem 8.7.} \textit{Suppose that the assumptions of Theorem} 8.3
\textit{are satisfied. Let $\O'$ be a strictly interior subdomain of the domain $\O$, and let
$\delta:= {\rm dist}\,\{ \O'; \partial \O \}$.}
\textit{Then for $0< \varepsilon \leq \varepsilon_1$ we have}
$$
\begin{aligned}
&\| \u_\eps - \check{\v}_\eps \|_{H^1(\O')}
\cr
&\leq
\left({\mathcal C}'_{30}\delta^{-1} (c(\psi) \rho_*(\zeta) + c(\psi)^{5/2} \rho_*(\zeta)^{3/4}) +
\check{\mathcal C}''_{30} c(\psi)^{1/2} \rho_*(\zeta)^{5/4}\right)
 \varepsilon \|\FF\|_{L_2(\O)},
\end{aligned}
\eqno(8.51)
$$
\textit{or, in operator terms,}
$$
\begin{aligned}
&\| (\A_{D,\varepsilon}- \zeta I)^{-1}  - (\A^0_{D}- \zeta I)^{-1} - \eps K_D^0(\eps;\zeta)\|_{L_2({\mathcal O}) \to H^1({\mathcal O}')}
\cr
&\leq
\left({\mathcal C}'_{30}\delta^{-1} (c(\psi) \rho_*(\zeta) + c(\psi)^{5/2} \rho_*(\zeta)^{3/4}) +
\check{\mathcal C}''_{30} c(\psi)^{1/2} \rho_*(\zeta)^{5/4}\right) \varepsilon.
\end{aligned}
$$
\textit{For the flux $\p_\eps= g^\eps b(\D)\u_\eps$ we have}
$$
\begin{aligned}
&\| \p_\eps - \widetilde{g}^\eps b(\D) {\u}_0 \|_{L_2(\O')}
\cr
&\leq
\left(\widetilde{\mathcal C}'_{30}\delta^{-1} (c(\psi) \rho_*(\zeta) + c(\psi)^{5/2} \rho_*(\zeta)^{3/4}) +
\widehat{\mathcal C}''_{30} c(\psi)^{1/2} \rho_*(\zeta)^{5/4}\right)
\varepsilon \|\FF\|_{L_2(\O)}
\end{aligned}
\eqno(8.52)
$$
\textit{for $0< \eps \leq \eps_1$. The constants ${\mathcal C}'_{30}$ and $\widetilde{\mathcal C}'_{30}$ are the same as in Theorem} 8.6.
\textit{The constants $\check{\mathcal C}''_{30}$ and $\widehat{\mathcal C}''_{30}$
depend on} $d$, $m$, $\alpha_0$, $\alpha_1$, $\| g \|_{L_\infty}$, $\| g^{-1} \|_{L_\infty}$,
\textit{the parameters of the lattice} $\Gamma$, \textit{the domain $\O$, and} $\|\Lambda\|_{L_\infty}$.

\smallskip\noindent\textbf{Proof.}
Relations (8.27) and (8.30) imply (8.51) with
$\check{\mathcal C}''_{30} = {\mathcal C}''_{30} + {\mathcal C}''' C_{\O}^{(2)} \widehat{c} \check{c}^{1/2}$.

From (8.51) it follows that
$$
\begin{aligned}
&\| \p_\eps - g^\eps b(\D) \check{\v}_\eps \|_{L_2(\O')} \leq
\|g\|_{L_\infty} (d \alpha_1)^{1/2}
\cr
&\times \left({\mathcal C}'_{30}\delta^{-1} (c(\psi) \rho_*(\zeta) + c(\psi)^{5/2} \rho(\zeta)^{3/4}) +
\check{\mathcal C}''_{30} c(\psi)^{1/2} \rho_*(\zeta)^{5/4}\right)
 \varepsilon \|\FF\|_{L_2(\O)}.
\end{aligned}
$$
Together with (8.23) and (8.29) this yields (8.52) with
$\widehat{\mathcal C}''_{30} = \|g\|_{L_\infty} (d \alpha_1)^{1/2}  \check{\mathcal C}''_{30}  +
\widetilde{\mathcal C}' \widehat{c} \check{c}^{1/2}$. $\bullet$

\section*{Chapter 3. The Neumann problem}

\section*{\S 9. The Neumann problem in a bounded domain: Preliminaries}

\smallskip\noindent\textbf{9.1. Coercivity.} As in Chapter 2, we assume that $\O \subset \R^d$ is a bounded domain with the boundary
of class $C^{1,1}$. We impose an additional condition on the symbol $b(\bxi)$ of the operator (1.2).

\smallskip\noindent\textbf{Condition 9.1.} \textit{The matrix-valued function $b(\bxi)=\sum_{l=1}^d b_l \xi_l$ is such that}
$$
\rank b(\bxi) =n, \quad 0 \ne \bxi \in \C^d.
\eqno(9.1)
$$

\smallskip
Note that condition (9.1) is more restrictive than (1.3).
According to [Ne] (see Theorem 7.8 in \S 3.7), Condition 9.1 \textit{is necessary and sufficient for coercivity}
of the form $\|b(\D)\u\|^2_{L_2(\O)}$ on $H^1(\O;\C^n)$.

\smallskip\noindent\textbf{Proposition 9.2. [Ne]}
\textit{Condition} 9.1 \textit{is necessary and sufficient for existence of constants
$k_1 >0$ and $k_2 \ge 0$ such that the G\"arding type inequality}
$$
\|b(\D)\u \|_{L_2(\O)}^2 + k_2 \|\u\|^2_{L_2(\O)} \geq k_1 \|\D \u\|^2_{L_2(\O)},\quad
\u \in H^1(\O;\C^n),
\eqno(9.2)
$$
\textit{is satisfied.}

\smallskip\noindent\textbf{Remark 9.3.} The constants $k_1$ and $k_2$ depend on
the matrix $b(\bxi)$ and the domain $\O$, but in the general case it is difficult to control these constants explicitly.
However, often for particular operators they are known. Therefore, in what follows we indicate the dependence
of other constants on $k_1$ and $k_2$.

\smallskip

Below \textit{we assume that Condition} 9.1 \textit{is satisfied}.

\smallskip\noindent\textbf{9.2. The operator $\A_{N,\eps}$.}
In $L_2 (\mathcal{O}; \mathbb{C}^{n})$, consider the operator
$\mathcal{A}_{N,\varepsilon}$ formally given by the differential expression
$b (\mathbf{D})^* g^{\varepsilon} (\mathbf{x}) b (\mathbf{D})$
with the Neumann condition on $\partial \mathcal{O}$.
Precisely, $\A_{N,\eps}$ is the selfadjoint operator in $L_2 (\mathcal{O}; \mathbb{C}^{n})$ generated by the quadratic form
$$
a_{N,\eps}[\u,\u] := \int_{\mathcal{O}} \left\langle g^{\varepsilon} (\mathbf{x}) b (\mathbf{D}) \mathbf{u},
b(\mathbf{D}) \mathbf{u} \right\rangle \, d \mathbf{x},
\quad \mathbf{u} \in H^1 (\mathcal{O}; \mathbb{C}^{n}).
\eqno(9.3)
$$

By (1.2) and (1.5),
$$
a_{N,\eps}[\u,\u] \le d \alpha_1 \|g\|_{L_\infty} \|\D\u\|^2_{L_2(\O)},\quad \u \in H^1(\O;\C^n).
\eqno(9.4)
$$
From (9.2) it follows that
$$
a_{N,\eps}[\u,\u] \ge \|g^{-1}\|^{-1}_{L_\infty} \left(
 k_1 \|\D \u\|^2_{L_2(\O)} - k_2 \| \u\|^2_{L_2(\O)}\right),
\quad \u \in H^1(\O;\C^n).
\eqno(9.5)
$$
By (9.3)--(9.5), the form (9.3) is closed and nonnegative.

The spectrum of $\A_{N,\eps}$ is contained in $\R_+$. The point $\zeta \in \C \setminus \R_+$ is regular for
this operator.
\textit{Our goal} in Chapter 3 is to approximate
the generalized solution $\u_\eps \in H^1(\O;\C^n)$ of the Neumann problem
$$
b (\mathbf{D})^* g^{\varepsilon} (\mathbf{x}) b (\mathbf{D}) \mathbf{u}_{\varepsilon}(\x)
- \zeta \u_\eps(\x) = \mathbf{F}(\x),
\ \ \x \in \O;
\quad \partial^\eps_{\bnu} \u_\eps\vert_{\partial \O}=0,
\eqno(9.6)
$$
for small $\eps$. Here $\mathbf{F} \in L_2 (\mathcal{O}; \mathbb{C}^{n})$.
Then $\u_\eps = (\mathcal{A}_{N,\varepsilon}-\zeta I)^{-1} \FF$.
Here the notation $\partial^\eps_{\bnu}$ for the "`conormal derivative"' was used.
Let $\bnu(\x)$ be the unit outer normal vector to $\partial \O$
at the point $\x \in \partial \O$. Then formally the conormal derivative is given by
$\partial^\eps_{\bnu} \u(\x) := b(\bnu(\x))^* g^\eps(\x) b(\nabla) \u(\x)$.

The following lemma is an analog of Lemma 3.1.

\smallskip\noindent\textbf{Lemma 9.4.}
\textit{Let} $\zeta \in \C \setminus \R_+$, $|\zeta|\geq 1$. \textit{Let} $\u_\eps$ \textit{be the solution of problem} (9.6).
\textit{Then for $\eps>0$ we have}
$$
\|\u_\eps \|_{L_2(\O)} \leq c(\varphi) |\zeta|^{-1} \|\FF\|_{L_2(\O)},
\eqno(9.7)
$$
$$
\|\D \u_\eps \|_{L_2(\O)} \leq {\mathfrak C}_0 c(\varphi) |\zeta|^{-1/2} \|\FF\|_{L_2(\O)},
\eqno(9.8)
$$
\textit{or, in operator terms},
$$
\| (\A_{N,\eps} - \zeta I)^{-1} \|_{L_2(\O)\to L_2(\O)} \leq c(\varphi) |\zeta|^{-1},
\eqno(9.9)
$$
$$
\|\D (\A_{N,\eps} - \zeta I)^{-1} \|_{L_2(\O)\to L_2(\O)} \leq {\mathfrak C}_0 c(\varphi) |\zeta|^{-1/2}.
$$
\textit{The constant ${\mathfrak C}_0$ depends only on $\|g^{-1}\|_{L_\infty}$ and the constants $k_1$ and $k_2$ from} (9.2).

 \smallskip\noindent\textbf{Proof.}
Since the norm of the resolvent $(\A_{N,\eps} - \zeta I)^{-1}$ does not exceed the inverse distance from $\zeta$ to $\R_+$,
we obtain (9.9).

In order to check (9.8), we write down the integral identity for the solution $\u_\eps \in H^1(\O;\C^n)$ of problem (9.6):
$$
(g^\eps b(\D)\u_\eps, b(\D) {\boldsymbol \eta})_{L_2(\O)} - \zeta (\u_\eps, {\boldsymbol \eta})_{L_2(\O)}=
(\FF, {\boldsymbol \eta})_{L_2(\O)},\ \  {\boldsymbol \eta}\in H^1(\O;\C^n).
\eqno(9.10)
$$
Substituting  ${\boldsymbol \eta}=\u_\eps$ in (9.10) and taking (9.7) into account, we obtain
$$
a_{N,\eps}[\u_\eps, \u_\eps] \leq 2 c(\varphi)^2  |\zeta|^{-1} \|\FF \|^2_{L_2(\O)}.
$$
Combining this with (9.5) and (9.7), we arrive at (9.8) with ${\mathfrak C}_0 = (2 k_1^{-1} \|g^{-1}\|_{L_\infty}+ k_2 k_1^{-1})^{1/2}$. $\bullet$

\smallskip\noindent\textbf{9.3. The effective operator $\A^0_{N}$.}
In $L_2(\O;\C^n)$, consider the selfadjoint operator $\A_N^0$ generated by the quadratic form
$$
a^0_N[\u,\u]= \int_{\mathcal{O}} \left\langle g^0 b (\mathbf{D}) \mathbf{u}, b(\mathbf{D}) \u \right\rangle \,d \mathbf{x},
 \quad \u \in H^1 (\mathcal{O}; \mathbb{C}^{n}).
\eqno(9.11)
$$
Here $g^0$ is the effective matrix defined by (1.8).
Taking (1.15) and (9.2) into account, we see that the form (9.11) satisfies estimates of the form
(9.4) and (9.5) with the same constants.

Since $\partial \mathcal{O} \in C^{1,1}$, the operator $\A_N^0$ is given by the differential expression
$b(\D)^* g^0 b(\D)$ on the domain $\{ \u \in H^2(\O;\C^n):\ \partial_{\bnu}^0 \u\vert_{\partial \O} =0\}$,
where $\partial_{\bnu}^0$ is the conormal derivative corresponding to the operator $b(\D)^* g^0 b(\D)$,
i.~e., $\partial_{\bnu}^0 \u (\x) = b(\bnu(\x))^* g^0 b(\nabla) \u(\x)$.
We have
$$
\|(\A_N^0 +I)^{-1}\|_{L_2(\O) \to H^2(\O)} \leq c^\circ.
\eqno(9.12)
$$
Here the constant $c^\circ$ depends only on the constants $k_1$, $k_2$ from (9.2),
$\alpha_0$, $\alpha_1$, $\|g\|_{L_\infty}$, $\|g^{-1}\|_{L_\infty}$, and the domain $\O$.
To justify this fact, we refer to the theorems about regularity of solutions
for strongly elliptic systems (see, e.~g., [McL, Chapter 4]).

Let $\mathbf{u}_0$ be the solution of the problem
$$
b(\D)^* g^0 b(\D) \u_0(\x) - \zeta \u_0(\x)= \FF(\x),\ \ \x \in \O;
\ \ \partial_{\bnu}^0 \u_0\vert_{\partial \O} =0,
\eqno(9.13)
$$
where $\FF \in L_2(\O;\C^n)$. Then $\u_0 = (\mathcal{A}^0_N - \zeta I)^{-1} \FF$.

\smallskip\noindent\textbf{Lemma 9.5.}
\textit{Let $\zeta \in \C \setminus \R_+$, $|\zeta|\geq 1$. Let $\u_0$ be the solution of problem} (9.13).
\textit{Then we have}
$$
\|\u_0 \|_{L_2(\O)} \leq c(\varphi) |\zeta|^{-1} \|\FF\|_{L_2(\O)},
$$
$$
\|\D \u_0 \|_{L_2(\O)} \leq {\mathfrak C}_0 c(\varphi) |\zeta|^{-1/2} \|\FF\|_{L_2(\O)},
$$
$$
\|\u_0 \|_{H^2(\O)} \leq 2 c^\circ  c(\varphi) \|\FF\|_{L_2(\O)},
\eqno(9.14)
$$
\textit{or, in operator terms},
$$
\| (\A^0_{N} - \zeta I)^{-1} \|_{L_2(\O)\to L_2(\O)} \leq c(\varphi) |\zeta|^{-1},
\eqno(9.15)
$$
$$
\|\D (\A^0_{N} - \zeta I)^{-1} \|_{L_2(\O)\to L_2(\O)} \leq {\mathfrak C}_0 c(\varphi) |\zeta|^{-1/2},
\eqno(9.16)
$$
$$
\|(\A^0_{N} - \zeta I)^{-1} \|_{L_2(\O)\to H^2(\O)} \leq 2 c^\circ c(\varphi).
\eqno(9.17)
$$

 \smallskip\noindent\textbf{Proof.}
Estimates (9.15) and (9.16) can be checked similarly to the proof of Lemma 9.4.
Estimate (9.17) is a consequence of (9.12) and the inequality
$$
\| (\A_N^0+I) (\A^0_{N} - \zeta I)^{-1} \|_{L_2(\O)\to L_2(\O)} \leq
\sup_{x\geq 0} (x+1) |x-\zeta|^{-1} \leq 2 c(\varphi),
$$
which is valid for $\zeta \in \C \setminus \R_+$, $|\zeta|\geq 1$.
$\bullet$

\section*{\S 10. The results for the Neumann problem}

\smallskip\noindent\textbf{10.1. Approximation of the resolvent $(\A_{N,\eps}-\zeta I)^{-1}$ for $|\zeta|\geq 1$.}
Now we formulate our main results for the operator  $\A_{N,\eps}$.

\smallskip\noindent\textbf{Theorem 10.1.} \textit{Let $\zeta \in {\mathbb C}\setminus {\mathbb R}_+$ and $|\zeta|\ge 1$.
Let $\u_\eps$ be the solution of problem} (9.6), \textit{and let $\u_0$ be the solution of problem} (9.13) \textit{with $\FF \in L_2(\O;\C^n)$.
Suppose that the number $\varepsilon_1$ satisfies Condition} 4.1. \textit{Then for $0< \varepsilon \leq \varepsilon_1$ we have}
$$
\| \u_\eps - \u_0 \|_{L_2(\O)} \leq {\mathfrak C}_1 c(\varphi)^5 ( |\zeta|^{-1/2}\varepsilon + \varepsilon^2) \|\FF\|_{L_2(\O)},
\eqno(10.1)
$$
\textit{or, in operator terms,}
$$
\| (\A_{N,\varepsilon}- \zeta I)^{-1}  - (\A^0_{N}- \zeta I)^{-1} \|_{L_2({\mathcal O}) \to L_2({\mathcal O})}
\le {\mathfrak C}_1 c(\varphi)^5 ( |\zeta|^{-1/2}\varepsilon + \varepsilon^2).
$$
\textit{The constant ${\mathfrak C}_1$ depends on
$d$, $m$, $\alpha_0$, $\alpha_1$, $\| g \|_{L_\infty}$, $\| g^{-1} \|_{L_\infty}$, the constants $k_1$, $k_2$ from} (9.2),
\textit{the parameters of the lattice} $\Gamma$, \textit{and the domain} $\O$.

\smallskip
To approximate the solution in $H^1(\O;\C^n)$, we introduce the corrector similar to (4.4):
$$
K_N(\varepsilon;\zeta)
= R_{\mathcal{O}} [\Lambda^\varepsilon] S_\varepsilon b(\mathbf{D}) P_{\mathcal{O}} (\A_N^0 -\zeta I)^{-1}.
\eqno(10.2)
$$
The operator $K_N(\varepsilon;\zeta)$ is a continuous mapping of $L_2(\mathcal{O};\mathbb{C}^n)$ to $H^1(\mathcal{O};\mathbb{C}^n)$.

Let $\u_0$ be the solution of problem (9.13). Denote $\widetilde{\u}_0:= P_\O \u_0$.
Similarly to (4.5)--(4.7), we define a function
$$
\widetilde{\v}_\eps(\x) = \widetilde{\u}_0(\x) + \eps \Lambda^\eps(\x) (S_\eps b(\D) \widetilde{\u}_0)(\x)
$$
in $\R^d$ and put
$$
 \v_\eps :=\widetilde{\v}_\eps\vert_\O.
 \eqno(10.3)
$$
Then
$$
\v_\eps = (\A_N^0 - \zeta I)^{-1}\FF + \eps K_N(\eps;\zeta) \FF.
\eqno(10.4)
$$

\smallskip\noindent\textbf{Theorem 10.2.} \textit{Suppose that the assumptions of Theorem} 10.1
 \textit{are satisfied. Let $\v_\eps$ be defined by} (10.2), (10.4).
\textit{Then for $0< \varepsilon \leq \varepsilon_1$ we have}
$$
\| \u_\eps - \v_\eps \|_{H^1(\O)} \leq \left({\mathfrak C}_2 c(\varphi)^2  |\zeta|^{-1/4}\varepsilon^{1/2} +
{\mathfrak C}_3 c(\varphi)^4 \varepsilon\right) \|\FF\|_{L_2(\O)},
\eqno(10.5)
$$
\textit{or, in operator terms},
$$
\begin{aligned}
&\| (\A_{N,\varepsilon}- \zeta I)^{-1}  - (\A^0_{N}- \zeta I)^{-1} - \eps K_N(\eps;\zeta) \|_{L_2({\mathcal O}) \to H^1({\mathcal O})}
\cr
&\leq  {\mathfrak C}_2 c(\varphi)^2 |\zeta|^{-1/4} \varepsilon^{1/2} + {\mathfrak C}_3 c(\varphi)^4 \varepsilon.
\end{aligned}
$$
\textit{For the flux} $\p_\eps := g^\eps b(\D) \u_\eps$ \textit{we have}
$$
\| \p_\eps - \widetilde{g}^\eps S_\eps b(\D) \widetilde{\u}_0 \|_{L_2(\O)}
\leq \left(\widetilde{\mathfrak C}_2 c(\varphi)^2 |\zeta|^{-1/4}\varepsilon^{1/2} +
\widetilde{\mathfrak C}_3 c(\varphi)^4 \varepsilon\right) \|\FF\|_{L_2(\O)}
\eqno(10.6)
$$
\textit{for $0< \eps \leq \eps_1$. The constants ${\mathfrak C}_2$, ${\mathfrak C}_3$, $\widetilde{\mathfrak C}_2$, and
$\widetilde{\mathfrak C}_3$ depend on}
$d$, $m$, $\alpha_0$, $\alpha_1$, $\| g \|_{L_\infty}$, $\| g^{-1} \|_{L_\infty}$,
\textit{the constants $k_1$, $k_2$ from} (9.2),
\textit{the parameters of the lattice} $\Gamma$, \textit{and the domain} $\O$.

\smallskip\noindent\textbf{10.2. The first step of the proof. Associated problem in $\R^d$.}
As in Subsection 4.2, we start with the associated problem in $\R^d$.
By Lemma 9.5 and (4.3), we have
$$
\| \widetilde{\u}_0 \|_{L_2(\R^d)} \leq C_\O^{(0)} c(\varphi) |\zeta|^{-1} \|\FF\|_{L_2(\O)},
\eqno(10.7)
$$
$$
\| \widetilde{\u}_0 \|_{H^1(\R^d)} \leq C_\O^{(1)} ({\mathfrak C}_0+1) c(\varphi) |\zeta|^{-1/2} \|\FF\|_{L_2(\O)},
\eqno(10.8)
$$
$$
\| \widetilde{\u}_0 \|_{H^2(\R^d)} \leq 2 C_\O^{(2)} c^\circ c(\varphi) \|\FF\|_{L_2(\O)}.
\eqno(10.9)
$$
We put
$$
\widetilde{\FF}:= \A^0 \widetilde{\u}_0 - \zeta \widetilde{\u}_0.
\eqno(10.10)
$$
Then $\widetilde{\FF}\in L_2(\R^d;\C^n)$ and $\widetilde{\FF}\vert_{\mathcal O} =\FF$.
Similarly to  (4.15), from (1.16), (10.7), and (10.9) it follows that
$$
\| \widetilde{\FF}\|_{L_2(\R^d)} \leq  {\mathfrak C}_4 c(\varphi) \|\FF\|_{L_2(\O)},
\eqno(10.11)
$$
where ${\mathfrak C}_4 = 2 c_1 C_\O^{(2)} c^\circ + C_\O^{(0)}$.

Let $\widetilde{\u}_\eps \in H^1(\R^d;\C^n)$ be the solution of the following equation in $\R^d$:
$$
\A_\eps \widetilde{\u}_\eps - \zeta \widetilde{\u}_\eps = \widetilde{\FF},
\eqno(10.12)
$$
i.~e., $\widetilde{\u}_\eps = (\A_\eps - \zeta I)^{-1} \widetilde{\FF}$.
We can apply theorems from  \S 2.
By Theorem 2.2 and (10.10)--(10.12), we see that
$$
\| \widetilde{\u}_\eps - \widetilde{\u}_0\|_{L_2(\R^d)} \leq
{\mathfrak C}_4 C_1 c(\varphi)^3 |\zeta|^{-1/2} \eps \|{\FF}\|_{L_2(\O)},\ \ \eps>0.
\eqno(10.13)
$$
From Theorem 2.4 and relations (10.10)--(10.12) it follows that
$$
\| \D (\widetilde{\u}_\eps - \widetilde{\v}_\eps) \|_{L_2(\R^d)} \leq
{\mathfrak C}_4 C_2 c(\varphi)^3 \eps \|{\FF}\|_{L_2(\O)},\ \ \eps >0,
\eqno(10.14)
$$
$$
\| \widetilde{\u}_\eps - \widetilde{\v}_\eps \|_{L_2(\R^d)} \leq
{\mathfrak C}_4 C_3 c(\varphi)^3 |\zeta|^{-1/2} \eps \|{\FF}\|_{L_2(\O)},\ \ \eps >0.
\eqno(10.15)
$$
Finally, by Theorem 2.6 and (10.10)--(10.12), we have
$$
\| g^\eps b(\D) \widetilde{\u}_\eps - \widetilde{g}^\eps S_\eps b(\D) \widetilde{\u}_0 \|_{L_2(\R^d)}
\leq
{\mathfrak C}_4 C_4  c(\varphi)^3 \eps  \|{\FF}\|_{L_2(\O)},\ \ \eps>0.
\eqno(10.16)
$$

\smallskip\noindent\textbf{10.3. The second step of the proof. Introduction of the correction term $\w_\eps$.}
Now we introduce the  "`correction term"' $\w_\eps \in H^1(\O;\C^n)$ as a function satisfying the
integral identity
$$
\begin{aligned}
&(g^\eps b(\D) \w_\eps,b(\D) {\boldsymbol \eta})_{L_2(\O)} -
\zeta (\w_\eps, {\boldsymbol \eta})_{L_2(\O)}
\cr
&=
(\widetilde{g}^\eps S_\eps b(\D) \widetilde{\u}_0, b(\D) {\boldsymbol \eta})_{L_2(\O)} - (\zeta \u_0 + \FF, {\boldsymbol \eta})_{L_2(\O)},
\ \ {\boldsymbol \eta} \in H^1(\O;\C^n).
\end{aligned}
\eqno(10.17)
$$
Since the right-hand side of  (10.17) is an antilinear continuous functional of ${\boldsymbol \eta} \in H^1(\O;\C^n)$,
one can check in a standard way that the solution $\w_\eps$ exists and is unique.

\smallskip\noindent\textbf{Lemma 10.3.} \textit{Let $\zeta \in \C \setminus \R_+$ and $|\zeta|\geq 1$.
Let $\u_\eps$ be the solution of problem} (9.6), \textit{and let $\widetilde{\u}_\eps$ be the solution of equation} (10.12).
\textit{Suppose that $\w_\eps$ satisfies} (10.17). \textit{Then for $\eps>0$ we have}
$$
\| \D (\u_\eps - \widetilde{\u}_\eps + \w_\eps) \|_{L_2({\mathcal O})} \leq {\mathfrak C}_5 c(\varphi)^4  \eps \|\FF\|_{L_2({\mathcal O})},
\eqno(10.18)
$$
$$
\| \u_\eps - \widetilde{\u}_\eps + \w_\eps \|_{L_2({\mathcal O})} \leq {\mathfrak C}_6 c(\varphi)^4 |\zeta|^{-1/2} \eps \|\FF\|_{L_2({\mathcal O})}.
\eqno(10.19)
$$
\textit{The constants ${\mathfrak C}_5$ and ${\mathfrak C}_6$ depend on} $d$, $m$, $\alpha_0$, $\alpha_1$, $\| g \|_{L_\infty}$,
$\| g^{-1} \|_{L_\infty}$,
\textit{the constants $k_1$, $k_2$ from} (9.2), \textit{the parameters of the lattice} $\Gamma$, \textit{and the domain} $\O$.

\smallskip\noindent\textbf{Proof.} Denote ${\mathbf U}_\eps:= \u_\eps - \widetilde{\u}_\eps  + \w_\eps$.
By (9.10) and (10.17), the function ${\mathbf U}_\eps \in H^1(\O;\C^n)$ satisfies the identity
$$
\begin{aligned}
&(g^\eps b(\D) {\mathbf U}_\eps, b(\D){\boldsymbol \eta})_{L_2(\O)} - \zeta({\mathbf U}_\eps, {\boldsymbol \eta})_{L_2(\O)}
\cr
&=
- (g^\eps b(\D) \widetilde{\u}_\eps - \widetilde{g}^\eps S_\eps b(\D) \widetilde{\u}_0, b(\D){\boldsymbol \eta})_{L_2(\O)}
+ \zeta( \widetilde{\u}_\eps - \u_0, {\boldsymbol \eta})_{L_2(\O)},
\cr
&\ \ {\boldsymbol \eta}\in H^1(\O;\C^n).
\end{aligned}
\eqno(10.20)
$$
From (10.13) and (10.16) it follows that the right-hand side in (10.20) does not exceed
${\mathfrak C}_4 c(\varphi)^3 \eps \left(C_4 \|b(\D){\boldsymbol \eta}\|_{L_2(\O)}
+ C_1 |\zeta|^{1/2} \| {\boldsymbol \eta}\|_{L_2(\O)} \right) \|\FF\|_{L_2(\O)}$.

Substituting ${\boldsymbol \eta} = {\mathbf U}_\eps$ in (10.20), similarly to (4.26)--(4.29) we deduce
$$
\begin{aligned}
\| {\mathbf U}_\eps \|^2_{L_2(\O)} &\leq 4 {\mathfrak C}_4 C_4 c(\varphi)^4 |\zeta|^{-1} \eps \| b(\D){\mathbf U}_\eps \|_{L_2(\O)} \|\FF\|_{L_2(\O)}
\cr
&+ 4 {\mathfrak C}_4^2 C_1^2 c(\varphi)^8 |\zeta|^{-1} \eps^2 \|\FF\|^2_{L_2(\O)}.
\end{aligned}
\eqno(10.21)
$$
Now, combining (10.20) with ${\boldsymbol \eta} = {\mathbf U}_\eps$ and (10.21), we obtain
$$
\| b(\D) {\mathbf U}_\eps \|_{L_2(\O)}
\leq  \widetilde{\mathfrak C}_5 c(\varphi)^4  \eps \|\FF\|_{L_2(\O)},
\eqno(10.22)
$$
where $\widetilde{\mathfrak C}_5^2 = 9{\mathfrak C}_4^2(2 C_1^2 \|g^{-1}\|_{L_\infty} + 9 C_4^2 \|g^{-1}\|^2_{L_\infty})$.
Relations (10.21) and (10.22) yield (10.19) with
${\mathfrak C}_6 = 2 ({\mathfrak C}_4 \widetilde{\mathfrak C}_5 C_4 + {\mathfrak C}_4^2 C_1^2)^{1/2}$.
Finally, (9.2), (10.19), and (10.22) imply (10.18) with
${\mathfrak C}_5^2 = k_1^{-1}(\widetilde{\mathfrak C}_5^2 + k_2 {\mathfrak C}_6^2)$. $\bullet$

\smallskip
\textbf{Conclusions.} 1) From (10.3), (10.14), (10.15), (10.18), and (10.19) it follows that
$$
\|\D(\u_\eps - \v_\eps)\|_{L_2(\O)} \leq
 {\mathfrak C}_7 c(\varphi)^4 \eps \|\FF\|_{L_2(\O)} + \|\D \w_\eps \|_{L_2(\O)},\ \ \eps >0,
\eqno(10.23)
$$
$$
\| \u_\eps - \v_\eps \|_{L_2(\O)} \leq
{\mathfrak C}_{8} c(\varphi)^4 |\zeta|^{-1/2} \eps \|\FF\|_{L_2(\O)} + \|\w_\eps \|_{L_2(\O)},\ \ \eps >0,
\eqno(10.24)
$$
where ${\mathfrak C}_7 = {\mathfrak C}_4 C_2 + {\mathfrak C}_5$ and ${\mathfrak C}_{8} = {\mathfrak C}_4 C_3 + {\mathfrak C}_6$.
Thus, in order to prove Theorem 10.2, we need to obtain an appropriate estimate for $\|\w_\eps\|_{H^1(\O)}$.

2) From (10.13) and (10.19) it follows that
$$
\| \u_\eps - \u_0 \|_{L_2(\O)} \leq
{\mathfrak C}_{9} c(\varphi)^4 |\zeta|^{-1/2} \eps \|\FF\|_{L_2(\O)} + \|\w_\eps \|_{L_2(\O)},\ \ \eps >0,
\eqno(10.25)
$$
where ${\mathfrak C}_{9} = {\mathfrak C}_6 + {\mathfrak C}_4 C_1$.
Hence, for the proof of Theorem 10.1 one has to find an appropriate estimate for $\|\w_\eps\|_{L_2(\O)}$.

\section*{\S 11. Estimates of the correction term. \\ Proof of Theorems 10.1 and 10.2}

As in \S 5, first we estimate the $H^1$-norm of  $\w_\eps$ and prove Theorem 10.2.
Next, using the already proved Theorem  10.2 and the duality arguments, we estimate the
$L_2$-norm of $\w_\eps$ and prove Theorem 10.1.

\smallskip\noindent\textbf{11.1. Estimate for the norm of $\w_\eps$ in $H^1(\O;\C^n)$.}
Denote
$$
{\mathcal I}_\eps[{\boldsymbol \eta}] :=
(\widetilde{g}^\eps S_\eps b(\D) \widetilde{\u}_0 - g^0 b(\D)\u_0, b(\D) {\boldsymbol \eta})_{L_2(\O)},\ \ {\boldsymbol \eta}\in H^1(\O;\C^n).
\eqno(11.1)
$$
Note that the solution $\u_0$ of problem (9.13) satisfies the identity
$$
(g^0 b(\D)\u_0, b(\D) {\boldsymbol \eta})_{L_2(\O)} - (\zeta \u_0 +\FF,{\boldsymbol \eta})_{L_2(\O)}=0,\ \ {\boldsymbol \eta}\in H^1(\O;\C^n).
\eqno(11.2)
$$
By (10.17) and (11.2), $\w_\eps$ satisfies the identity
$$
(g^\eps b(\D)\w_\eps, b(\D) {\boldsymbol \eta})_{L_2(\O)} - \zeta (\w_\eps, {\boldsymbol \eta})_{L_2(\O)}
={\mathcal I}_\eps[{\boldsymbol \eta}],\ \ {\boldsymbol \eta} \in H^1(\O;\C^n).
\eqno(11.3)
$$

\smallskip\noindent\textbf{Lemma 11.1.} \textit{Let $\zeta \in \C \setminus \R_+$ and $|\zeta|\geq 1$.
Suppose that the number $\eps_1$ satisfies Condition} 4.1.
\textit{Then for $0< \eps \leq \eps_1$ the functional} (11.1) \textit{satisfies the following estimate}:
$$
| {\mathcal I}_\eps[{\boldsymbol \eta}] | \leq c(\varphi) \left( {\mathfrak C}_{10}|\zeta|^{-1/4} \eps^{1/2} + {\mathfrak C}_{11}\eps\right)\|\FF\|_{L_2(\O)} \| \D{\boldsymbol \eta}\|_{L_2(\O)}, \ \ {\boldsymbol \eta} \in H^1(\O;\C^n).
\eqno(11.4)
$$
\textit{The constants ${\mathfrak C}_{10}$ and ${\mathfrak C}_{11}$
depend on} $d$, $m$, $\alpha_0$, $\alpha_1$, $\| g \|_{L_\infty}$, $\| g^{-1} \|_{L_\infty}$,
\textit{the constants $k_1$, $k_2$ from} (9.2), \textit{the parameters of the lattice} $\Gamma$, \textit{and the domain} $\O$.

\smallskip\noindent\textbf{Proof.}
The functional (11.1) can be represented as
$$
{\mathcal I}_\eps[{\boldsymbol \eta}] = {\mathcal I}_\eps^{(1)}[{\boldsymbol \eta}] + {\mathcal I}_\eps^{(2)}[{\boldsymbol \eta}],
\eqno(11.5)
$$
where
$$
{\mathcal I}_\eps^{(1)}[{\boldsymbol \eta}]=
({g}^0 S_\eps b(\D) \widetilde{\u}_0 - g^0 b(\D)\u_0, b(\D) {\boldsymbol \eta})_{L_2(\O)},
\eqno(11.6)
$$
$$
{\mathcal I}_\eps^{(2)}[{\boldsymbol \eta}]=
((\widetilde{g}^\eps - {g}^0) S_\eps b(\D) \widetilde{\u}_0, b(\D) {\boldsymbol \eta})_{L_2(\O)}.
\eqno(11.7)
$$

The term (11.6) is estimated with the help of Proposition 1.4 and relations (1.2), (1.4), (1.5), (1.15), and (10.9):
$$
\begin{aligned}
&|{\mathcal I}_\eps^{(1)}[{\boldsymbol \eta}]| \leq
|{g}^0| \| (S_\eps -I)b(\D) \widetilde{\u}_0\|_{L_2(\R^d)} \| b(\D) {\boldsymbol \eta}\|_{L_2(\O)}
\cr
&\leq {\mathfrak C}^{(1)} c(\varphi) \eps \|\FF\|_{L_2(\O)} \| \D {\boldsymbol \eta}\|_{L_2(\O)},\ \ \eps >0,
\end{aligned}
\eqno(11.8)
$$
where ${\mathfrak C}^{(1)} = 2 \|g\|_{L_\infty} r_1 \alpha_1 d^{1/2} C_\O^{(2)} c^\circ$.

Using (1.2), we transform the term (11.7):
$$
{\mathcal I}_\eps^{(2)}[{\boldsymbol \eta}]=
\sum_{l=1}^d ( f_l^\eps S_\eps b(\D) \widetilde{\u}_0, D_l {\boldsymbol \eta})_{L_2(\O)},
\eqno(11.9)
$$
where $f_l(\x):= b_l^* (\widetilde{g}(\x) - {g}^0)$, $l=1,\dots,d$. According to [Su3, (5.13)], we have
$$
\begin{aligned}
&\|f_l\|_{L_2(\Omega)}\leq |\Omega|^{1/2} \check{\mathfrak C},\ \ l=1,\dots,d,
\cr
& \check{\mathfrak C}= \alpha_1^{1/2} \|g\|_{L_\infty} (1+ (d m)^{1/2} \alpha_1^{1/2}\alpha_0^{-1/2} \|g\|_{L_\infty}^{1/2} \|g^{-1}\|_{L_\infty}^{1/2}).
\end{aligned}
\eqno(11.10)
$$
As was checked in  [Su3, Subsection 5.2], there exist $\Gamma$-periodic $(n\times m)$-matrix-valued functions $M_{lj}(\x)$
in $\R^d$, $l,j=1,\dots,d$, such that
$$
\begin{aligned}
&M_{lj}\in \widetilde{H}^1(\Omega),\ \ \int_\Omega M_{lj}(\x)\,d\x=0,\ \ M_{lj}(\x)=-M_{jl}(\x),\ \ l,j=1,\dots,d,
\cr
&f_l(\x) = \sum_{j=1}^d \partial_j M_{lj}(\x),\ \ l=1,\dots,d,
\end{aligned}
\eqno(11.11)
$$
and
$$
\|M_{lj} \|_{L_2(\Omega)}\leq r_0^{-1} |\Omega|^{1/2} \check{\mathfrak C},\ \ l,j =1,\dots,d.
\eqno(11.12)
$$
By (11.11),
$f_l^\eps(\x) = \eps \sum_{j=1}^d \partial_j  M^\eps_{lj}(\x)$, $l=1,\dots,d$, whence the term (11.9) can be written as
$$
{\mathcal I}_\eps^{(2)}[{\boldsymbol \eta}]=
\widetilde{\mathcal I}_\eps^{(2)}[{\boldsymbol \eta}]+
\widehat{\mathcal I}_\eps^{(2)}[{\boldsymbol \eta}],
\eqno(11.13)
$$
where
$$
\widetilde{\mathcal I}_\eps^{(2)}[{\boldsymbol \eta}]
= \eps \sum_{l,j=1}^d \left( \partial_j(M_{lj}^\eps S_\eps b(\D) \widetilde{\u}_0), D_l {\boldsymbol \eta}\right)_{L_2(\O)},
\eqno(11.14)
$$
$$
\widehat{\mathcal I}_\eps^{(2)}[{\boldsymbol \eta}]
= - \eps \sum_{l,j=1}^d \left( M_{lj}^\eps S_\eps b(\D) \partial_j \widetilde{\u}_0, D_l {\boldsymbol \eta}\right)_{L_2(\O)}.
\eqno(11.15)
$$
The term (11.15) is estimated by using Proposition 1.5:
$$
| \widehat{\mathcal I}_\eps^{(2)}[{\boldsymbol \eta}] | \leq \eps \sum_{l,j=1}^d |\Omega|^{-1/2} \| M_{lj}\|_{L_2(\Omega)}
\| b(\D) \partial_j \widetilde{\u}_0 \|_{L_2(\R^d)}  \| D_l {\boldsymbol \eta} \|_{L_2(\O)}.
$$
Combining this with (1.4), (10.9), and (11.12), we obtain
$$
| \widehat{\mathcal I}_\eps^{(2)}[{\boldsymbol \eta}] | \leq
{\mathfrak C}^{(2)} c(\varphi) \eps
\| \FF \|_{L_2(\O)}  \| \D {\boldsymbol \eta} \|_{L_2(\O)},\ \ \eps>0,
\eqno(11.16)
$$
where ${\mathfrak C}^{(2)}= 2 d r_0^{-1} \check{\mathfrak C} \alpha_1^{1/2} C_\O^{(2)} c^\circ$.

Now consider the term (11.14). Assume that  $0< \eps \leq \eps_1$.
Let $\theta_\eps$ be a cut-off function satisfying (5.1). We have 
$$
\sum_{l,j=1}^d \left( \partial_j( (1-\theta_\eps)M_{lj}^\eps S_\eps b(\D) \widetilde{\u}_0), D_l {\boldsymbol \eta}\right)_{L_2(\O)} =0,
\ \ {\boldsymbol \eta} \in H^1(\O;\C^n),
$$
which can be checked by integration by parts and using (11.11). Consequently, the term (11.14) can be written as
$$
\widetilde{\mathcal I}_\eps^{(2)}[{\boldsymbol \eta}]
= \eps \sum_{l,j=1}^d \left( \partial_j( \theta_\eps M_{lj}^\eps S_\eps b(\D) \widetilde{\u}_0), D_l {\boldsymbol \eta}\right)_{L_2(\O)}.
\eqno(11.17)
$$
Consider the following expression
$$
\begin{aligned}
&\eps \sum_{j=1}^d \partial_j( \theta_\eps M_{lj}^\eps S_\eps b(\D) \widetilde{\u}_0) =
\eps \theta_\eps \sum_{j=1}^d  M_{lj}^\eps S_\eps (b(\D) \partial_j \widetilde{\u}_0)
\cr
&+ \eps \sum_{j=1}^d (\partial_j \theta_\eps) M_{lj}^\eps S_\eps b(\D) \widetilde{\u}_0 +
\theta_\eps f^\eps_l  S_\eps b(\D) \widetilde{\u}_0.
\end{aligned}
$$
We have
$$
\eps \left\|\sum_{j=1}^d \partial_j( \theta_\eps M_{lj}^\eps S_\eps b(\D) \widetilde{\u}_0) \right\|_{L_2(\O)}
\leq J_l^{(1)}(\eps) + J_l^{(2)}(\eps) + J_l^{(3)}(\eps),
\eqno(11.18)
$$
where
$$
J_l^{(1)}(\eps) =
\eps \sum_{j=1}^d \| \theta_\eps M_{lj}^\eps S_\eps (b(\D) \partial_j \widetilde{\u}_0)\|_{L_2(\R^d)},
\eqno(11.19)
$$
$$
J_l^{(2)}(\eps) =
\eps \sum_{j=1}^d \| (\partial_j \theta_\eps) M_{lj}^\eps S_\eps b(\D) \widetilde{\u}_0\|_{L_2(\R^d)},
\eqno(11.20)
$$
$$
J_l^{(3)}(\eps) =
 \| \theta_\eps f_l^\eps S_\eps b(\D) \widetilde{\u}_0\|_{L_2(\R^d)}.
\eqno(11.21)
$$
To estimate the term (11.19), we apply (5.1), Proposition 1.5, and (1.4), (10.9), (11.12):
$$
J_l^{(1)}(\eps) \leq
\eps \sum_{j=1}^d |\Omega|^{-1/2} \| M_{lj}\|_{L_2(\Omega)} \| b(\D) \partial_j \widetilde{\u}_0\|_{L_2(\R^d)}
\leq \nu_1 c(\varphi) \eps \|\FF\|_{L_2(\O)},
\eqno(11.22)
$$
where $\nu_1 = 2 r_0^{-1} \check{\mathfrak C} (d \alpha_1)^{1/2} C_\O^{(2)} c^\circ$.
The term (11.20) is estimated with the help of  (5.1) and Lemma 3.4:
$$
\begin{aligned}
&(J_l^{(2)}(\eps))^2 \leq  d \kappa^2 \sum_{j=1}^d \int_{(\partial \O)_\eps} |M_{lj}^\eps S_\eps b(\D) \widetilde{\u}_0|^2 \,d\x
\cr
&\leq \eps d \kappa^2 \beta_* |\Omega|^{-1}  \sum_{j=1}^d \|M_{lj}\|^2_{L_2(\Omega)}
\|b(\D) \widetilde{\u}_0\|_{H^1(\R^d)} \| b(\D) \widetilde{\u}_0\|_{L_2(\R^d)}
\end{aligned}
$$
for $0< \eps \leq \eps_1$. Combining this with relations (1.4), (10.8), (10.9), and (11.12), we obtain
$$
J_l^{(2)}(\eps) \leq \nu_2 c(\varphi) |\zeta|^{-1/4} \eps^{1/2} \|\FF \|_{L_2(\O)},\ \ 0< \eps \leq \eps_1,
\eqno(11.23)
$$
where $\nu_2 = d \kappa r_0^{-1} \check{\mathfrak C} \left(2 \beta_* \alpha_1 C_\O^{(2)} C_\O^{(1)}
c^\circ ({\mathfrak C}_0+1)\right)^{1/2}$.
Similarly, using (5.1) and Lemma 3.4, we see that the term (11.21) satisfies
$$
\begin{aligned}
&(J_l^{(3)}(\eps))^2 \leq   \int_{(\partial \O)_\eps} |f_{l}^\eps S_\eps b(\D) \widetilde{\u}_0|^2 \,d\x
\cr
&\leq \eps \beta_* |\Omega|^{-1}  \|f_{l}\|^2_{L_2(\Omega)}
\|b(\D) \widetilde{\u}_0\|_{H^1(\R^d)} \| b(\D) \widetilde{\u}_0\|_{L_2(\R^d)}
\end{aligned}
$$
for $0< \eps \leq \eps_1$. Taking (1.4), (10.8), (10.9), and (11.10) into account, we deduce that
$$
J_l^{(3)}(\eps) \leq \nu_3 c(\varphi) |\zeta|^{-1/4} \eps^{1/2} \|\FF \|_{L_2(\O)},\ \ 0< \eps \leq \eps_1,
\eqno(11.24)
$$
where $\nu_3 = \check{\mathfrak C} \left(2  \beta_* \alpha_1 C_\O^{(2)} C_\O^{(1)}
c^\circ ({\mathfrak C}_0+1)\right)^{1/2}$.

Now relations (11.18) and (11.22)--(11.24) lead to the inequality
$$
\begin{aligned}
&\eps \left\|\sum_{j=1}^d \partial_j( \theta_\eps M_{lj}^\eps S_\eps b(\D) \widetilde{\u}_0) \right\|_{L_2(\O)}
\cr
&\leq c(\varphi) \left(\nu_1 \eps + (\nu_2+\nu_3) |\zeta|^{-1/4} \eps^{1/2}\right) \|\FF\|_{L_2(\O)},\ \ 0< \eps \leq \eps_1.
\end{aligned}
\eqno(11.25)
$$
Combining (11.17) and (11.25), we arrive at
$$
|\widetilde{\mathcal I}_\eps^{(2)}[{\boldsymbol \eta}] |\leq
d^{1/2} c(\varphi) \left(\nu_1 \eps + (\nu_2+\nu_3) |\zeta|^{-1/4} \eps^{1/2}\right)
\|\FF\|_{L_2(\O)} \|\D{\boldsymbol \eta}\|_{L_2(\O)}
\eqno(11.26)
$$
for $\ 0< \eps \leq \eps_1$.
As a result, relations (11.5), (11.8), (11.13), (11.16), and (11.26) imply the required estimate (11.4) with
${\mathfrak C}_{10}= d^{1/2}(\nu_2+\nu_3)$ and
${\mathfrak C}_{11}= {\mathfrak C}^{(1)} + {\mathfrak C}^{(2)} + d^{1/2} \nu_1$. $\bullet$

\smallskip\noindent\textbf{Lemma 11.2.} \textit{Let $\zeta \in \C \setminus \R_+$ and $|\zeta| \geq 1$.
Suppose that $\w_\eps$ satisfies} (10.17). \textit{Suppose that the number $\eps_1$ satisfies Condition} 4.1.
\textit{Then for $0< \eps \leq \eps_1$ we have}
$$
\| \w_\eps \|_{H^1(\O)} \leq c(\varphi)^2 \left( {\mathfrak C}_{12}|\zeta|^{-1/4} \eps^{1/2} +
{\mathfrak C}_{13}\eps \right)\|\FF\|_{L_2(\O)}.
\eqno(11.27)
$$
\textit{The constants ${\mathfrak C}_{12}$ and ${\mathfrak C}_{13}$ depend on}
$d$, $m$, $\alpha_0$, $\alpha_1$, $\| g \|_{L_\infty}$, $\| g^{-1} \|_{L_\infty}$,
\textit{the constants $k_1$, $k_2$ from} (9.2), \textit{the parameters of the lattice} $\Gamma$, \textit{and the domain} $\O$.

\smallskip\noindent\textbf{Proof.}
We substitute ${\boldsymbol \eta} =\w_\eps$ in (11.3) and take the imaginary part of the corresponding relation.
Taking (11.4) into account, for $0< \eps \leq \eps_1$ we have
$$
| {\rm Im}\, \zeta| \| \w_\eps \|^2_{L_2(\O)} \leq c(\varphi)  \left(
{\mathfrak C}_{10} |\zeta|^{-1/4} \eps^{1/2} + {\mathfrak C}_{11} \eps \right)
\|\FF\|_{L_2(\O)} \| \D \w_\eps \|_{L_2(\O)}.
\eqno(11.28)
$$
If ${\rm Re}\,\zeta \geq 0$ (and then ${\rm Im}\,\zeta \ne 0$), we deduce
$$
 \| \w_\eps \|^2_{L_2(\O)} \leq c(\varphi)^2   |\zeta|^{-1} \left( {\mathfrak C}_{10} |\zeta|^{-1/4} \eps^{1/2} + {\mathfrak C}_{11} \eps \right)
\|\FF\|_{L_2(\O)} \| \D \w_\eps \|_{L_2(\O)}.
\eqno(11.29)
$$
If ${\rm Re}\,\zeta < 0$, we take the real part of the corresponding relation and obtain
$$
| {\rm Re}\, \zeta| \| \w_\eps \|^2_{L_2(\O)} \leq c(\varphi)
\left( {\mathfrak C}_{10} |\zeta|^{-1/4} \eps^{1/2} + {\mathfrak C}_{11} \eps \right)
\|\FF\|_{L_2(\O)} \| \D \w_\eps \|_{L_2(\O)}
\eqno(11.30)
$$
for $0 < \eps \leq \eps_1$.
Summing up (11.28) and (11.30), we deduce the inequality similar to (11.29).
As a result, for all values of $\zeta$ under consideration we have
$$
\| \w_\eps \|^2_{L_2(\O)} \leq 2 c(\varphi)^2   |\zeta|^{-1} \left( {\mathfrak C}_{10} |\zeta|^{-1/4} \eps^{1/2} + {\mathfrak C}_{11} \eps \right)
\|\FF\|_{L_2(\O)} \| \D \w_\eps \|_{L_2(\O)}
\eqno(11.31)
$$
for $0< \eps\leq \eps_1$. Now from (11.3) with ${\boldsymbol \eta} = \w_\eps$, (11.4), and (11.31) it follows that
$$
a_{N,\eps}[\w_\eps, \w_\eps] \leq 3 c(\varphi)^2   \left( {\mathfrak C}_{10} |\zeta|^{-1/4} \eps^{1/2} +
{\mathfrak C}_{11} \eps \right)
\|\FF\|_{L_2(\O)} \| \D \w_\eps \|_{L_2(\O)}
$$
for $0< \eps\leq \eps_1$. Together with (9.5) and (11.31) this implies
$$
\| \D \w_\eps \|_{L_2(\O)} \leq  c(\varphi)^2  {\mathfrak C}_{14}
\left( {\mathfrak C}_{10} |\zeta|^{-1/4} \eps^{1/2} + {\mathfrak C}_{11} \eps \right)
\|\FF\|_{L_2(\O)},\ \ 0< \eps\leq \eps_1,
\eqno(11.32)
$$
where ${\mathfrak C}_{14} = k_1^{-1} (3 \|g^{-1}\|_{L_\infty} + 2 k_2)$.
Comparing (11.31) and (11.32), we arrive at (11.27) with
${\mathfrak C}_{12} = (2 {\mathfrak C}_{14} + {\mathfrak C}_{14}^2)^{1/2} {\mathfrak C}_{10}$ and
${\mathfrak C}_{13} = (2 {\mathfrak C}_{14} + {\mathfrak C}_{14}^2)^{1/2}{\mathfrak C}_{11}$. $\bullet$

\smallskip\noindent\textbf{11.2. Completion of the proof of Theorem 10.2.}
Relations (10.23), (10.24), and (11.27) imply the required estimate (10.5) with
${\mathfrak C}_2 = \sqrt{2} {\mathfrak C}_{12}$
and ${\mathfrak C}_3 = {\mathfrak C}_7 + {\mathfrak C}_{8} + \sqrt{2} {\mathfrak C}_{13}$.

It remains to check (10.6). From (10.5) and (1.2), (1.5) it follows that
$$
\begin{aligned}
&\| \p_\eps - g^\eps b(\D) \v_\eps \|_{L_2(\O)}
\cr
&\leq \|g\|_{L_\infty} (d \alpha_1)^{1/2}
({\mathfrak C}_2 c(\varphi)^2 |\zeta|^{-1/4} \eps^{1/2}  + {\mathfrak C}_3 c(\varphi)^4 \eps) \|\FF\|_{L_2(\O)}
\end{aligned}
\eqno(11.33)
$$
for $0< \eps\leq \eps_1$.
Similarly to (5.19)--(5.22), taking (10.9) into account, we have
$$
\|  g^\eps b(\D) \v_\eps - \widetilde{g}^\eps S_\eps b(\D) \widetilde{\u}_0 \|_{L_2(\O)} \leq
{\mathfrak C}_{15} c(\varphi) \eps \|\FF\|_{L_2(\O)},
\eqno(11.34)
$$
where ${\mathfrak C}_{15} = 2 ({\mathcal C}' + {\mathcal C}'') C_\O^{(2)} c^\circ$.
Relations (11.33) and (11.34) imply (10.6) with
$\widetilde{{\mathfrak C}}_2 = \|g\|_{L_\infty} (d \alpha_1)^{1/2} {\mathfrak C}_2$ and
$\widetilde{{\mathfrak C}}_3 = \|g\|_{L_\infty} (d \alpha_1)^{1/2} {\mathfrak C}_3 + {\mathfrak C}_{15}$. $\bullet$

\smallskip\noindent\textbf{11.3. Proof of Theorem 10.1.}

\smallskip\noindent\textbf{Lemma 11.3.} \textit{Under the assumptions of Lemma} 11.2 \textit{for $0< \eps \leq \eps_1$ we have}
$$
\| \w_\eps \|_{L_2(\O)} \leq c(\varphi)^5 \left( {\mathfrak C}_{16} |\zeta|^{-1/2} \eps +
{\mathfrak C}_{17}\eps^2 \right)\|\FF\|_{L_2(\O)}.
\eqno(11.35)
$$
\textit{The constants ${\mathfrak C}_{16}$ and ${\mathfrak C}_{17}$ depend on} $d$, $m$, $\alpha_0$, $\alpha_1$,
$\| g \|_{L_\infty}$, $\| g^{-1} \|_{L_\infty}$,
\textit{the constants $k_1$, $k_2$ from} (9.2), \textit{the parameters of the lattice} $\Gamma$, \textit{and the domain} $\O$.

\smallskip\noindent\textbf{Proof.}
Consider the identity (11.3) and substitute ${\boldsymbol \eta}=
{\boldsymbol \eta}_\eps = (\A_{N,\eps} - \overline{\zeta} I)^{-1}{\boldsymbol \Phi}$, where ${\boldsymbol \Phi} \in L_2(\O;\C^n)$.
Then the left-hand side of (11.3) coincides with $(\w_\eps,{\boldsymbol \Phi})_{L_2(\O)}$, and the identity can be written as
$$
(\w_\eps,{\boldsymbol \Phi})_{L_2(\O)}= {\mathcal I}_\eps[{\boldsymbol \eta}_\eps].
\eqno(11.36)
$$

Assume that $0< \eps \leq \eps_1$. To estimate ${\mathcal I}_\eps[{\boldsymbol \eta}_\eps]$, we use relations
(11.5), (11.8), (11.13), and (11.16). We have
$$
|{\mathcal I}_\eps[{\boldsymbol \eta}_\eps]| \leq ({\mathfrak C}^{(1)}+ {\mathfrak C}^{(2)})
c(\varphi) \eps \|\FF\|_{L_2(\O)} \| \D {\boldsymbol \eta}_\eps \|_{L_2(\O)} + |\widetilde{\mathcal I}^{(2)}_\eps[{\boldsymbol \eta}_\eps]|.
$$
Applying Lemma 9.4 to estimate $\| \D {\boldsymbol \eta}_\eps \|_{L_2(\O)}$, we obtain
$$
|{\mathcal I}_\eps[{\boldsymbol \eta}_\eps]| \leq {\mathfrak C}_0({\mathfrak C}^{(1)}+ {\mathfrak C}^{(2)})
c(\varphi)^2 |\zeta|^{-1/2} \eps \|\FF\|_{L_2(\O)} \| {\boldsymbol \Phi}\|_{L_2(\O)} + |\widetilde{\mathcal I}^{(2)}_\eps[{\boldsymbol \eta}_\eps]|.
\eqno(11.37)
$$
Consider the term $\widetilde{\mathcal I}^{(2)}_\eps[{\boldsymbol \eta}_\eps]$. According to (11.17), we have
$$
\widetilde{\mathcal I}^{(2)}_\eps[{\boldsymbol \eta}_\eps]
= \eps \sum_{l,j=1}^d \left( \partial_j( \theta_\eps M_{lj}^\eps S_\eps b(\D) \widetilde{\u}_0), D_l {\boldsymbol \eta}_\eps \right)_{L_2(\O)}.
\eqno(11.38)
$$
To approximate ${\boldsymbol \eta}_\eps$, we apply the already proved Theorem 10.2. We put
${\boldsymbol \eta}_0 = (\A_N^0 - \overline{\zeta} I)^{-1} {\boldsymbol \Phi}$,
$\widetilde{\boldsymbol \eta}_0 = P_\O {\boldsymbol \eta}_0$.
The approximation of ${\boldsymbol \eta}_\eps$ is given by
${\boldsymbol \eta}_0 + \eps \Lambda^\eps S_\eps b(\D) \widetilde{\boldsymbol \eta}_0$.
We rewrite the term (11.38) as
$$
\widetilde{\mathcal I}^{(2)}_\eps[{\boldsymbol \eta}_\eps]
= {\mathcal J}_1(\eps) + {\mathcal J}_2(\eps) + {\mathcal J}_3(\eps),
\eqno(11.39)
$$
where
$$
{\mathcal J}_1(\eps)=
\eps \sum_{l,j=1}^d \left( \partial_j( \theta_\eps M_{lj}^\eps S_\eps b(\D) \widetilde{\u}_0), D_l ({\boldsymbol \eta}_\eps
- {\boldsymbol \eta}_0 - \eps \Lambda^\eps S_\eps b(\D) \widetilde{\boldsymbol \eta}_0) \right)_{L_2(\O)},
$$
$$
{\mathcal J}_2(\eps)=
\eps \sum_{l,j=1}^d \left( \partial_j( \theta_\eps M_{lj}^\eps S_\eps b(\D) \widetilde{\u}_0), D_l {\boldsymbol \eta}_0 \right)_{L_2(\O)},
\eqno(11.40)
$$
$$
{\mathcal J}_3(\eps)=
\eps \sum_{l,j=1}^d \left( \partial_j( \theta_\eps M_{lj}^\eps S_\eps b(\D) \widetilde{\u}_0), D_l
(\eps \Lambda^\eps S_\eps b(\D) \widetilde{\boldsymbol \eta}_0) \right)_{L_2(\O)}.
\eqno(11.41)
$$

Applying Theorem 10.2 and (11.25), we arrive at
$$
\begin{aligned}
&|{\mathcal J}_1(\eps) | \leq c(\varphi) \left( \nu_1 \eps + (\nu_2+\nu_3) |\zeta|^{-1/4} \eps^{1/2}\right)
\cr
&\times d^{1/2} \left( {\mathfrak C}_2 c(\varphi)^2 |\zeta|^{-1/4} \eps^{1/2}
+ {\mathfrak C}_3 c(\varphi)^4 \eps\right)
\|\FF\|_{L_2(\O)} \| {\boldsymbol \Phi} \|_{L_2(\O)}.
\end{aligned}
$$
Hence,
$$
|{\mathcal J}_1(\eps) | \leq c(\varphi)^5 \left( \check{\nu}_1 |\zeta|^{-1/2} \eps + \check{\nu}_2 \eps^2 \right)
\|\FF\|_{L_2(\O)} \| {\boldsymbol \Phi} \|_{L_2(\O)},
\eqno(11.42)
$$
where $\check{\nu}_1 = d^{1/2} \left((\nu_1 + \nu_2 + \nu_3) {\mathfrak C}_2 + (\nu_2+\nu_3){\mathfrak C}_3 \right)$ and
$\check{\nu}_2 = d^{1/2}\left( \nu_1 ({\mathfrak C}_2 + {\mathfrak C}_3) + (\nu_2+\nu_3) {\mathfrak C}_3 \right)$.

By (5.1) and (11.25), the term (11.40) satisfies the following estimate
$$
|{\mathcal J}_2(\eps) | \leq c(\varphi)
\left( \nu_1 \eps + (\nu_2+\nu_3) |\zeta|^{-1/4} \eps^{1/2}\right) \| \FF \|_{L_2(\O)} d^{1/2}
\biggl( \int_{B_\eps} |\D {\boldsymbol \eta}_0 |^2\,d\x \biggr)^{1/2}.
$$
Applying Lemma 3.3($1^\circ$) and Lemma 9.5, we obtain
$$
|{\mathcal J}_2(\eps) | \leq c(\varphi)^2 \left( \check{\nu}_3 |\zeta|^{-1/2} \eps + \check{\nu}_4 \eps^2 \right)
\|\FF\|_{L_2(\O)} \| {\boldsymbol \Phi} \|_{L_2(\O)},
\eqno(11.43)
$$
where $\check{\nu}_3 = (2 d \beta c^\circ {\mathfrak C}_0)^{1/2}(\nu_1 + \nu_2 + \nu_3)$ and
$\check{\nu}_4 = (2 d \beta c^\circ {\mathfrak C}_0)^{1/2}\nu_1$.

It remains to estimate the term (11.41) which can be represented as
$$
{\mathcal J}_3(\eps)=
{\mathcal J}_3^{(1)}(\eps) + {\mathcal J}_3^{(2)}(\eps),
\eqno(11.44)
$$
where
$$
{\mathcal J}_3^{(1)}(\eps)
= \eps \sum_{l,j=1}^d \left( \partial_j( \theta_\eps M_{lj}^\eps S_\eps b(\D) \widetilde{\u}_0),
(D_l \Lambda)^\eps S_\eps b(\D) \widetilde{\boldsymbol \eta}_0 \right)_{L_2(\O)},
\eqno(11.45)
$$
$$
{\mathcal J}_3^{(2)}(\eps)
= \eps \sum_{l,j=1}^d \left( \partial_j( \theta_\eps M_{lj}^\eps S_\eps b(\D) \widetilde{\u}_0),
\eps \Lambda^\eps S_\eps b(\D) D_l \widetilde{\boldsymbol \eta}_0 \right)_{L_2(\O)}.
\eqno(11.46)
$$
By (5.1) and (11.25), the term (11.45) satisfies the estimate
$$
\begin{aligned}
&| {\mathcal J}_3^{(1)}(\eps) | \leq c(\varphi)
\left( \nu_1 \eps + (\nu_2+\nu_3) |\zeta|^{-1/4} \eps^{1/2}\right) \| \FF \|_{L_2(\O)}
\cr
&\times d^{1/2} \left( \int_{(\partial\O)_\eps}
| (\D\Lambda)^\eps S_\eps b(\D) \widetilde{\boldsymbol \eta}_0 |^2\,d\x \right)^{1/2}.
\end{aligned}
$$
Applying Lemma 3.4, (1.4), (1.10), and the analogs of estimates (10.8), (10.9) for $\widetilde{\boldsymbol \eta}_0$, we obtain
$$
|{\mathcal J}_3^{(1)}(\eps) | \leq c(\varphi)^2 \left( \check{\nu}_5 |\zeta|^{-1/2} \eps + \check{\nu}_6 \eps^2 \right)
\|\FF\|_{L_2(\O)} \| {\boldsymbol \Phi} \|_{L_2(\O)},
\eqno(11.47)
$$
where $\check{\nu}_5 = M_2 (2d \beta_* \alpha_1 C_\O^{(1)} C_\O^{(2)} ({\mathfrak C}_0+1)
c^\circ )^{1/2} (\nu_1 + \nu_2 + \nu_3)$ and
 $\check{\nu}_6 = M_2 (2d \beta_* \alpha_1 C_\O^{(1)} C_\O^{(2)} ({\mathfrak C}_0+1) c^\circ )^{1/2} \nu_1$.

Finally, it is easy to estimate the term (11.46) by using (11.25), (1.4), (1.19), and the analog of (10.9) for $\widetilde{\boldsymbol \eta}_0$:
$$
|{\mathcal J}_3^{(2)}(\eps) | \leq c(\varphi)^2 \left( \check{\nu}_7 |\zeta|^{-1/2} \eps + \check{\nu}_8 \eps^2 \right)
\|\FF\|_{L_2(\O)} \| {\boldsymbol \Phi} \|_{L_2(\O)}.
\eqno(11.48)
$$
Here $\check{\nu}_7 = 2 M_1 (d \alpha_1)^{1/2} C_\O^{(2)} c^\circ (\nu_2 + \nu_3),$
 $\check{\nu}_8 = 2 M_1 (d \alpha_1)^{1/2} C_\O^{(2)} c^\circ (\nu_1+ \nu_2 + \nu_3)$.

As a result, relations (11.39), (11.42)--(11.44), (11.47), and (11.48) imply that
$$
| \widetilde{\mathcal I}_\eps^{(2)}[{\boldsymbol \eta}_\eps] | \leq c(\varphi)^5 \left( \check{\nu}_9 |\zeta|^{-1/2} \eps + \check{\nu}_{10} \eps^2 \right)
\|\FF\|_{L_2(\O)} \| {\boldsymbol \Phi} \|_{L_2(\O)},\ \ 0< \eps \leq \eps_1,
$$
where $\check{\nu}_9 = \check{\nu}_1 + \check{\nu}_3 + \check{\nu}_5 + \check{\nu}_7$ and
 $\check{\nu}_{10} = \check{\nu}_2 + \check{\nu}_4 + \check{\nu}_6 + \check{\nu}_8$.
Together with (11.37) this yields
$$
|{\mathcal I}_\eps [{\boldsymbol \eta}_\eps] | \leq c(\varphi)^5 \left( {\mathfrak C}_{16} |\zeta|^{-1/2} \eps
+ {\mathfrak C}_{17} \eps^2 \right)
\|\FF\|_{L_2(\O)} \| {\boldsymbol \Phi} \|_{L_2(\O)},\ \ 0< \eps \leq \eps_1,
\eqno(11.49)
$$
 where ${\mathfrak C}_{16} = {\mathfrak C}_0({\mathfrak C}^{(1)} + {\mathfrak C}^{(2)}) + \check{\nu}_9$ and
${\mathfrak C}_{17} = \check{\nu}_{10}$.

From (11.36) and (11.49) it follows that
$$
|( \w_\eps, {\boldsymbol \Phi})_{L_2(\O)}| \leq c(\varphi)^5 \left( {\mathfrak C}_{16} |\zeta|^{-1/2} \eps +
{\mathfrak C}_{17} \eps^2\right)
\|\FF\|_{L_2(\O)} \| {\boldsymbol \Phi} \|_{L_2(\O)}
$$
for $0< \eps \leq \eps_1$ and any ${\boldsymbol \Phi} \in L_2(\O;\C^n)$,
which implies the required estimate (11.35). $\bullet$

\smallskip\noindent\textbf{Completion of the proof of Theorem 10.1.}
Estimate (10.25) and Lemma 11.3 show that (10.1) is true for $0< \eps \leq \eps_1$ with ${\mathfrak C}_1 = \max \{ {\mathfrak C}_{9}
+ {\mathfrak C}_{16}, {\mathfrak C}_{17}\}$. $\bullet$

\section*{\S 12. The results for the Neumann problem: \\ the case where $\Lambda \in L_\infty$, special cases,
\\ estimates in a strictly interior subdomain}

\smallskip\noindent\textbf{12.1. The case where $\Lambda \in L_\infty$.}
As for the problem in $\R^d$ (see Subsection 2.3) and for the Dirichlet problem (see Subsection 6.1),
under Condition 2.8 the operator $S_\eps$ in the corrector can be removed.
Now instead of the corrector (10.2) we will use the operator
$$
K^0_N(\eps;\zeta) = [\Lambda^\eps] b(\D) (\A^0_N - \zeta I)^{-1}.
\eqno(12.1)
$$
Under Condition 2.8 the operator (12.1) is a continuous mapping of $L_2(\O;\C^n)$ to $H^1(\O;\C^n)$.
Instead of (10.4) we use another approximation of the solution $\u_\eps$ of problem (9.6):
$$
\check{\v}_\eps := (\A^0_N - \zeta I)^{-1}\FF +\eps K^0_N(\eps;\zeta)\FF = \u_0 + \eps \Lambda^\eps b(\D)\u_0.
\eqno(12.2)
$$

\smallskip\noindent\textbf{Theorem 12.1.} \textit{Suppose that the assumptions of Theorem} 10.1 \textit{and Condition} 2.8
 \textit{are satisfied. Let $\check{\v}_\eps$ be defined by} (12.2).
\textit{Then for $0< \varepsilon \leq \varepsilon_1$ we have}
$$
\| \u_\eps - \check{\v}_\eps \|_{H^1(\O)} \leq
\left( {\mathfrak C}_2 c(\varphi)^2 |\zeta|^{-1/4} \varepsilon^{1/2} +
{\mathfrak C}_3^\circ c(\varphi)^4 \varepsilon\right) \|\FF\|_{L_2(\O)},
\eqno(12.3)
$$
\textit{or, in operator terms},
$$
\begin{aligned}
&\| (\A_{N,\varepsilon}- \zeta I)^{-1}  - (\A^0_{N}- \zeta I)^{-1} - \eps K^0_N(\eps;\zeta) \|_{L_2({\mathcal O}) \to H^1({\mathcal O})}
\cr
&\leq {\mathfrak C}_2 c(\varphi)^2  |\zeta|^{-1/4}\varepsilon^{1/2} + {\mathfrak C}_3^\circ c(\varphi)^4 \varepsilon.
\end{aligned}
$$
\textit{For the flux} $\p_\eps := g^\eps b(\D) \u_\eps$ \textit{we have}
$$
\| \p_\eps - \widetilde{g}^\eps b(\D) {\u}_0 \|_{L_2(\O)}
\leq \left(\widetilde{{\mathfrak C}}_2 c(\varphi)^2 |\zeta|^{-1/4}\varepsilon^{1/2} +
\widetilde{{\mathfrak C}}_3^\circ c(\varphi)^4 \varepsilon\right) \|\FF\|_{L_2(\O)}
\eqno(12.4)
$$
\textit{for $0< \eps \leq \eps_1$. The constants ${\mathfrak C}_2$ and $\widetilde{{\mathfrak C}}_2$ are the same as in Theorem} 10.2.
\textit{The constants ${\mathfrak C}_3^\circ$ and $\widetilde{{\mathfrak C}}_3^\circ$ depend on}
$d$, $m$, $\alpha_0$, $\alpha_1$, $\| g \|_{L_\infty}$, $\| g^{-1} \|_{L_\infty}$,
\textit{the constants $k_1$, $k_2$ from} (9.2), \textit{the parameters of the lattice} $\Gamma$, \textit{the domain $\O$, and the norm}
$\|\Lambda\|_{L_\infty}$.

\smallskip\noindent\textbf{Proof.} In order to deduce (12.3) from (10.5), we need to estimate the difference of the functions (10.4) and (12.2).
Similarly to (6.5)--(6.9), it is easy to check that
$$
\| \v_\eps - \check{\v}_\eps \|_{H^1(\O)} \leq {\mathcal C}'''  \eps \| \widetilde{\u}_0 \|_{H^2(\R^d)}.
$$
Combining this with (10.5) and (10.9), we obtain (12.3) with
${\mathfrak C}_3^\circ = {\mathfrak C}_3 + 2 {\mathcal C}''' C_{\O}^{(2)} c^\circ$.

Let us check (12.4). From (12.3) and (1.2), (1.5) it follows that
$$
\begin{aligned}
&\| \p_\eps - g^\eps b(\D) \check{\v}_\eps \|_{L_2(\O)}
\cr
&\leq
\|g\|_{L_\infty} (d \alpha_1)^{1/2} \left( {\mathfrak C}_{2} c(\varphi)^2 |\zeta|^{-1/4} \varepsilon^{1/2} +
{\mathfrak C}_3^\circ c(\varphi)^4 \varepsilon\right) \|\FF\|_{L_2(\O)}.
\end{aligned}
\eqno(12.5)
$$
Similarly to  (6.11)--(6.13), we check that
$$
\| g^\eps b(\D) \check{\v}_\eps - \widetilde{g}^\eps b(\D) \u_0 \|_{L_2(\O)} \leq
\widetilde{\mathcal C}'  \varepsilon \| \u_0 \|_{H^2(\O)}.
\eqno(12.6)
$$
Relations (9.14), (12.5), and (12.6) imply (12.4) with
$\widetilde{\mathfrak C}_3^\circ = \|g\|_{L_\infty} (d \alpha_1)^{1/2} {\mathfrak C}_3^\circ +2\widetilde{\mathcal C}' c^\circ$. $\bullet$

\smallskip\noindent\textbf{12.2. Special cases.} The following two statements are similar to
Propositions 6.2 and 6.3.

\smallskip\noindent\textbf{Proposition 12.2.} \textit{Suppose that the assumptions of Theorem} 10.2
 \textit{are satisfied. If $g^0 = \overline{g}$, i.~e., relations} (1.13) \textit{are satisfied, then
for $0< \varepsilon \leq \varepsilon_1$ we have}
$$
\| \u_\eps - {\u}_0 \|_{H^1(\O)} \leq \left({\mathfrak C}_2 c(\varphi)^2 |\zeta|^{-1/4}\varepsilon^{1/2} +
 {\mathfrak C}_3  c(\varphi)^4 \varepsilon\right) \|\FF\|_{L_2(\O)}.
$$

\smallskip\noindent\textbf{Proposition 12.3.} \textit{Suppose that the assumptions of Theorem} 10.2
 \textit{are satisfied. If $g^0 = \underline{g}$, i.~e., relations} (1.14) \textit{are satisfied, then
for $0< \varepsilon \leq \varepsilon_1$ we have}
$$
\| \p_\eps - g^0 b(\D){\u}_0 \|_{L_2(\O)} \leq \left(\widetilde{\mathfrak C}_2 c(\varphi)^2 |\zeta|^{-1/4}\varepsilon^{1/2} +
\widetilde{{\mathfrak C}}_3^\circ c(\varphi)^4 \varepsilon\right) \|\FF\|_{L_2(\O)}.
$$

\smallskip\noindent\textbf{12.3. Approximation of the solutions of the Neumann problem in a strictly interior
subdomain.}
By analogy with the proof of Theorem 7.1, it is easy to prove the following theorem on the basis of
Theorem 10.1 and the results for the problem in $\R^d$. We omit the proof.

\smallskip\noindent\textbf{Theorem 12.4.} \textit{Suppose that the assumptions of Theorem} 10.2
\textit{are satisfied. Let $\O'$ be a strictly interior subdomain of the domain $\O$, and let}
$\delta:= {\rm dist}\,\{ \O'; \partial \O \}$. \textit{Then for $0< \eps \leq \eps_1$ we have}
$$
\|\u_\eps - \v_\eps \|_{H^1(\O')} \leq ( {\mathfrak C}'_{18} \delta^{-1} + {\mathfrak C}''_{18}) c(\varphi)^6 \eps \| \FF \|_{L_2(\O)},
$$
\textit{or, in operator terms,}
$$
\begin{aligned}
&\| (\A_{N,\eps} - \zeta I)^{-1} - (\A_{N}^0 - \zeta I)^{-1} - \eps K_N(\eps;\zeta) \|_{L_2(\O) \to H^1(\O')}
\cr
&\leq ( {\mathfrak C}'_{18} \delta^{-1} + {\mathfrak C}''_{18}) c(\varphi)^6 \eps.
\end{aligned}
$$
\textit{For the flux $\p_\eps = g^\eps b(\D) \u_\eps$ we have}
$$
\| \p_\eps - \widetilde{g}^\eps S_\eps b(\D) \widetilde{\u}_0
\|_{L_2(\O')} \leq (\widetilde{\mathfrak C}'_{18} \delta^{-1} + \widetilde{\mathfrak C}''_{18})
c(\varphi)^6 \eps \| \FF \|_{L_2(\O)}
$$
\textit{for $0< \eps \leq \eps_1$. The constants ${\mathfrak C}'_{18}$, ${\mathfrak C}''_{18}$, $\widetilde{{\mathfrak C}}'_{18}$, and
$\widetilde{{\mathfrak C}}''_{18}$ depend on}
$d$, $m$, $\alpha_0$, $\alpha_1$, $\| g \|_{L_\infty}$, $\| g^{-1} \|_{L_\infty}$, \textit{the constants $k_1$, $k_2$ from} (9.2),
\textit{the parameters of the lattice} $\Gamma$, \textit{and the domain $\O$}.

\smallskip
Similarly to the proof of Theorem 7.2 it is easy to check the following result.

\smallskip\noindent\textbf{Theorem 12.5.} \textit{Suppose that the assumptions of Theorem} 12.1
\textit{are satisfied. Let $\O'$ be a strictly interior subdomain of the domain $\O$, and let}
$\delta:= {\rm dist}\,\{ \O'; \partial \O \}$. \textit{Then for $0< \eps \leq \eps_1$ we have}
$$
\|\u_\eps - \check{\v}_\eps \|_{H^1(\O')} \leq ({\mathfrak C}'_{18} \delta^{-1} + \check{{\mathfrak C}}''_{18}) c(\varphi)^6 \eps \| \FF \|_{L_2(\O)},
$$
\textit{or, in operator terms,}
$$
\| (\A_{N,\eps} - \zeta I)^{-1} - (\A_{N}^0 - \zeta I)^{-1} - \eps K_N^0(\eps;\zeta) \|_{L_2(\O) \to H^1(\O')}
\leq ({\mathfrak C}'_{18} \delta^{-1} + \check{\mathfrak C}''_{18}) c(\varphi)^6 \eps.
$$
\textit{For the flux $\p_\eps = g^\eps b(\D) \u_\eps$ we have}
$$
\| \p_\eps - \widetilde{g}^\eps b(\D) {\u}_0
\|_{L_2(\O')} \leq (\widetilde{{\mathfrak C}}'_{18} \delta^{-1} + \widehat{{\mathfrak C}}''_{18}) c(\varphi)^6 \eps \| \FF \|_{L_2(\O)}
$$
\textit{for $0< \eps \leq \eps_1$. The constants ${\mathfrak C}'_{18}$ and $\widetilde{{\mathfrak C}}'_{18}$ are the same as in Theorem} 12.4.
\textit{The constants $\check{\mathfrak C}''_{18}$ and $\widehat{\mathfrak C}''_{18}$ depend on}
$d$, $m$, $\alpha_0$, $\alpha_1$, $\| g \|_{L_\infty}$, $\| g^{-1} \|_{L_\infty}$, \textit{the constants $k_1$, $k_2$ from} (9.2),
\textit{the parameters of the lattice} $\Gamma$, \textit{the domain $\O$, and the norm $\|\Lambda\|_{L_\infty}$}.

\section*{\S 13. Another approximation of the resolvent $(\A_{N,\eps} - \zeta I)^{-1}$}

\smallskip\noindent\textbf{13.1. Approximation of the resolvent $(\A_{N,\eps} - \zeta I)^{-1}$ for
$\zeta \in \C \setminus \R_+$.}
The next result is similar to Theorem 8.1.

\smallskip\noindent\textbf{Theorem 13.1.} \textit{Let $\zeta = |\zeta|e^{i\varphi} \in \C \setminus \R_+$.
Denote}
$$
\rho_0(\zeta) =
\begin{cases}
c(\varphi)^2 |\zeta|^{-2}, & |\zeta|<1, \\
c(\varphi)^2, & |\zeta| \ge 1.
\end{cases}
$$
\textit{Let $\u_\eps$ be the solution of problem} (9.6),
 \textit{and let $\u_0$  be the solution of problem} (9.13).  \textit{Let $\v_\eps$ be defined by} (10.4).
\textit{Suppose that the number $\varepsilon_1$ satisfies Condition} 4.1.
\textit{Then for $0< \varepsilon \leq \varepsilon_1$ we have}
$$
\| \u_\eps - \u_0 \|_{L_2(\O)} \leq {{\mathfrak C}}_{19} \rho_0(\zeta) \varepsilon \|\FF\|_{L_2(\O)},
\eqno(13.1)
$$
$$
\| \u_\eps - \v_\eps \|_{H^1(\O)} \leq
{{\mathfrak C}}_{20} \rho_0(\zeta) \varepsilon^{1/2} \|\FF\|_{L_2(\O)},
\eqno(13.2)
$$
\textit{or, in operator terms,}
$$
\| (\A_{N,\varepsilon}- \zeta I)^{-1}  - (\A^0_{N}- \zeta I)^{-1} \|_{L_2({\mathcal O}) \to L_2({\mathcal O})}
\leq    {\mathfrak C}_{19} \rho_0(\zeta) \varepsilon,
\eqno(13.3)
$$
$$
\| (\A_{N,\varepsilon}- \zeta I)^{-1}  - (\A^0_{N}- \zeta I)^{-1} - \eps K_N(\eps;\zeta)\|_{L_2({\mathcal O}) \to H^1({\mathcal O})}
\leq    {\mathfrak C}_{20} \rho_0(\zeta) \varepsilon^{1/2}.
\eqno(13.4)
$$
\textit{For the flux $\p_\eps= g^\eps b(\D)\u_\eps$ we have}
$$
\| \p_\eps - \widetilde{g}^\eps S_\eps b(\D) \widetilde{\u}_0 \|_{L_2(\O)}
\leq    {\mathfrak C}_{21} \rho_0(\zeta) \varepsilon^{1/2} \|\FF\|_{L_2(\O)}, \ \ 0< \eps \leq \eps_1.
\eqno(13.5)
$$
\textit{The constants ${\mathfrak C}_{19}$, ${\mathfrak C}_{20}$, and ${\mathfrak C}_{21}$ depend on}
$d$, $m$, $\alpha_0$, $\alpha_1$, $\| g \|_{L_\infty}$, $\| g^{-1} \|_{L_\infty}$,
\textit{the constants $k_1$, $k_2$ from} (9.2), \textit{the parameters of the lattice} $\Gamma$, \textit{and the domain} $\O$.

\smallskip\noindent\textbf{Proof.}
Applying Theorem 10.1 with $\zeta=-1$ and the analog of the identity (2.10), for $0< \eps \leq \eps_1$ we obtain
$$
\| (\A_{N,\varepsilon}- \zeta I)^{-1}  - (\A^0_{N}- \zeta I)^{-1} \|_{L_2({\mathcal O}) \to L_2({\mathcal O})}
\leq   2 {\mathfrak C}_1 \eps \sup_{x\geq 0}(x+1)^2 |x-\zeta|^{-2}.
\eqno(13.6)
$$
A calculation shows that
$$
\sup_{x\geq 0}(x+1)^2 |x-\zeta|^{-2} \leq 4 \rho_0(\zeta).
\eqno(13.7)
$$
Relations (13.6) and (13.7) imply (13.3) with the constant ${\mathfrak C}_{19} = 8 {\mathfrak C}_1$.

Applying Theorem 10.2 with $\zeta =-1$, for $0< \eps \leq \eps_1$ we have
$$
\| (\A_{N,\varepsilon}+ I)^{-1}  - (\A^0_{N}+ I)^{-1} - \eps K_N(\eps;-1)\|_{L_2({\mathcal O}) \to H^1({\mathcal O})}
\leq    ({\mathfrak C}_2+{\mathfrak C}_3) \varepsilon^{1/2}.
\eqno(13.8)
$$
We put $\lambda := \|g^{-1}\|_{L_\infty}^{-1}(k_1+k_2)$. Note that, by (9.5),
$$
\|\v\|^2_{H^1(\O)} \leq k_1^{-1} \|g^{-1}\|_{L_\infty} \| (\A_{N,\eps} + \lambda I)^{1/2}\v\|_{L_2(\O)}^2,\ \ \v \in H^1(\O;\C^n).
\eqno(13.9)
$$
Using the analog of the identity (2.34) and (13.7), we obtain
$$
\begin{aligned}
&\| (\A_{N,\eps} + \lambda I)^{1/2} \left( (\A_{N,\varepsilon}- \zeta I)^{-1}  - (\A^0_{N}- \zeta I)^{-1} - \eps K_N(\eps;\zeta)\right)
\|_{L_2 \to L_2}
\cr
&\leq  4 \rho_0(\zeta)
\| (\A_{N,\eps} + \lambda I)^{1/2} \left( (\A_{N,\varepsilon} + I)^{-1}  - (\A^0_{N} + I)^{-1} - \eps K_N(\eps;-1)\right)
\|_{L_2 \to L_2}
\cr
&+ \eps |\zeta +1| \sup_{x \geq 0} (x+\lambda)^{1/2} |x-\zeta|^{-1} \| K_N(\eps;\zeta)\|_{L_2\to L_2}.
\end{aligned}
\eqno(13.10)
$$
We denote the summands in the right-hand side of (13.10) by $L_1(\eps)$ and $L_2(\eps)$.
By (9.4) and (13.8), the first term satisfies
$$
L_1(\eps) \leq {\mathfrak C}_{22} \rho_0(\zeta) \eps^{1/2},\ \ 0< \eps \leq \eps_1,
\eqno(13.11)
$$
where ${\mathfrak C}_{22} = 4 (\max\{ d \alpha_1 \|g\|_{L_\infty}, \lambda\})^{1/2}
({\mathfrak C}_2 + {\mathfrak C}_3)$.

Next, the operator $K_N(\eps;\zeta)$ can be represented as
$K_N(\eps;\zeta) = R_\O [\Lambda^\eps] S_\eps b(\D) P_\O (\A_N^0 + \lambda I)^{-1/2} (\A_N^0 + \lambda I)^{1/2} (\A_N^0 -\zeta I)^{-1}$.
Combining this with (1.4), (1.19), and (4.3), we obtain
$$
L_2(\eps) \leq \eps |\zeta+1| \left(\sup_{x\geq 0} (x+\lambda)|x-\zeta|^{-2}\right) M_1 \alpha_1^{1/2} C_\O^{(1)}
\| (\A_N^0 + \lambda I)^{-1/2}\|_{L_2\to H^1}.
\eqno(13.12)
$$
A calculation shows that
$$
|\zeta+1| \sup_{x\geq 0} (x+\lambda)|x-\zeta|^{-2} \leq 2(\lambda+1) \rho_0(\zeta).
\eqno(13.13)
$$
Similarly to (13.9),
$$
\| (\A_N^0 + \lambda I)^{-1/2}\|_{L_2(\O)\to H^1(\O)} \leq k_1^{-1/2} \|g^{-1}\|_{L_\infty}^{1/2}.
\eqno(13.14)
$$
From (13.12)--(13.14) it follows that
$$
L_2(\eps) \leq {\mathfrak C}_{23} \rho_0(\zeta) \eps,
\eqno(13.15)
$$
where ${\mathfrak C}_{23} = 2 (\lambda+1 ) M_1  \alpha_1^{1/2} C_\O^{(1)} k_1^{-1/2} \|g^{-1}\|^{1/2}_{L_\infty}$.
Now relations (13.10), (13.11), (13.15), and (13.9) imply (13.4) with the constant
${\mathfrak C}_{20} = k_1^{-1/2} \|g^{-1}\|_{L_\infty}^{1/2}({\mathfrak C}_{22} + {\mathfrak C}_{23})$.

Let us check (13.5). From (13.2) it follows that
$$
\| \p_\eps - g^\eps b(\D) \v_\eps\|_{L_2(\O)} \leq \| g \|_{L_\infty} (d \alpha_1)^{1/2}
{\mathfrak C}_{20} \rho_0(\zeta) \eps^{1/2} \| \FF \|_{L_2(\O)}
\eqno(13.16)
$$
for $0< \eps \leq \eps_1$. Next, similarly to (5.19)--(5.21) we have
$$
\| g^\eps b(\D) \v_\eps - \widetilde{g}^\eps S_\eps b(\D) \widetilde{\u}_0 \|_{L_2(\O)}
\leq ({\mathcal C}' + {\mathcal C}'') C_\O^{(2)} \eps \| \u_0 \|_{H^2(\O)}.
\eqno(13.17)
$$
By (9.12) and  (13.7),
$$
\| (\A_N^0 - \zeta I)^{-1}\|_{L_2(\O) \to H^2(\O)} \leq 2 c^\circ \rho_0(\zeta)^{1/2}.
$$
Hence,
$$
\| \u_0 \|_{H^2(\O)} \leq  2 c^\circ \rho_0(\zeta)^{1/2} \| \FF \|_{L_2(\O)}.
\eqno(13.18)
$$
Relations (13.16), (13.17), and (13.18) imply (13.5) with
${\mathfrak C}_{21} =  {\mathfrak C}_{20} \| g \|_{L_\infty} (d \alpha_1)^{1/2}+
2 ({\mathcal C}' + {\mathcal C}'') C_\O^{(2)}  c^\circ$. $\bullet$

\smallskip
Next theorem is an analog of Theorem 8.3.

\smallskip\noindent\textbf{Theorem 13.2.} \textit{Suppose that the assumptions of Theorem} 13.1 \textit{and Condition} 2.8
\textit{are satisfied. Let $\check{\v}_\eps$ be defined by} (12.2).
\textit{Then for $0< \varepsilon \leq \varepsilon_1$ we have}
$$
\| \u_\eps - \check{\v}_\eps \|_{H^1(\O)} \leq
{\mathfrak C}_{20}^\circ  \rho_0(\zeta) \varepsilon^{1/2} \|\FF\|_{L_2(\O)},
\eqno(13.19)
$$
\textit{or, in operator terms,}
$$
\| (\A_{N,\varepsilon}- \zeta I)^{-1}  - (\A^0_{N}- \zeta I)^{-1} - \eps K_N^0(\eps;\zeta)\|_{L_2({\mathcal O}) \to H^1({\mathcal O})}
\leq   {\mathfrak C}_{20}^\circ \rho_0(\zeta) \varepsilon^{1/2}.
$$
\textit{For the flux $\p_\eps= g^\eps b(\D)\u_\eps$ we have}
$$
\| \p_\eps - \widetilde{g}^\eps  b(\D) {\u}_0 \|_{L_2(\O)}
\leq    {\mathfrak C}_{21}^\circ \rho_0(\zeta) \varepsilon^{1/2} \|\FF\|_{L_2(\O)}, \ \ 0< \eps \leq \eps_1.
\eqno(13.20)
$$
\textit{The constants ${\mathfrak C}_{20}^\circ$ and ${\mathfrak C}_{21}^\circ$ depend on}
$d$, $m$, $\alpha_0$, $\alpha_1$, $\| g \|_{L_\infty}$, $\| g^{-1} \|_{L_\infty}$,
\textit{the constants $k_1$, $k_2$ from} (9.2), \textit{the parameters of the lattice} $\Gamma$, \textit{the domain $\O$, and
the norm $\|\Lambda\|_{L_\infty}$}.

\smallskip\noindent\textbf{Proof.}
Similarly to (6.5)--(6.9), we have
$$
\| \v_\eps - \check{\v}_\eps\|_{H^1(\O)} \leq {\mathcal C}''' \eps \| \widetilde{\u}_0\|_{H^2(\R^d)}.
\eqno(13.21)
$$
Relations (4.3), (13.2), (13.18), and (13.21) imply (13.19) with
${\mathfrak C}_{20}^\circ = {\mathfrak C}_{20} + 2 {\mathcal C}''' C^{(2)}_{\O} c^\circ$.

Now we check (13.20). From (13.19) it follows that
$$
\| \p_\eps - g^\eps b(\D) \check{\v}_\eps \|_{L_2(\O)}
\leq \|g\|_{L_\infty} (d \alpha_1)^{1/2}  {\mathfrak C}_{20}^\circ \rho_0(\zeta)  \eps^{1/2} \| \FF \|_{L_2(\O)}
\eqno(13.22)
$$
 for $0< \eps \leq \eps_1$. Similarly to (6.11) and (6.12), we have
$$
\|  g^\eps b(\D) \check{\v}_\eps - \widetilde{g}^\eps b(\D) \u_0 \|_{L_2(\O)}
\leq \widetilde{\mathcal C}' \eps \| \u_0 \|_{H^2(\O)}.
\eqno(13.23)
$$
Combining (13.22), (13.23), and (13.18), we obtain (13.20) with
${\mathfrak C}_{21}^\circ = \|g\|_{L_\infty}  (d \alpha_1)^{1/2} {\mathfrak C}_{20}^\circ +
2 c^\circ \widetilde{\mathcal C}'$. $\bullet$

\smallskip\noindent\textbf{13.2. Special cases.}
The following statement concerns the case where the corrector is equal to zero (cf. Proposition 8.4).

\smallskip\noindent\textbf{Proposition 13.3.} \textit{Suppose that the assumptions of Theorem} 13.1
\textit{are satisfied. If $g^0 = \overline{g}$, i.~e., relations} (1.13) \textit{are satisfied, then
$\Lambda=0$, $\v_\eps = \u_0$, and for $0< \eps \leq \eps_1$ we have}
$$
\| \u_\eps - \u_0 \|_{H^1(\O)} \leq {\mathfrak C}_{20} \rho_0(\zeta) \eps^{1/2} \| \FF \|_{L_2(\O)}.
$$

\smallskip
Next statement is similar to Proposition 8.5.

\smallskip\noindent\textbf{Proposition 13.4.} \textit{Suppose that the assumptions of Theorem} 13.1
\textit{are satisfied. If $g^0 = \underline{g}$, i.~e., relations} (1.14) \textit{are satisfied, then
for $0< \eps \leq \eps_1$ we have}
$$
\| \p_\eps - g^0 b(\D)\u_0 \|_{L_2(\O)} \leq {\mathfrak C}_{21}^\circ \rho_0(\zeta) \eps^{1/2} \| \FF \|_{L_2(\O)}.
$$

\smallskip\noindent\textbf{13.3. Approximation of the resolvent $(\A_{N,\eps} - \zeta I)^{-1}$ for $\zeta \in \C \setminus \R_+$
in a strictly interior subdomain.}
The following result is an analog of Theorem 8.6.

\smallskip\noindent\textbf{Theorem 13.5.} \textit{Suppose that the assumptions of Theorem} 13.1
\textit{are satisfied. Let $\O'$ be a strictly interior subdomain of the domain $\O$, and let $\delta = {\rm dist}\, \{\O'; \partial \O\}$.}
\textit{Then for $0< \varepsilon \leq \varepsilon_1$ we have}
$$
\begin{aligned}
&\| \u_\eps - \v_\eps \|_{H^1(\O')}
\cr
&\leq ({\mathcal C}_{23} \delta^{-1} +1) \left(
{\mathfrak C}_{19} c(\varphi) \rho_0(\zeta)+
{\mathfrak C}_{24} c(\varphi)^{5/2} \rho_0(\zeta)^{3/4}\right)  \varepsilon \|\FF\|_{L_2(\O)},
\end{aligned}
\eqno(13.24)
$$
\textit{or, in operator terms},
$$
\begin{aligned}
&\| (\A_{N,\eps} - \zeta I)^{-1} - (\A_{N}^0  - \zeta I)^{-1} - \eps K_N(\eps;\zeta) \|_{L_2(\O) \to H^1(\O')}
\cr
&\leq ({\mathcal C}_{23} \delta^{-1} +1) \left({\mathfrak C}_{19} c(\varphi) \rho_0(\zeta)+
{\mathfrak C}_{24} c(\varphi)^{5/2} \rho_0(\zeta)^{3/4}\right)  \varepsilon.
\end{aligned}
$$
\textit{For the flux $\p_\eps = g^\eps b(\D) \u_\eps$ we have}
$$
\begin{aligned}
&\| \p_\eps - \widetilde{g}^\eps S_\eps b(\D) \widetilde{\u}_0 \|_{L_2(\O')}
\cr
&\leq ({\mathcal C}_{23} \delta^{-1} +1) ({\mathfrak C}_{25} c(\varphi) \rho_0(\zeta)+
 {\mathfrak C}_{26} c(\varphi)^{5/2} \rho_0(\zeta)^{3/4})  \varepsilon \|\FF\|_{L_2(\O)}
\end{aligned}
\eqno(13.25)
$$
\textit{for $0< \eps \leq \eps_1$. Here the constant ${\mathcal C}_{23}$ is the same as in} (7.14),
${\mathfrak C}_{19}$ \textit{is the constant from} (13.1). \textit{The constants 
${\mathfrak C}_{24}$, ${\mathfrak C}_{25}$, and ${\mathfrak C}_{26}$ depend on}
$d$, $m$, $\alpha_0$, $\alpha_1$, $\| g \|_{L_\infty}$, $\| g^{-1} \|_{L_\infty}$,
\textit{the constants $k_1$, $k_2$ from} (9.2), \textit{the parameters of the lattice} $\Gamma$, \textit{and the domain} $\O$.

\smallskip\noindent\textbf{Proof.} The first part of the proof is similar to that of Theorem 7.1.
Let $\chi(\x)$ be a cut-off function satisfying
(7.3). Let $\u_\eps$ be the solution of problem (9.6), and let $\widetilde{\u}_\eps$ be the solution of equation (10.12). Then
$(\A_\eps - \zeta) ( \u_\eps - \widetilde{\u}_\eps)=0$ in the domain $\O$. Similarly to (7.4)--(7.14), we deduce
$$
\| \D (\chi ( \u_\eps - \widetilde{\u}_\eps))\|_{L_2(\O)} \leq {\mathcal C}_{23} c(\varphi) \delta^{-1}
\| \u_\eps - \widetilde{\u}_\eps\|_{L_2(\O)}.
\eqno(13.26)
$$
Next, by Theorem 2.2,
$$
 \| \widetilde{\u}_\eps  - \widetilde{\u}_0 \|_{L_2(\R^d)} \leq C_1 c(\varphi)^2 |\zeta|^{-1/2} \eps \| \widetilde{\FF} \|_{L_2(\R^d)}.
\eqno(13.27)
$$
Relations (1.16) and (4.3) imply the following estimate for the function (10.10):
$$
\begin{aligned}
&\| \widetilde{\FF} \|_{L_2(\R^d)} \leq c_1 \| \widetilde{\u}_0 \|_{H^2(\R^d)} + |\zeta| \| \widetilde{\u}_0 \|_{L_2(\R^d)}
\cr
&\leq c_1 C_\O^{(2)} \| {\u}_0 \|_{H^2(\O)} + |\zeta| C_\O^{(0)} \| {\u}_0 \|_{L_2(\O)}.
\end{aligned}
$$
Combining this with (13.18) and the estimate $\|\u_0\|_{L_2(\O)} \leq |\zeta|^{-1} c(\varphi) \|\FF\|_{L_2(\O)}$, we obtain
$$
\| \widetilde{\FF} \|_{L_2(\R^d)} \leq {\mathfrak C}_{27} \rho_0(\zeta)^{1/2} \|\FF \|_{L_2(\O)},
\eqno(13.28)
$$
where ${\mathfrak C}_{27}= 2 c_1 C_\O^{(2)} c^\circ + C_\O^{(0)}$. From (13.27) and (13.28) it follows that
$$
 \| \widetilde{\u}_\eps  - \widetilde{\u}_0 \|_{L_2(\R^d)} \leq C_1
 {\mathfrak C}_{27} c(\varphi)^2 |\zeta|^{-1/2} \rho_0(\zeta)^{1/2} \eps \| {\FF} \|_{L_2(\O)}.
$$
Together with (13.1) this yields
$$
 \| {\u}_\eps  - \widetilde{\u}_\eps \|_{L_2(\O)} \leq
 \left({\mathfrak C}_{19} \rho_0(\zeta) + C_1 {\mathfrak C}_{27} c(\varphi)^{3/2} \rho_0(\zeta)^{3/4} \right)
\eps \| {\FF} \|_{L_2(\O)}.
\eqno(13.29)
$$
Now, comparing (13.26) and (13.29), we have
$$
\begin{aligned}
&\| \D (\chi ( \u_\eps - \widetilde{\u}_\eps))\|_{L_2(\O)}
\cr
&\leq {\mathcal C}_{23}  \delta^{-1}
\left( {\mathfrak C}_{19} c(\varphi) \rho_0(\zeta) + C_1 {\mathfrak C}_{27} c(\varphi)^{5/2} \rho_0(\zeta)^{3/4} \right)
\eps \| \FF \|_{L_2(\O)}.
\end{aligned}
$$
Hence,
$$
 \begin{aligned}
 &\| {\u}_\eps  - \widetilde{\u}_\eps \|_{H^1(\O')}
 \cr
 &\leq
({\mathcal C}_{23} \delta^{-1} +1) \left({\mathfrak C}_{19} c(\varphi) \rho_0(\zeta) +
C_1 {\mathfrak C}_{27} c(\varphi)^{5/2} \rho_0(\zeta)^{3/4} \right)
\eps \| {\FF} \|_{L_2(\O)}.
\end{aligned}
\eqno(13.30)
$$
By Corollary 2.5 and (13.28),
$$
\begin{aligned}
&\| \widetilde{\u}_\eps - \widetilde{\v}_\eps \|_{H^1(\R^d)} \leq (C_2 + C_3 |\zeta|^{-1/2}) c(\varphi)^2 \eps \| \widetilde{\FF}\|_{L_2(\R^d)}
\cr
&\leq (C_2 + C_3){\mathfrak C}_{27} c(\varphi)^{3/2} \rho_0(\zeta)^{3/4} \eps \| {\FF}\|_{L_2(\O)}.
\end{aligned}
\eqno(13.31)
$$
Relations (13.30), (13.31), and (10.3) imply (13.24) with the constant
${\mathfrak C}_{24} = {\mathfrak C}_{27} (C_1+C_2+C_3)$.

Inequality (13.25) is deduced from (13.24) and (13.17), (13.18). $\bullet$

\smallskip
The following result is an analog of Theorem 8.7.

\smallskip\noindent\textbf{Theorem 13.6.} \textit{Suppose that the assumptions of Theorem} 13.2
\textit{are satisfied. Let $\O'$ be a strictly interior subdomain of the domain $\O$, and let $\delta = {\rm dist}\, \{\O'; \partial \O\}$.}
\textit{Then for $0< \varepsilon \leq \varepsilon_1$ we have}
$$
\| \u_\eps - \check{\v}_\eps \|_{H^1(\O')} \leq ({\mathcal C}_{23} \delta^{-1} +1)
({\mathfrak C}_{19} c(\varphi) \rho_0(\zeta)+
\widetilde{\mathfrak C}_{24} c(\varphi)^{5/2} \rho_0(\zeta)^{3/4})  \varepsilon \|\FF\|_{L_2(\O)},
\eqno(13.32)
$$
\textit{or, in operator terms},
$$
\begin{aligned}
&\| (\A_{N,\eps} - \zeta I)^{-1} - (\A_{N}^0  - \zeta I)^{-1} - \eps K_N^0(\eps;\zeta) \|_{L_2(\O) \to H^1(\O')}
\cr
&\leq ({\mathcal C}_{23} \delta^{-1} +1) ({\mathfrak C}_{19} c(\varphi) \rho_0(\zeta)+
\widetilde{\mathfrak C}_{24} c(\varphi)^{5/2} \rho_0(\zeta)^{3/4})  \varepsilon.
\end{aligned}
$$
\textit{For the flux $\p_\eps = g^\eps b(\D) \u_\eps$ we have}
$$
\begin{aligned}
&\| \p_\eps - \widetilde{g}^\eps b(\D) {\u}_0 \|_{L_2(\O')}
\cr
&\leq ({\mathcal C}_{23} \delta^{-1} +1) ({\mathfrak C}_{28}  c(\varphi) \rho_0(\zeta)+
{\mathfrak C}_{29} c(\varphi)^{5/2} \rho_0(\zeta)^{3/4})  \varepsilon \|\FF\|_{L_2(\O)}
\end{aligned}
\eqno(13.33)
$$
\textit{for $0< \eps \leq \eps_1$. The constants $\widetilde{\mathfrak C}_{24}$, ${\mathfrak C}_{28}$, and ${\mathfrak C}_{29}$ depend on}
$d$, $m$, $\alpha_0$, $\alpha_1$, $\| g \|_{L_\infty}$, $\| g^{-1} \|_{L_\infty}$,
\textit{the constants $k_1$, $k_2$ from} (9.2), \textit{the parameters of the lattice} $\Gamma$, \textit{the domain} $\O$,
\textit{and the norm} $\|\Lambda\|_{L_\infty}$.

\smallskip\noindent\textbf{Proof.}
Relations (4.3), (13.18), (13.21), and (13.24) imply (13.32) with $\widetilde{\mathfrak C}_{24} =
{\mathfrak C}_{24} + 2 {\mathcal C}''' C_\O^{(2)} c^\circ$.

Inequality (13.33) is deduced from (13.18), (13.23), and (13.32). $\bullet$

\section*{\S 14. Approximation of the resolvent $(\B_{N,\eps} - \zeta I)^{-1}$}

\smallskip\noindent\textbf{14.1. The operator $\B_{N,\eps}$.}
We denote
$$
Z = {\rm Ker}\, b(\D) = \{ \z \in H^1(\O;\C^n):\ b(\D)\z=0\}.
$$
From (9.2) it follows that
$$
\| \D \z \|_{L_2(\O)}^2 \leq k_1^{-1} k_2 \| \z \|^2_{L_2(\O)},\ \ \z \in Z.
$$
Due to compactness of the embedding of $H^1(\O;\C^n)$ into $L_2(\O;\C^n)$, this shows that
$Z$ is finite-dimensional.
Denote ${\rm dim}\,Z =p$. Obviously, $Z$ contains the $n$-dimensional subspace of constant $\C^n$-valued functions in $\O$.
We put
$$
{\mathcal H}(\O) = L_2(\O;\C^n) \ominus Z,\ \ H^1_\perp(\O;\C^n) = H^1(\O;\C^n) \cap {\mathcal H}(\O).
$$
As was checked in [Su3, Proposition 9.1], the form $\|b(\D)\u\|_{L_2(\O)}$ defines a norm in $H^1_\perp(\O;\C^n)$ equivalent to
the standard $H^1$-norm: there exists a constant $\widetilde{k}_1$ such that
$$
\widetilde{k}_1 \|\u\|^2_{H^1(\O)} \leq \|b(\D)\u\|_{L_2(\O)}^2,\ \ \u \in H^1_\perp(\O;\C^n).
\eqno(14.1)
$$

Let $\A_{N,\eps}$ be the operator in $L_2(\O;\C^n)$ generated by the form (9.3). Obviously,
${\rm Ker}\,\A_{N,\eps} = Z$. The orthogonal decomposition $L_2(\O;\C^n)= {\mathcal H}(\O) \oplus Z$
reduces the operator $\A_{N,\eps}$. Denote by $\B_{N,\eps}$ the part of $\A_{N,\eps}$
in the subspace ${\mathcal H}(\O)$. In other words,  $\B_{N,\eps}$ is the selfadjoint operator in
${\mathcal H}(\O)$ generated by the quadratic form
$$
b_{N,\eps}[\u,\u] = (g^\eps b(\D)\u, b(\D)\u)_{L_2(\O)},\ \ \u \in H^1_\perp(\O;\C^n).
$$
By (1.2), (1.5), and (14.1), we have
$$
\|g^{-1}\|^{-1}_{L_\infty}  \widetilde{k}_1 \|\u \|^2_{H^1(\O)}
\leq b_{N,\eps}[\u,\u] \leq \|g\|_{L_\infty} d \alpha_1 \| \D\u \|^2_{L_2(\O)},
\ \ \u \in H^1_\perp(\O;\C^n).
\eqno(14.2)
$$

Let $\A_N^0$ be the effective operator generated by the form (9.11) in $L_2(\O;\C^n)$.
Obviously, ${\rm Ker}\,\A_N^0 =Z$. The decomposition $L_2(\O;\C^n)= {\mathcal H}(\O) \oplus Z$
reduces the operator $\A_{N}^0$. Let $\B_{N}^0$ be the part of $\A_{N}^0$
in ${\mathcal H}(\O)$. In other words,  $\B_{N}^0$ is the operator generated in
${\mathcal H}(\O)$ by the quadratic form
$$
b_{N}^0[\u,\u] = (g^0 b(\D)\u, b(\D)\u)_{L_2(\O)},\ \ \u \in H^1_\perp(\O;\C^n).
$$
Similarly to (14.2), we have
$$
\|g^{-1}\|^{-1}_{L_\infty}  \widetilde{k}_1 \|\u \|^2_{H^1(\O)}
\leq b_{N}^0[\u,\u] \leq \|g\|_{L_\infty} d \alpha_1 \|\D\u \|^2_{L_2(\O)},
\ \ \u \in H^1_\perp(\O;\C^n).
\eqno(14.3)
$$
\textit{Let  ${\mathcal P}$ denote the orthogonal projection of $L_2(\O;\C^n)$ onto ${\mathcal H}(\O)$,
and let  ${\mathcal P}_Z$ denote the orthogonal projection onto} $Z$; then ${\mathcal P}= I - {\mathcal P}_Z$.

Let $\zeta \in \C \setminus [c_\flat,\infty)$, where $c_\flat >0$ is a common lower bound of the operators
$\B_{N,\eps}$ and $\B_N^0$. In other words, $0< c_\flat \leq \min \{\lambda_{2,\eps}(N), \lambda_2^0(N)\}$,
where $\lambda_{2,\eps}(N)$ (respectively, $\lambda_{2}^0(N)$) is the first nonnegative eigenvalue of the operator
$\A_{N,\eps}$ (respectively, $\A_N^0$). (Counting multiplicities, they are $(p+1)$-th eigenvalues.)

\smallskip\noindent\textbf{Remark 14.1.} 1) By (14.2) and (14.3), one can take $c_\flat$ equal to
$\|g^{-1}\|^{-1}_{L_\infty} \widetilde{k}_1$.
2) Let $\nu >0$ be arbitrarily small. If $\eps$ is sufficiently small, one can take
$c_\flat$ equal to $c_\flat = \lambda_2^0(N) - \nu$.
3)~It is easy to find an upper bound for $c_\flat$: from (14.2), (14.3), and the variational principle it follows that
$c_\flat \leq \|g\|_{L_\infty} d \alpha_1 \mu^0_{p+1}(N)$, where $\mu^0_{p+1}(N)>0$ is the $(p+1)$-th eigenvalue of the operator $-\Delta$
with the Neumann condition in $L_2(\O;\C^n)$ (the eigenvalues are enumerated in the increasing order with multiplicities taken into account).
Thus, $c_\flat$ does not exceed a number depending only on $d$, $n$, $p$,
$\alpha_1$, $\|g\|_{L_\infty}$, and the domain $\O$.

Let ${\boldsymbol \varphi}_\eps= (\B_{N,\eps} -\zeta I)^{-1} \FF$, where $\FF \in {\mathcal H}(\O)$,
i.~e., ${\boldsymbol \varphi}_\eps$ is the generalized solution of the Neumann problem
$$
\begin{aligned}
&b(\D)^* g^\eps b(\D){\boldsymbol \varphi}_\eps - \zeta {\boldsymbol \varphi}_\eps = \FF \ \ \textrm{in}\ \O;
\cr
&\partial^\eps_{\boldsymbol \nu} {\boldsymbol \varphi}_\eps \vert_{\partial \O} =0;\ \
({\boldsymbol \varphi}_\eps, \z)_{L_2(\O)}=0\ \ \forall \ \z \in Z.
\end{aligned}
\eqno(14.4)
$$
Here the condition that ${\boldsymbol \varphi}_\eps$ is orthogonal to $Z$ is valid automatically if $\zeta \ne 0$,
and for $\zeta=0$ this is an additional condition.

Let ${\boldsymbol \varphi}_0 = (\B^0_N - \zeta I)^{-1}\FF$, $\FF \in {\mathcal H}(\O)$.
Then ${\boldsymbol \varphi}_0$ is the generalized solution of the Neumann problem
$$
\begin{aligned}
&b(\D)^* g^0 b(\D){\boldsymbol \varphi}_0 - \zeta {\boldsymbol \varphi}_0 = \FF \ \ \textrm{in}\ \O;
\cr
&\partial^0_{\boldsymbol \nu} {\boldsymbol \varphi}_0 \vert_{\partial \O} =0;\ \
({\boldsymbol \varphi}_0, \z)_{L_2(\O)}=0\ \ \forall \ \z \in Z.
\end{aligned}
\eqno(14.5)
$$
Here also the orthogonality condition is valid automatically if
$\zeta \ne 0$, and for $\zeta=0$ this is an additional condition.

Denote
$$
{\mathcal K}_N(\eps;\zeta) := R_\O [\Lambda^\eps] S_\eps b(\D) P_\O (\B_{N,\eps} -\zeta I)^{-1}
$$
and put $ \widetilde{\boldsymbol \varphi}_0 = P_\O {\boldsymbol \varphi}_0$,
$$
 {\boldsymbol \psi}_\eps  ={\boldsymbol \varphi}_0 + \eps \Lambda^\eps S_\eps b(\D) \widetilde{\boldsymbol \varphi}_0
= (\B^0_N - \zeta I)^{-1}\FF + \eps {\mathcal K}_N(\eps;\zeta) \FF.
\eqno(14.6)
$$

\smallskip\noindent\textbf{Theorem 14.2.} \textit{Let $\zeta \in \C \setminus [c_\flat,\infty)$,
where $c_\flat >0$ is a common lower bound of the operators $\B_{N,\eps}$ and $\B_N^0$.
We put $\zeta - c_\flat = |\zeta - c_\flat|e^{i \vartheta}$ and denote
$$
\rho_\flat(\zeta) =
\begin{cases}
c(\vartheta)^2 |\zeta - c_\flat |^{-2}, & |\zeta - c_\flat|<1, \\
c(\vartheta)^2, & |\zeta - c_\flat | \ge 1.
\end{cases}
\eqno(14.7)
$$
Let ${\boldsymbol \varphi}_\eps$ be the solution of problem} (14.4), \textit{and let ${\boldsymbol \varphi}_0$
be the solution of problem} (14.5) \textit{for} $\FF \in {\mathcal H}(\O)$.
\textit{Suppose that ${\boldsymbol \psi}_\eps$ is given by} (14.6). \textit{Suppose that the number $\eps_1$
satisfies Condition} 4.1. \textit{Then for $0< \eps \leq \eps_1$ we have}
$$
\| {\boldsymbol \varphi}_\eps - {\boldsymbol \varphi}_0 \|_{L_2(\O)}
\leq {\mathfrak C}_{30} \rho_\flat(\zeta) \eps \| \FF \|_{L_2(\O)},
$$
$$
\| {\boldsymbol \varphi}_\eps - {\boldsymbol \psi}_\eps \|_{H^1(\O)} \leq {\mathfrak C}_{31}
\rho_\flat(\zeta) \eps^{1/2} \|\FF\|_{L_2(\O)},
\eqno(14.8)
$$
\textit{or, in operator terms},
$$
\| (\B_{N,\eps} - \zeta I)^{-1} - (\B_{N}^0  - \zeta I)^{-1} \|_{{\mathcal H}(\O) \to {\mathcal H}(\O)}
\leq {\mathfrak C}_{30} \rho_\flat(\zeta) \eps,
\eqno(14.9)
$$
$$
\| (\B_{N,\eps} - \zeta I)^{-1} - (\B_{N}^0  - \zeta I)^{-1} - \eps {\mathcal K}_N(\eps;\zeta)
\|_{{\mathcal H}(\O) \to H^1(\O)} \leq {\mathfrak C}_{31} \rho_\flat(\zeta) \eps^{1/2}.
\eqno(14.10)
$$
\textit{For the flux $g^\eps b(\D) {\boldsymbol \varphi}_\eps$ we have}
$$
\| g^\eps b(\D) {\boldsymbol \varphi}_\eps - \widetilde{g}^\eps S_\eps
b(\D) \widetilde{{\boldsymbol \varphi}}_0 \|_{L_2(\O)}
\leq {\mathfrak C}_{32} \rho_\flat(\zeta) \eps^{1/2} \|\FF\|_{L_2(\O)}
\eqno(14.11)
$$
\textit{for $0< \eps \leq \eps_1$. The constants ${\mathfrak C}_{30}$, ${\mathfrak C}_{31}$, and ${\mathfrak C}_{32}$ depend on}
$d$, $n$, $m$, $p$, $\alpha_0$, $\alpha_1$, $\| g \|_{L_\infty}$, $\| g^{-1} \|_{L_\infty}$,
\textit{the constants $k_1$, $k_2$ from} (9.2), \textit{the constant $\widetilde{k}_1$ from} (14.1), \textit{the parameters
of the lattice} $\Gamma$, \textit{and the domain} $\O$.

\smallskip\noindent\textbf{Proof.} Applying Theorem 10.1 with $\zeta=-1$, for $0< \eps \leq \eps_1$ we obtain
$$
\begin{aligned}
&\| (\B_{N,\eps} + I)^{-1} - (\B_{N}^0  + I)^{-1} \|_{{\mathcal H}(\O) \to {\mathcal H}(\O)}
\cr
&= \| \left((\A_{N,\eps} + I)^{-1} - (\A_{N}^0  + I)^{-1}\right){\mathcal P}
\|_{L_2(\O) \to L_2(\O)}
\leq {\mathfrak C}_{1} (\eps + \eps^2) \leq 2 {\mathfrak C}_{1} \eps.
\end{aligned}
$$
Then, using the analog of identity (2.10) for the operators ${\mathcal B}_{N,\eps}$
and ${\mathcal B}_{N}^0$, we arrive at
$$
\| (\B_{N,\eps} - \zeta I)^{-1} - (\B_{N}^0 -\zeta I)^{-1} \|_{{\mathcal H}(\O) \to {\mathcal H}(\O)}
\leq 2 {\mathfrak C}_1 \eps \sup_{x\geq c_\flat} (x+1)^2 |x-\zeta|^{-2}
\eqno(14.12)
$$
for $0< \eps \leq \eps_1$. Similarly to (8.8),
$$
\sup_{x\geq c_\flat} (x+1)^2 |x-\zeta|^{-2} \leq \check{c}_\flat \rho_\flat(\zeta),
\ \ \zeta \in \C \setminus [c_\flat,\infty),
\eqno(14.13)
$$
where $\check{c}_\flat = (c_\flat +2)^2$.
By Remark 14.1(3), $\check{c}_\flat$ is bounded by a number depending only on $d$, $n$, $p$, $\alpha_1$,
$\|g\|_{L_\infty}$, and the domain $\O$. Relations (14.12) and (14.13) imply (14.9) with ${\mathfrak C}_{30}=
2 {\mathfrak C}_1 \check{c}_\flat$.

Now we apply Theorem 10.2 with $\zeta = -1$. For $0< \eps \leq \eps_1$ we have:
$$
\| (\A_{N,\eps} + I)^{-1} - (\A_{N}^0  + I)^{-1} -\eps K_N(\eps;-1)
\|_{L_2(\O) \to H^1(\O)}
\leq ({\mathfrak C}_{2} + {\mathfrak C}_{3}) \eps^{1/2}.
\eqno(14.14)
$$
It follows that
$$
\begin{aligned}
&\| \A_{N,\eps}^{1/2} \left((\A_{N,\eps} + I)^{-1} - (\A_{N}^0  + I)^{-1} -\eps K_N(\eps;-1)\right)
\|_{L_2(\O) \to L_2(\O)}
\cr
&\leq \|g\|_{L_\infty}^{1/2} (d \alpha_1)^{1/2} ({\mathfrak C}_{2} + {\mathfrak C}_{3}) \eps^{1/2},
\ \ 0< \eps \leq \eps_1.
\end{aligned}
\eqno(14.15)
$$
Multiplying the operator under the norm-sign in (14.15) by the projection  $\mathcal P$ from both sides, we obtain
$$
\begin{aligned}
&\| \B_{N,\eps}^{1/2} \left((\B_{N,\eps} + I)^{-1} - (\B_{N}^0  + I)^{-1} -\eps {\mathcal P}
{\mathcal K}_N(\eps;-1)\right)
\|_{{\mathcal H}(\O) \to {\mathcal H}(\O)}
\cr
&\leq \|g\|_{L_\infty}^{1/2} (d \alpha_1)^{1/2} ({\mathfrak C}_{2} + {\mathfrak C}_{3}) \eps^{1/2},
\ \ 0< \eps \leq \eps_1.
\end{aligned}
\eqno(14.16)
$$
We have
$$
\begin{aligned}
&(\B_{N,\eps} -\zeta I)^{-1} - (\B_{N}^0 -\zeta I)^{-1} -\eps {\mathcal P} {\mathcal K}_N(\eps;\zeta)
\cr
&=
(\B_{N,\eps} + I)(\B_{N,\eps} -\zeta I)^{-1}
\left((\B_{N,\eps} + I)^{-1} - (\B_{N}^0  + I)^{-1} -\eps {\mathcal P} {\mathcal K}_N(\eps;-1)\right)
\cr
&\times(\B_{N}^0 + I)(\B_{N}^0 -\zeta I)^{-1}
+ (\zeta +1) (\B_{N,\eps} -\zeta I)^{-1} \eps {\mathcal P}
{\mathcal K}_N(\eps;\zeta).
\end{aligned}
\eqno(14.17)
$$
Multiplying (14.17) by $\B_{N,\eps}^{1/2}$ from the left and using (14.16), we arrive at
$$
\begin{aligned}
&\|\B_{N,\eps}^{1/2} \left((\B_{N,\eps} -\zeta I)^{-1} - (\B_{N}^0 -\zeta I)^{-1} -\eps {\mathcal P} {\mathcal K}_N(\eps;\zeta)\right) \|_{{\mathcal H}(\O) \to {\mathcal H}(\O)}
\cr
&\leq \|g\|_{L_\infty}^{1/2} (d \alpha_1)^{1/2} ({\mathfrak C}_{2} + {\mathfrak C}_{3}) \eps^{1/2}
\sup_{x \ge c_\flat} (x+1)^2 |x-\zeta|^{-2}
\cr
&+ |\zeta +1| \eps
\sup_{x \ge c_\flat} x |x-\zeta|^{-2}
\|\Lambda^\eps S_\eps b(\D) P_\O (\B_N^0)^{-1/2}\|_{{\mathcal H}(\O) \to L_2(\R^d)}.
 \end{aligned}
 \eqno(14.18)
 $$
Denote the summands in the right-hand side of (14.18) by ${\mathcal T}_1(\eps)$ and ${\mathcal T}_2(\eps)$.
Taking (14.13) into account, we have
$$
{\mathcal T}_1(\eps) \leq {\mathfrak C}_{33} \eps^{1/2} \rho_\flat(\zeta),
\ \ 0< \eps \leq \eps_1,
\eqno(14.19)
$$
where ${\mathfrak C}_{33}= \check{c}_\flat \|g\|_{L_\infty} (d \alpha_1)^{1/2} ({\mathfrak C}_{2} + {\mathfrak C}_{3})$.
Next, from (14.3) it follows that
$$
\| (\B_N^0)^{-1/2} \|_{{\mathcal H}(\O) \to H^1(\O)} \leq \|g^{-1}\|^{1/2}_{L_\infty} \widetilde{k}_1^{-1/2}.
$$
Combining this with (1.4), (1.19), and (4.3), we obtain
$$
\|\Lambda^\eps S_\eps b(\D) P_\O (\B_N^0)^{-1/2}\|_{{\mathcal H}(\O) \to L_2(\R^d)} \leq
M_1 \alpha_1^{1/2} C_\O^{(1)} \|g^{-1}\|^{1/2}_{L_\infty} \widetilde{k}_1^{-1/2}.
\eqno(14.20)
$$
By (14.20), the second term in (14.18) admits the estimate
$$
{\mathcal T}_2(\eps) \leq \eps |\zeta +1|
M_1 \alpha_1^{1/2} C_\O^{(1)} \|g^{-1}\|^{1/2}_{L_\infty} \widetilde{k}_1^{-1/2}
\sup_{x \ge c_\flat} x |x-\zeta|^{-2}.
\eqno(14.21)
$$
Similarly to (8.17), we have
$$
|\zeta +1| \sup_{x \ge c_\flat} x |x-\zeta|^{-2} \leq (c_\flat +2) (c_\flat +1) \rho_\flat(\zeta),
\ \ \zeta \in \C \setminus [c_\flat,\infty).
$$
Together with (14.21) this implies
$$
{\mathcal T}_2(\eps) \leq {\mathfrak C}_{34} \eps \rho_\flat(\zeta),
\eqno(14.22)
$$
where ${\mathfrak C}_{34} = (c_\flat +2) (c_\flat +1) M_1 \alpha_1^{1/2} C_\O^{(1)} \|g^{-1}\|^{1/2}_{L_\infty} \widetilde{k}_1^{-1/2}.$

As a result, relations (14.18), (14.19), and (14.22) imply that
$$
\begin{aligned}
&\|\B_{N,\eps}^{1/2} \left((\B_{N,\eps} -\zeta I)^{-1} - (\B_{N}^0 -\zeta I)^{-1} -\eps {\mathcal P} {\mathcal K}_N(\eps;\zeta)\right) \|_{{\mathcal H}(\O) \to {\mathcal H}(\O)}
\cr
&\leq ({\mathfrak C}_{33}+{\mathfrak C}_{34}) \eps^{1/2} \rho_\flat(\zeta),
\ \ 0< \eps \leq \eps_1.
\end{aligned}
$$
Together with (14.2) this yields
$$
\| (\B_{N,\eps} -\zeta I)^{-1} - (\B_{N}^0 -\zeta I)^{-1} -\eps {\mathcal P} {\mathcal K}_N(\eps;\zeta)
 \|_{{\mathcal H}(\O) \to H^1(\O)}
\leq {\mathfrak C}_{35} \eps^{1/2} \rho_\flat(\zeta)
\eqno(14.23)
$$
for $0< \eps \leq \eps_1$, where ${\mathfrak C}_{35}= \|g^{-1}\|^{1/2}_{L_\infty} \widetilde{k}_1^{-1/2}
({\mathfrak C}_{33}+{\mathfrak C}_{34})$.

Now we show that the last term under the norm-sign in (14.23) can be replaced by
$\eps {\mathcal K}_N(\eps;\zeta)$; this will only change a constant in estimate. Multiplying
the operator under the norm-sign in (14.14) by $\mathcal P$ from the right, we obtain
$$
\| (\B_{N,\eps} + I)^{-1} - (\B_{N}^0  + I)^{-1} -\eps {\mathcal K}_N(\eps;-1)
\|_{{\mathcal H}(\O) \to H^1(\O)}
\leq ({\mathfrak C}_{2} + {\mathfrak C}_{3}) \eps^{1/2}
\eqno(14.24)
$$
for $0< \eps \leq \eps_1$.
On the other side, from (14.16) and (14.2) it follows that
$$
\| (\B_{N,\eps} + I)^{-1} - (\B_{N}^0  + I)^{-1} -\eps {\mathcal P}{\mathcal K}_N(\eps;-1)
\|_{{\mathcal H}(\O) \to H^1(\O)}
\leq {\mathfrak C}_{36} \eps^{1/2}
\eqno(14.25)
$$
for $0< \eps \leq \eps_1$, where ${\mathfrak C}_{36} = \widetilde{k}_1^{-1/2} \|g^{-1}\|^{1/2}_{L_\infty} \|g\|^{1/2}_{L_\infty}
(d \alpha_1)^{1/2} ({\mathfrak C}_2 + {\mathfrak C}_3)$. Comparing (14.24) and (14.25), we see that
$$
\eps \| {\mathcal P}_Z{\mathcal K}_N(\eps;-1)
\|_{{\mathcal H}(\O) \to H^1(\O)}
\leq ({\mathfrak C}_{2} + {\mathfrak C}_{3} +{\mathfrak C}_{36}) \eps^{1/2},\ \ 0< \eps \leq \eps_1.
$$
Hence, by (14.13), we have
$$
\begin{aligned}
&\eps \| {\mathcal P}_Z {\mathcal K}_N(\eps;\zeta)\|_{{\mathcal H}(\O) \to H^1(\O)}\leq
\eps \| {\mathcal P}_Z {\mathcal K}_N(\eps;-1)\|_{{\mathcal H}(\O) \to H^1(\O)}
\cr
&\times \| ({\mathcal B}_N^0+I) ({\mathcal B}_N^0-\zeta I)^{-1}\|_{{\mathcal H}(\O) \to {\mathcal H}(\O)}
\leq ({\mathfrak C}_{2} + {\mathfrak C}_{3} +{\mathfrak C}_{36})
\check{c}_\flat^{1/2} \eps^{1/2} \rho_\flat(\zeta)^{1/2}
\end{aligned}
\eqno(14.26)
$$
for $0< \eps \leq \eps_1$.
As a result, relations (14.23) and (14.26) together with the identity ${\mathcal P}+ {\mathcal P}_Z=I$
imply (14.10) with ${\mathfrak C}_{31}= {\mathfrak C}_{35}+ ({\mathfrak C}_{2} + {\mathfrak C}_{3} +{\mathfrak C}_{36})
\check{c}_\flat^{1/2}$.

It remains to check (14.11). From (14.8) and (1.2), (1.5) it follows that for $0< \eps \leq \eps_1$ we have
$$
\| g^\eps b(\D) {\boldsymbol \varphi}_\eps - g^\eps b(\D) {\boldsymbol \psi}_\eps\|_{L_2(\O)}
\leq \|g\|_{L_\infty} (d \alpha_1)^{1/2} {\mathfrak C}_{31} \rho_\flat(\zeta) \eps^{1/2} \|\FF\|_{L_2(\O)}.
\eqno(14.27)
$$
Similarly to (5.19)--(5.21), 
$$
\|  g^\eps b(\D) {\boldsymbol \psi}_\eps - \widetilde{g}^\eps S_\eps b(\D)
\widetilde{\boldsymbol \varphi}_0 \|_{L_2(\O)}
\leq  ({\mathcal C}'+{\mathcal C}'')  \eps \|\widetilde{\boldsymbol \varphi}_0 \|_{H^2(\R^d)}.
\eqno(14.28)
$$
Now we estimate the $H^2(\O)$-norm of the function ${\boldsymbol \varphi}_0 = (\B_N^0 - \zeta I)^{-1} \FF$.
By (9.12) and (14.13), we have
$$
\begin{aligned}
&\|(\B_N^0 - \zeta I)^{-1}\|_{{\mathcal H}(\O) \to H^2(\O)} \leq
\|(\B_N^0 + I)^{-1}\|_{{\mathcal H}(\O) \to H^2(\O)}
\cr
&\times \|(\B_N^0 + I) (\B_N^0 - \zeta I)^{-1}\|_{{\mathcal H}(\O) \to {\mathcal H}(\O)}
\leq c^\circ \check{c}_\flat^{1/2} \rho_\flat(\zeta)^{1/2}.
\end{aligned}
$$
Hence,
$$
\|{\boldsymbol \varphi}_0\|_{H^2(\O)} \leq
c^\circ \check{c}_\flat^{1/2} \rho_\flat(\zeta)^{1/2} \|\FF\|_{L_2(\O)}.
\eqno(14.29)
$$
By (4.3), (14.28), and (14.29),
$$
\|  g^\eps b(\D) {\boldsymbol \psi}_\eps - \widetilde{g}^\eps S_\eps b(\D)
\widetilde{\boldsymbol \varphi}_0 \|_{L_2(\O)}
\leq  {\mathfrak C}_{37} \eps \rho_\flat(\zeta)^{1/2} \|\FF\|_{L_2(\O)},
\eqno(14.30)
$$
where ${\mathfrak C}_{37}= ({\mathcal C}'+{\mathcal C}'') C_{\O}^{(2)} c^\circ \check{c}_\flat^{1/2}$.
Relations (14.27) and (14.30) imply (14.11) with ${\mathfrak C}_{32}= \|g\|_{L_\infty} (d\alpha_1)^{1/2}
{\mathfrak C}_{31} + {\mathfrak C}_{37}$. $\bullet$

\smallskip
Now we distinguish the case where $\Lambda \in L_\infty$. Denote
$$
{\mathcal K}_N^0(\eps;\zeta) := [\Lambda^\eps] b(\D) (\B_N^0 - \zeta I)^{-1}
\eqno(14.31)
$$
and put
$$
\check{\boldsymbol \psi}_\eps = {\boldsymbol \varphi}_0 + \eps \Lambda^\eps b(\D) {\boldsymbol \varphi}_0
= (\B_N^0 -\zeta I)^{-1}\FF + \eps {\mathcal K}_N^0(\eps;\zeta) \FF.
\eqno(14.32)
$$

\smallskip\noindent\textbf{Theorem 14.3.} \textit{Suppose that the assumptions of Theorem} 14.2
\textit{and Condition} 2.8 \textit{are satisfied. Let $\check{\boldsymbol \psi}_\eps$ be given by} (14.32).
\textit{Then for $0< \eps \leq \eps_1$ we have}
$$
\| {\boldsymbol \varphi}_\eps - \check{\boldsymbol \psi}_\eps \|_{H^1(\O)} \leq {\mathfrak C}_{31}^\circ
\rho_\flat(\zeta) \eps^{1/2} \|\FF\|_{L_2(\O)},
\eqno(14.33)
$$
\textit{or, in operator terms},
$$
\| (\B_{N,\eps} - \zeta I)^{-1} - (\B_{N}^0  - \zeta I)^{-1} - \eps {\mathcal K}_N^0(\eps;\zeta)
\|_{{\mathcal H}(\O) \to H^1(\O)} \leq {\mathfrak C}_{31}^\circ \rho_\flat(\zeta) \eps^{1/2}.
\eqno(14.34)
$$
\textit{For the flux $g^\eps b(\D) {\boldsymbol \varphi}_\eps$ we have}
$$
\| g^\eps b(\D) {\boldsymbol \varphi}_\eps -
\widetilde{g}^\eps b(\D) {\boldsymbol \varphi}_0 \|_{L_2(\O)}
\leq {\mathfrak C}_{32}^\circ \rho_\flat(\zeta) \eps^{1/2} \|\FF\|_{L_2(\O)}
\eqno(14.35)
$$
\textit{for $0< \eps \leq \eps_1$. The constants ${\mathfrak C}_{31}^\circ$ and ${\mathfrak C}_{32}^\circ$ depend on}
$d$, $n$, $m$, $p$, $\alpha_0$, $\alpha_1$, $\| g \|_{L_\infty}$, $\| g^{-1} \|_{L_\infty}$,
\textit{the constants $k_1$, $k_2$ from} (9.2), \textit{the constant $\widetilde{k}_1$ from} (14.1), \textit{the parameters of the lattice} $\Gamma$, \textit{the domain $\O$, and the norm} $\|\Lambda\|_{L_\infty}$.

\smallskip\noindent\textbf{Proof.} Similarly to (6.5)--(6.9),
$$
\| {\boldsymbol \psi}_\eps - \check{\boldsymbol \psi}_\eps \|_{H^1(\O)} \leq {\mathcal C}'''
\eps \| \widetilde{\boldsymbol \varphi}_0 \|_{H^2(\R^d)}.
$$
Together with (4.3) and (14.29) this yields
$$
\| {\boldsymbol \psi}_\eps - \check{\boldsymbol \psi}_\eps \|_{H^1(\O)} \leq
{\mathcal C}'''C_{\O}^{(2)} c^\circ \check{c}_\flat^{1/2} \eps \rho_\flat(\zeta)^{1/2} \| \FF \|_{L_2(\O)}.
\eqno(14.36)
$$
Now relations (14.8) and (14.36) imply the required estimate (14.33) with
${\mathfrak C}^\circ_{31}= {\mathfrak C}_{31} + {\mathcal C}'''C_{\O}^{(2)} c^\circ \check{c}_\flat^{1/2}$.

Let us check (14.35). From (14.33) it follows that
$$
\| g^\eps b(\D){\boldsymbol \varphi}_\eps -
g^\eps b(\D) \check{\boldsymbol \psi}_\eps \|_{L_2(\O)} \leq
\|g\|_{L_\infty} (d \alpha_1)^{1/2}{\mathfrak C}_{31}^\circ
\rho_\flat(\zeta) \eps^{1/2} \|\FF\|_{L_2(\O)}
\eqno(14.37)
$$
for $0< \eps \leq \eps_1$. Similarly to  (6.11) and (6.12), we have
$$
\| g^\eps b(\D) \check{\boldsymbol \psi}_\eps -
\widetilde{g}^\eps b(\D) {\boldsymbol \varphi}_0
\|_{L_2(\O)} \leq    \widetilde{\mathcal C}' \eps  \|{\boldsymbol \varphi}_0\|_{H^2(\O)}.
\eqno(14.38)
$$
Relations (14.29), (14.37), and (14.38) imply (14.35) with
${\mathfrak C}_{32}^\circ = \|g\|_{L_\infty} (d\alpha_1)^{1/2}{\mathfrak C}_{31}^\circ
+\widetilde{\mathcal C}' c^\circ \check{c}_\flat^{1/2}$. $\bullet$

\smallskip\noindent\textbf{14.3. Special cases.} The following statement concerns the case where the corrector is equal to zero.

\smallskip\noindent\textbf{Proposition 14.4.} \textit{Suppose that the assumptions of Theorem} 14.2
\textit{are satisfied. If $g^0 = \overline{g}$, i.~e., relations} (1.13) \textit{are satisfied, then
$\Lambda=0$, ${\boldsymbol \psi}_\eps = {\boldsymbol \varphi}_0$, and for $0< \eps \leq \eps_1$ we have}
$$
\| {\boldsymbol \varphi}_\eps - {\boldsymbol \varphi}_0 \|_{H^1(\O)} \leq
{\mathfrak C}_{31} \rho_\flat(\zeta) \eps^{1/2} \| \FF \|_{L_2(\O)}.
$$

\smallskip
Next statement is similar to Proposition 8.5.

\smallskip\noindent\textbf{Proposition 14.5.} \textit{Suppose that the assumptions of Theorem} 14.2
\textit{are satisfied. If $g^0 = \underline{g}$, i.~e., relations} (1.14) \textit{are satisfied, then
for $0< \eps \leq \eps_1$ we have}
$$
\| g^\eps b(\D) {\boldsymbol \varphi}_\eps - g^0 b(\D) {\boldsymbol \varphi}_0
\|_{L_2(\O)} \leq {\mathfrak C}_{32}^\circ \rho_\flat(\zeta) \eps^{1/2} \| \FF \|_{L_2(\O)}.
$$

\smallskip\noindent\textbf{14.4. Approximation of the resolvent of the operator $\B_{N,\eps}$ in a strictly interior subdomain.}
Let $\O'$ be a strictly interior subdomain of the domain $\O$.
Using Theorem 14.2 and the results for the problem in $\R^d$, we obtain approximation
of the solution ${\boldsymbol \varphi}_\eps$ in $H^1(\O')$ of sharp order with respect to $\eps$.

\smallskip\noindent\textbf{Theorem 14.6.} \textit{Suppose that the assumptions of Theorem} 14.2
\textit{are satisfied. Let $\O'$ be a strictly interior subdomain of the domain $\O$, and let $\delta :=
{\rm dist}\,\{ \O'; \partial \O \}$. Then for $0< \eps \leq \eps_1$ we have}
$$
\begin{aligned}
&\| {\boldsymbol \varphi}_\eps - {\boldsymbol \psi}_\eps \|_{H^1(\O')}
\cr
&\leq
\left( {\mathfrak C}_{\flat}' \delta^{-1} (c(\vartheta) \rho_\flat(\zeta) +
c(\vartheta)^{5/2} \rho_\flat(\zeta)^{3/4}) + {\mathfrak C}_{\flat}''
c(\vartheta)^{1/2} \rho_\flat(\zeta)^{5/4} \right)
\eps \|\FF\|_{L_2(\O)},
\end{aligned}
$$
\textit{or, in operator terms},
$$
\begin{aligned}
&\| (\B_{N,\eps} - \zeta I)^{-1} - (\B_{N}^0  - \zeta I)^{-1} - \eps {\mathcal K}_N(\eps;\zeta)
\|_{{\mathcal H}(\O) \to H^1(\O')}
\cr
&\leq
\left( {\mathfrak C}_{\flat}' \delta^{-1} (c(\vartheta) \rho_\flat(\zeta) +
c(\vartheta)^{5/2} \rho_\flat(\zeta)^{3/4}) + {\mathfrak C}_{\flat}''
c(\vartheta)^{1/2} \rho_\flat(\zeta)^{5/4} \right) \eps.
\end{aligned}
\eqno(14.39)
$$
\textit{For the flux $g^\eps b(\D) {\boldsymbol \varphi}_\eps$ we have}
$$
\begin{aligned}
&\| g^\eps b(\D) {\boldsymbol \varphi}_\eps -
\widetilde{g}^\eps S_\eps b(\D) \widetilde{\boldsymbol \varphi}_0 \|_{L_2(\O')}
\cr
&\leq \left( \widetilde{\mathfrak C}_{\flat}' \delta^{-1} (c(\vartheta) \rho_\flat(\zeta) +
c(\vartheta)^{5/2} \rho_\flat(\zeta)^{3/4}) + \widetilde{\mathfrak C}_{\flat}''
c(\vartheta)^{1/2} \rho_\flat(\zeta)^{5/4} \right)
\eps \|\FF\|_{L_2(\O)}
\end{aligned}
\eqno(14.40)
$$
\textit{for $0< \eps \leq \eps_1$. The constants ${\mathfrak C}_{\flat}'$, ${\mathfrak C}_{\flat}''$,
$\widetilde{\mathfrak C}_{\flat}'$, and $\widetilde{\mathfrak C}_{\flat}''$ depend on}
$d$, $n$, $m$, $p$, $\alpha_0$, $\alpha_1$, $\| g \|_{L_\infty}$, $\| g^{-1} \|_{L_\infty}$,
\textit{the constants $k_1$, $k_2$ from} (9.2), \textit{the constant $\widetilde{k}_1$ from} (14.1),
\textit{the parameters of the lattice} $\Gamma$, \textit{and the domain} $\O$.

\smallskip\noindent\textbf{Proof} is quite similar to the proof of Theorem 8.6; we omit the details.
$\bullet$

\smallskip\noindent\textbf{Theorem 14.7.} \textit{Suppose that the assumptions of Theorem} 14.3
\textit{are satisfied. Let $\O'$ be a strictly interior subdomain of the domain $\O$, and let $\delta :=
{\rm dist}\,\{ \O'; \partial \O \}$. Then for $0< \eps \leq \eps_1$ we have}
$$
\begin{aligned}
&\| {\boldsymbol \varphi}_\eps - \check{\boldsymbol \psi}_\eps \|_{H^1(\O')}
\cr
&\leq
\left( {\mathfrak C}_{\flat}' \delta^{-1} (c(\vartheta) \rho_\flat(\zeta) +
c(\vartheta)^{5/2} \rho_\flat(\zeta)^{3/4}) + \check{\mathfrak C}_{\flat}''
c(\vartheta)^{1/2} \rho_\flat(\zeta)^{5/4} \right)
\eps \|\FF\|_{L_2(\O)},
\end{aligned}
$$
\textit{or, in operator terms},
$$
\begin{aligned}
&\| (\B_{N,\eps} - \zeta I)^{-1} - (\B_{N}^0  - \zeta I)^{-1} - \eps {\mathcal K}_N^0(\eps;\zeta)
\|_{{\mathcal H}(\O) \to H^1(\O')}
\cr
&\leq
\left( {\mathfrak C}_{\flat}' \delta^{-1} (c(\vartheta) \rho_\flat(\zeta) +
c(\vartheta)^{5/2} \rho_\flat(\zeta)^{3/4}) + \check{\mathfrak C}_{\flat}''
c(\vartheta)^{1/2} \rho_\flat(\zeta)^{5/4} \right) \eps.
\end{aligned}
\eqno(14.41)
$$
\textit{For the flux $g^\eps b(\D) {\boldsymbol \varphi}_\eps$ we have}
$$
\begin{aligned}
&\| g^\eps b(\D) {\boldsymbol \varphi}_\eps -
\widetilde{g}^\eps b(\D) {\boldsymbol \varphi}_0 \|_{L_2(\O')}
\cr
&\leq \left( \widetilde{\mathfrak C}_{\flat}' \delta^{-1} (c(\vartheta) \rho_\flat(\zeta) +
c(\vartheta)^{5/2} \rho_\flat(\zeta)^{3/4}) + \widehat{\mathfrak C}_{\flat}''
c(\vartheta)^{1/2} \rho_\flat(\zeta)^{5/4} \right)
\eps \|\FF\|_{L_2(\O)}
\end{aligned}
\eqno(14.42)
$$
\textit{for $0< \eps \leq \eps_1$. The constants ${\mathfrak C}_{\flat}'$, $\widetilde{\mathfrak C}_{\flat}'$ are the same as in
Theorem} 14.6.
\textit{The constants  $\check{\mathfrak C}_{\flat}''$ and $\widehat{\mathfrak C}_{\flat}''$ depend on}
$d$, $n$, $m$, $p$, $\alpha_0$, $\alpha_1$, $\| g \|_{L_\infty}$, $\| g^{-1} \|_{L_\infty}$,
\textit{the constants $k_1$, $k_2$ from} (9.2), \textit{the constant $\widetilde{k}_1$ from} (14.1),
\textit{the parameters of the lattice} $\Gamma$, \textit{the domain $\O$, and} $\|\Lambda\|_{L_\infty}$.

\smallskip\noindent\textbf{Proof.} It is easy to prove this theorem, using Theorem 14.6 and relations (14.29), (14.36), and (14.38). $\bullet$

\smallskip\noindent\textbf{14.5. Application of the resuts for $\B_{N,\eps}$ to the operator $\A_{N,\eps}$.}
Theorem 14.2 allows us to obtain approximation of the resolvent $(\A_{N,\eps} - \zeta I)^{-1}$ in a regular point
$\zeta \in \C \setminus [c_\flat,\infty)$, $\zeta \ne 0$.

\smallskip\noindent\textbf{Theorem 14.8.} \textit{Let $\zeta \in \C \setminus [c_\flat,\infty)$
and $\zeta \ne 0$. Let $\u_\eps$ be the solution of problem} (9.6), \textit{and let $\u_0$ be the solution of problem} (9.13)
\textit{with $\FF \in L_2(\O;\C^n)$.
Suppose that the number $\eps_1$ satisfies Condition} 4.1.
\textit{Then for $0< \eps \leq \eps_1$ we have}
$$
\| \u_\eps - \u_0 \|_{L_2(\O)}
\leq {\mathfrak C}_{30} \rho_\flat(\zeta) \eps \| \FF \|_{L_2(\O)},
$$
\textit{where $\rho_\flat(\zeta)$ is defined by} (14.7). \textit{In operator terms,
$$
\| (\A_{N,\eps} - \zeta I)^{-1} - (\A_{N}^0  - \zeta I)^{-1} \|_{L_2(\O) \to L_2(\O)}
\leq {\mathfrak C}_{30} \rho_\flat(\zeta) \eps.
\eqno(14.43)
$$
Denote $\widehat{\v}_\eps = \u_0 + \eps \Lambda^\eps S_\eps b(\D) \widehat{\u}_0$, where
$\widehat{\u}_0 = P_\O (\A_N^0 - \zeta I)^{-1} {\mathcal P} \FF$.
Then for $0< \eps \leq \eps_1$ we have}
$$
\| \u_\eps - \widehat{\v}_\eps \|_{H^1(\O)} \leq {\mathfrak C}_{31}
\rho_\flat(\zeta) \eps^{1/2} \|\FF\|_{L_2(\O)},
$$
\textit{or, in operator terms},
$$
\| (\A_{N,\eps} - \zeta I)^{-1} - (\A_{N}^0  - \zeta I)^{-1} -
\eps \widehat{K}_N(\eps;\zeta) \|_{L_2(\O) \to H^1(\O)}
\leq {\mathfrak C}_{31} \rho_\flat(\zeta) \eps^{1/2}.
\eqno(14.44)
$$
\textit{Here} $\widehat{K}_N(\eps;\zeta):= R_\O [\Lambda^\eps] S_\eps b(\D) P_\O (\A_N^0 - \zeta I)^{-1} {\mathcal P}$.
\textit{For the flux $\p_\eps = g^\eps b(\D) \u_\eps$ we have}
$$
\| \p_\eps - \widetilde{g}^\eps S_\eps b(\D) \widehat{\u}_0 \|_{L_2(\O)}
\leq {\mathfrak C}_{32} \rho_\flat(\zeta) \eps^{1/2} \|\FF\|_{L_2(\O)}
\eqno(14.45)
$$
\textit{for $0< \eps \leq \eps_1$. The constants ${\mathfrak C}_{30}$, ${\mathfrak C}_{31}$, and ${\mathfrak C}_{32}$ are the same as in
Theorem} 14.2.

\smallskip\noindent\textbf{Proof.} Note that for $\zeta \in \C \setminus [c_\flat, \infty)$, $\zeta \ne 0$, we have
$$
(\A_{N,\eps} -\zeta I)^{-1} = (\A_{N,\eps} -\zeta I)^{-1}{\mathcal P}
+ (\A_{N,\eps} -\zeta I)^{-1}{\mathcal P}_Z,
$$
and $(\A_{N,\eps} -\zeta I)^{-1}{\mathcal P}= (\B_{N,\eps} -\zeta I)^{-1}{\mathcal P}$,
$(\A_{N,\eps} -\zeta I)^{-1}{\mathcal P}_Z = - \zeta^{-1} {\mathcal P}_Z$.
Hence,
$$
(\A_{N,\eps} -\zeta I)^{-1} = (\B_{N,\eps} -\zeta I)^{-1}{\mathcal P}
- \zeta^{-1} {\mathcal P}_Z.
$$
Similarly,
$$
(\A_{N}^0 -\zeta I)^{-1} = (\B_{N}^0 -\zeta I)^{-1}{\mathcal P}
- \zeta^{-1} {\mathcal P}_Z.
$$
Therefore,
$$
(\A_{N,\eps} -\zeta I)^{-1} - (\A_{N}^0 -\zeta I)^{-1}
= \left((\B_{N,\eps} -\zeta I)^{-1} - (\B_{N}^0 -\zeta I)^{-1}\right){\mathcal P}.
\eqno(14.46)
$$

Estimate (14.43) follows directly from (14.9) and (14.46).
Note that $\widehat{K}_N(\eps;\zeta) = {\mathcal K}_N(\eps;\zeta) {\mathcal P}$.
Then
$$
\begin{aligned}
&(\A_{N,\eps} - \zeta I)^{-1} - (\A_{N}^0  - \zeta I)^{-1} -
\eps \widehat{K}_N(\eps;\zeta)
\cr
&= \left( (\B_{N,\eps} - \zeta I)^{-1} - (\B_{N}^0  - \zeta I)^{-1} -
\eps {\mathcal K}_N(\eps;\zeta)\right) {\mathcal P}.
\end{aligned}
\eqno(14.47)
$$
Together with (14.10) this implies (14.44).

Next, since $b(\D) (\A_{N,\eps} -\zeta I)^{-1}{\mathcal P}_Z =0$, then
$$
\begin{aligned}
&g^\eps b(\D) (\A_{N,\eps} -\zeta I)^{-1} - \widetilde{g}^\eps S_\eps b(\D) P_\O
(\A_N^0 - \zeta I)^{-1} {\mathcal P}
\cr
&=
\left( g^\eps b(\D) (\B_{N,\eps} -\zeta I)^{-1} - \widetilde{g}^\eps S_\eps b(\D) P_\O
(\B_N^0 - \zeta I)^{-1} \right) {\mathcal P}.
\end{aligned}
\eqno(14.48)
$$
Obviously, in operator terms (14.11) means that
$$
\|g^\eps b(\D) (\B_{N,\eps} -\zeta I)^{-1} - \widetilde{g}^\eps S_\eps b(\D) P_\O
(\B_N^0 - \zeta I)^{-1} \|_{{\mathcal H}(\O) \to L_2(\O)} \leq {\mathfrak C}_{32} \rho_\flat(\zeta) \eps^{1/2}
\eqno(14.49)
$$
for $0< \eps \leq \eps_1$.
Now relations (14.48) and (14.49) imply that
$$
\begin{aligned}
&\| g^\eps b(\D) (\A_{N,\eps} -\zeta I)^{-1} - \widetilde{g}^\eps S_\eps b(\D) P_\O
(\A_N^0 - \zeta I)^{-1} {\mathcal P}\|_{L_2(\O) \to L_2(\O)}
\cr
&\leq {\mathfrak C}_{32} \rho_\flat(\zeta) \eps^{1/2},\ \ 0< \eps \leq \eps_1,
\end{aligned}
$$
which is equivalent to  (14.45). $\bullet$

\smallskip
Theorem 14.3 implies the following result.

\smallskip\noindent\textbf{Theorem 14.9.} \textit{Let $\zeta \in \C \setminus [c_\flat,\infty)$
and $\zeta \ne 0$. Let $\u_\eps$ be the solution of problem} (9.6), \textit{and let $\u_0$ be the solution of problem} (9.13)
\textit{with $\FF \in L_2(\O;\C^n)$. Suppose that Condition} 2.8
\textit{is satisfied.  Suppose that $K_N^0(\eps;\zeta)$ is defined by} (12.1)
\textit{and $\check{\v}_\eps$ is given by} (12.2).
 \textit{Suppose that the number $\eps_1$ satisfies Condition} 4.1.
\textit{Then for $0< \eps \leq \eps_1$ we have}
$$
\| \u_\eps - \check{\v}_\eps \|_{H^1(\O)}
\leq {\mathfrak C}_{31}^\circ \rho_\flat(\zeta) \eps^{1/2} \| \FF \|_{L_2(\O)},
$$
\textit{or, in operator terms,
$$
\| (\A_{N,\eps} - \zeta I)^{-1} - (\A_{N}^0  - \zeta I)^{-1} - \eps K_N^0(\eps;\zeta)
\|_{L_2(\O) \to H^1(\O)}
\leq {\mathfrak C}_{31}^\circ \rho_\flat(\zeta) \eps^{1/2}.
\eqno(14.50)
$$
For the flux $\p_\eps = g^\eps b(\D) \u_\eps$ we have}
$$
\| \p_\eps - \widetilde{g}^\eps b(\D) {\u}_0 \|_{L_2(\O)}
\leq {\mathfrak C}_{32}^\circ \rho_\flat(\zeta) \eps^{1/2} \|\FF\|_{L_2(\O)}
\eqno(14.51)
$$
\textit{for $0< \eps \leq \eps_1$. The constants ${\mathfrak C}_{31}^\circ$ and ${\mathfrak C}_{32}^\circ$ are the same
as in Theorem} 14.3.

\smallskip\noindent\textbf{Proof.} By (12.1), (14.31), and the identity $b(\D) {\mathcal P}_Z =0$, we have
$K_N^0(\eps;\zeta) = {\mathcal K}_N^0(\eps;\zeta) {\mathcal P}$. Together with (14.46) this yields
$$
\begin{aligned}
&(\A_{N,\eps} - \zeta I)^{-1} - (\A_{N}^0  - \zeta I)^{-1} - \eps K_N^0(\eps;\zeta)
\cr
&=
\left((\B_{N,\eps} - \zeta I)^{-1} - (\B_{N}^0  - \zeta I)^{-1} - \eps {\mathcal K}_N^0(\eps;\zeta) \right)
{\mathcal P}.
\end{aligned}
\eqno(14.52)
$$
Relations (14.34) and (14.52) imply (14.50).

Next, since $b(\D) {\mathcal P}_Z=0$, then
$$
\begin{aligned}
&g^\eps b(\D) (\A_{N,\eps} - \zeta I)^{-1} - \widetilde{g}^\eps b(\D) (\A_{N}^0 - \zeta I)^{-1}
\cr
&=
\left(g^\eps b(\D) (\B_{N,\eps} - \zeta I)^{-1} - \widetilde{g}^\eps b(\D) (\B_{N}^0 - \zeta I)^{-1}\right)
{\mathcal P}.
\end{aligned}
\eqno(14.53)
$$
Relations (14.35) and (14.53) yield (14.51). $\bullet$

\smallskip
The next result follows from Theorem 14.6.

\smallskip\noindent\textbf{Theorem 14.10.} \textit{Suppose that the assumptions of Theorem} 14.8
\textit{are satisfied. Let $\O'$ be a strictly interior subdomain of the domain $\O$, and let $\delta :=
{\rm dist}\,\{ \O'; \partial \O \}$. Then for any $0< \eps \leq \eps_1$ we have}
$$
\begin{aligned}
&\|\u_\eps - \widehat{\v}_\eps \|_{H^1(\O')}
\cr
&\leq
\left( {\mathfrak C}_{\flat}' \delta^{-1} (c(\vartheta) \rho_\flat(\zeta) +
c(\vartheta)^{5/2} \rho_\flat(\zeta)^{3/4}) + {\mathfrak C}_{\flat}''
c(\vartheta)^{1/2} \rho_\flat(\zeta)^{5/4} \right)
\eps \|\FF\|_{L_2(\O)},
\end{aligned}
$$
\textit{or, in operator terms},
$$
\begin{aligned}
&\| (\A_{N,\eps} - \zeta I)^{-1} - (\A_{N}^0  - \zeta I)^{-1} - \eps \widehat{K}_N(\eps;\zeta)
\|_{{\mathcal H}(\O) \to H^1(\O')}
\cr
&\leq
\left( {\mathfrak C}_{\flat}' \delta^{-1} (c(\vartheta) \rho_\flat(\zeta) +
c(\vartheta)^{5/2} \rho_\flat(\zeta)^{3/4}) + {\mathfrak C}_{\flat}''
c(\vartheta)^{1/2} \rho_\flat(\zeta)^{5/4} \right) \eps.
\end{aligned}
\eqno(14.54)
$$
\textit{For the flux $\p_\eps = g^\eps b(\D) \u_\eps$ we have}
$$
\begin{aligned}
&\| \p_\eps  -
\widetilde{g}^\eps S_\eps b(\D) \widehat{\u}_0 \|_{L_2(\O')}
\cr
&\leq \left( \widetilde{\mathfrak C}_{\flat}' \delta^{-1} (c(\vartheta) \rho_\flat(\zeta) +
c(\vartheta)^{5/2} \rho_\flat(\zeta)^{3/4}) + \widetilde{\mathfrak C}_{\flat}''
c(\vartheta)^{1/2} \rho_\flat(\zeta)^{5/4} \right)
\eps \|\FF\|_{L_2(\O)}
\end{aligned}
\eqno(14.55)
$$
\textit{for $0< \eps \leq \eps_1$. The constants ${\mathfrak C}_{\flat}'$, ${\mathfrak C}_{\flat}''$,
$\widetilde{\mathfrak C}_{\flat}'$, and $\widetilde{\mathfrak C}_{\flat}''$ are the same as in Theorem} 14.6.

\smallskip\noindent\textbf{Proof.} Estimate (14.54) is a consequence of (14.39) and (14.47).
Inequality (14.55) follows from (14.40) and (14.48). $\bullet$

\smallskip
The next result is deduced from Theorem 14.7.

\smallskip\noindent\textbf{Theorem 14.11.} \textit{Suppose that the assumptions of Theorem} 14.9
\textit{are satisfied. Let $\O'$ be a strictly interior subdomain of the domain $\O$, and let $\delta :=
{\rm dist}\,\{ \O'; \partial \O \}$. Then for $0< \eps \leq \eps_1$ we have}
$$
\begin{aligned}
&\|\u_\eps - \check{\v}_\eps \|_{H^1(\O')}
\cr
&\leq
\left( {\mathfrak C}_{\flat}' \delta^{-1} (c(\vartheta) \rho_\flat(\zeta) +
c(\vartheta)^{5/2} \rho_\flat(\zeta)^{3/4}) + \check{\mathfrak C}_{\flat}''
c(\vartheta)^{1/2} \rho_\flat(\zeta)^{5/4} \right)
\eps \|\FF\|_{L_2(\O)},
\end{aligned}
$$
\textit{or, in operator terms},
$$
\begin{aligned}
&\| (\A_{N,\eps} - \zeta I)^{-1} - (\A_{N}^0  - \zeta I)^{-1} - \eps {K}^0_N(\eps;\zeta)
\|_{L_2(\O) \to H^1(\O')}
\cr
&\leq
\left( {\mathfrak C}_{\flat}' \delta^{-1} (c(\vartheta) \rho_\flat(\zeta) +
c(\vartheta)^{5/2} \rho_\flat(\zeta)^{3/4}) + \check{\mathfrak C}_{\flat}''
c(\vartheta)^{1/2} \rho_\flat(\zeta)^{5/4} \right) \eps.
\end{aligned}
\eqno(14.56)
$$
\textit{For the flux $\p_\eps = g^\eps b(\D) \u_\eps$ we have}
$$
\begin{aligned}
&\| \p_\eps  -
\widetilde{g}^\eps b(\D) {\u}_0 \|_{L_2(\O')}
\cr
&\leq \left( \widetilde{\mathfrak C}_{\flat}' \delta^{-1} (c(\vartheta) \rho_\flat(\zeta) +
c(\vartheta)^{5/2} \rho_\flat(\zeta)^{3/4}) + \widehat{\mathfrak C}_{\flat}''
c(\vartheta)^{1/2} \rho_\flat(\zeta)^{5/4} \right)
\eps \|\FF\|_{L_2(\O)}
\end{aligned}
\eqno(14.57)
$$
\textit{for $0< \eps \leq \eps_1$. The constants ${\mathfrak C}_{\flat}'$, $\check{\mathfrak C}_{\flat}''$,
$\widetilde{\mathfrak C}_{\flat}'$, and $\widehat{\mathfrak C}_{\flat}''$ are the same as in Theorem} 14.7.

\smallskip\noindent\textbf{Proof.} Estimate (14.56) follows from (14.41) and (14.52).
Inequality (14.57) is a consequence of (14.42) and (14.53). $\bullet$

\end{document}